\documentclass{amsart}
\usepackage{amssymb}
\usepackage{amsfonts}
\usepackage{amsmath}
\usepackage{amsthm}
\usepackage{enumerate}
\usepackage[colorlinks,citecolor=blue,urlcolor=blue, linkcolor=blue]{hyperref}

\usepackage[nocompress]{cite}

\usepackage{tabularx}
\usepackage{centernot}
\usepackage{mathtools}
\usepackage{stmaryrd}
\usepackage{amsthm,amssymb}
\usepackage{etoolbox}
\usepackage{url}
\usepackage{tikz}
\usepackage{amssymb}
\usetikzlibrary{matrix}
\usepackage{tikz-cd}
\usepackage{tikz}
\usepackage{marginnote}
\definecolor{mygray}{gray}{0.85}
\usepackage[backgroundcolor=mygray,colorinlistoftodos,prependcaption,textsize=small]{todonotes}
\usepackage{xargs}
\usepackage{subfigure}
\usepackage[capitalise]{cleveref}

\allowdisplaybreaks[3]

\newcommand{\mrm}[1]{\mathrm{#1}}

\renewcommand{\leq}{\leqslant}
\renewcommand{\geq}{\geqslant}

\newcommand{\VdashR}{\mathrel{\reflectbox{$\Vdash$}}}

\makeatletter
\def\subsection{\@startsection{subsection}{3}%
	\z@{.5\linespacing\@plus.7\linespacing}{.3\linespacing}%
	{\bfseries\centering}}
\makeatother

\makeatletter
\def\subsubsection{\@startsection{subsubsection}{3}%
	\z@{.5\linespacing\@plus.7\linespacing}{.3\linespacing}%
	{\centering}}
\makeatother

\makeatletter
\def\myfnt{\ifx\protect\@typeset@protect\expandafter\footnote\else\expandafter\@gobble\fi}
\makeatother

\renewcommand{\restriction}{ {\upharpoonright} }
\newtheorem{theorem}{Theorem}[section]
\newtheorem*{theorem*}{Main Theorem}
\theoremstyle{definition}
\newtheorem{corollary}[theorem]{Corollary}
\newtheorem{definition}[theorem]{Definition}
\newtheorem{lemma}[theorem]{Lemma}
\newtheorem{proposition}[theorem]{Proposition}
\newtheorem{example}[theorem]{Example}

\newtheorem{question}[theorem]{Question}

\newtheorem{fact}[theorem]{Fact}

\newtheorem{remark}[theorem]{Remark}
\newtheorem{notation}[theorem]{Notation}
\newtheorem{context}[theorem]{Context}

\usepackage{xcolor}

\newcounter{claimcounter}
\numberwithin{claimcounter}{theorem}

\newcommand\dna{\mathtt{DNA}}

\newcommand\cpc{\mathtt{CPC}}
\newcommand\ipc{\mathtt{IPC}}

\newcommand\Fm{\mathsf{Fm}}
\newcommand\Eq{\mathsf{Eq}}
\newcommand\QEq{\mathsf{QEq}}

\newcommand{\atoms}{\mathsf{Var}}
\newcommand{\formulas}{\mathsf{Fm}}
\newcommand{\subst}{\mathsf{Subst}}
\newcommand{\atomsubst}{\mathsf{At}}
\newcommand{\admsubst}{\mathsf{AS}}
\newcommand{\id}{\mathsf{id}}

\newcommand\inqI{\mathtt{InqI}}
\newcommand\inqB{\mathtt{InqB}}

\newcommand\PAL{\mathtt{PAL}}

\newcommand\vari{\mathbf{V}}

\newcommand{\core}{\mathsf{core}}
\newcommand{\dom}{\mathsf{dom}}

\newcommand{\truth}{\mathsf{truth}}

\newcommand{\logic}{\mathsf{Log}}

\newcommand{\homset}{\mathsf{Hom}}

\newcommand{\strongend}{\mathsf{sEnd}}
\newcommand{\enset}{\mathsf{End}}
\newcommand{\horn}{\mathsf{Horn}}

\newcommand{\cA}{A}
\newcommand{\cB}{B}

\newcommand{\cD}{D}
\newcommand{\cH}{H}

\newcommand{\cL}{\mathcal{L}}

\newcommand{\cM}{M}

\newcommand{\bK}{\mathbf{K}}
\newcommand{\bQ}{\mathbf{Q}}
\newcommand{\bV}{\mathbf{V}}
\newcommand{\bC}{\mathbf{C}}

\newcommand{\BA}{\mathbf{BA}}
\newcommand{\HA}{\mathbf{HA}}

\newcommand{\varopCore}{\mathbb{C}_\Sigma}

\newcommand{\coregen}[1]{#1_{\mrm{CG}}}
\newcommand{\coremodels}[1]{\mathrm{Mod}^c(#1)}
\newcommand{\coregenmodels}[1]{\mathrm{Mod}_{\mrm{CG}}^c(#1)}
\newcommand{\coretheories}[1]{\mathrm{Th}_{\mrm{qe}}^c(#1)}
\newcommand{\coreequations}[1]{\mathrm{Th}_{\mrm{e}}^c(#1)}

\newcommand{\foltheory}[1]{\mathrm{Th}_{\mrm{fol}}(#1)}
\newcommand{\theory}[1]{\mathrm{Th}_{\mrm{qe}}(#1)}
\newcommand{\equations}[1]{\mathrm{Th}_{\mrm{e}}(#1)}
\newcommand{\horntheory}[1]{\mathrm{Th}_{\mrm{h}}(#1)}

\newcommand{\corevDash}{\models^c}
\newcommand{\ncorevDash}{\not\models^c}

\newcommand{\Fi}{\mathsf{Fi}}
\newcommand{\Con}{\mathsf{Con}}
\newcommand{\Th}{\mathsf{Th}}

\newcommand{\xbar}{\bar{x}}

\newcommand\intlog{\mathtt{L}}
\newcommand\ml{\mathtt{ML}}

\newcommand\langInqI{\mathcal{L}_{\ipc}^\otimes}
\newcommand\langInt{\mathcal{L}_{\ipc}}
\newcommand\langIntsta{\mathcal{L}_{\mathtt{CL}}}

\newcommand{\Schm}{\mathrm{Schm}}
\newcommand{\Cn}{\mathrm{Cn}}
\newcommand{\Mod}{\mathrm{Mod}}
\renewcommand{\restriction}{{\upharpoonright}}

\usepackage{tikz}
\usetikzlibrary{cd}
\usetikzlibrary{calc}
\tikzstyle{dot}=[inner sep=1pt, fill, black, circle, draw, minimum size = 5pt]
\tikzstyle{highlight}=[inner sep=1pt, draw=black, fill=white, circle, minimum size = 9pt]

\tikzstyle{index on}=[inner sep=2pt, white, circle, fill=black]
\tikzstyle{index off}=[inner sep=2pt, black, circle, draw]
\tikzstyle{index gray}=[inner sep=2pt, black, circle, fill=lightgray]
\tikzstyle{opaque}=[fill=gray,fill opacity=.1]
\tikzstyle{counter}=[densely dashed]

\title{Algebraizable Weak Logics}

\author{Georgi Nakov}
\author{Davide Emilio Quadrellaro}

\address{Department of Computer and Information Sciences,
	University of Strathclyde,
	Glasgow, United Kingdom.}
\email{georgi.nakov@strath.ac.uk}

\address{Department of Mathematics “Giuseppe Peano”,
	University of Torino,
	Via Carlo Alberto 10,
	10123 Torino, Italy.}
\email{davideemilio.quadrellaro@unito.it}

\address{Istituto Nazionale di Alta Matematica ``Francesco Severi'',
	Piazzale Aldo Moro 5
	00185 Roma, Italy.}
\email{quadrellaro@altamatematica.it}

\begin{document}

\begin{abstract}
  We extend the framework of abstract algebraic logic to weak logics,
  namely logical systems that are not necessarily closed under uniform
  substitution. We interpret weak logics by algebras expanded with an
  additional predicate, and we introduce a loose and strict version of
  algebraizability for weak logics. We study this framework by
  investigating the connection between the algebraizability of a weak
  logic and the algebraizability of its schematic fragment, and we
  then prove a version of Blok and Pigozzi's Isomorphism Theorem in
  our setting. We apply this framework to logics in team semantics and
  show that the classical versions of inquisitive and dependence logic
  are strictly algebraizable, while their intuitionistic versions are
  only loosely so.
\end{abstract}

\subjclass[2020]{03G27, 03C05, 03B55, 03B50}

\date{\today}
\maketitle

%
\section*{Introduction}

At least since Tarski's 1936 article on logical notions, logic has
often been understood as the subject studying those notions which are
``invariant under all possible one-one transformations of the world
onto itself'' \cite[p. 149]{Tarski1986-TARWAL}. Tarski's view was
inspired by Felix Klein's 1872 \emph{Erlanger Programm}, and expresses
very neatly the philosophical idea of logic as the discipline with the
most general character. While, for instance, metric geometry studies
notions that are invariant under transformations that preserve
distances, and topology studies notions invariant under continuous
maps, logic deals with notions which are invariant under arbitrary
bijections. Formally, this led to the technical definition of logics
as \emph{consequence relations which are additionally closed under
  uniform substitutions}.

As a matter of fact, Tarski's approach can be seen as a key step in
the transition from the ``symbolic'' algebraic logic of the
19\emph{th} century (exemplified by the works of Boole, De Morgan,
Jevons, Peirce, etc.) to the contemporary field of abstract algebraic
logic. Using the general notion of logic as a closure operator on the
term algebra, Rasiowa made the first steps into the ``mathematical''
version of algebraic logic, and in particular she developed in
\emph{An Algebraic Approach to Non-Classical Logics}
\cite{Rasiowa1974-RASAAA-2} a general theory of algebraization for
implicative logics. Finally, the algebraic approach was put in its
contemporary formulation by Blok and Pigozzi, who introduced in
\emph{Algebraizable Logics} \cite{Blok1989-BLOAL} the notion of
algebraizable logics and developed their general theory.

In recent years, however, there has been an increasing interest into
systems which do possess a logical nature but fail nonetheless to be
closed under uniform substitution in Tarski's original strong
sense. The field of modal logic is particularly rich of such examples:
Buss' pure provability logic \cite{buss1990modal}, public announcement
logic \cite{Holliday2013-HOLIDA-5, 10.1007/978-3-642-24130-7_6} and
other epistemic logics are all examples of this
behaviour. Furthermore, propositional logics based on team semantics,
such as inquisitive \cite{Ciardelli2011-CIAIL} and dependence logic
\cite{Yang2016-YANPLO}, also do not satisfy Tarski's requirement of
closure under uniform substitution. We believe that this state of
affairs is not a mistake requiring correction, but that it rather
reflects the increasing plurality of logics and their aptness for
applications. At the same time, however, the anomalous behaviour of
these logical systems prevented so far a uniform abstract study, and
did not allow for immediately applying facts and results from abstract
algebraic logic, thus forcing scholars to reprove abstract results in
these settings, or adapt standard techniques to their specific
situation.

Motivated by these facts, we propose in this article a generalisation
of the notion of algebraizable logic which breaks apart from Tarski's
original view of logics as being invariant under \emph{all}
substitutions. In other words, we want to study consequence operators
that are invariant under \emph{some} substitutions, but not
necessarily all. We introduce the notion of \emph{weak logic},
generalising previous definitions of Ciardelli and Roelofsen
\cite{Ciardelli2011-CIAIL} and Punčochář \cite{Puncochar2016-PUNAGO},
as a consequence relation which is invariant under all substitutions
which map atomic variables to atomic variables. In other words, we
relax Tarski's constraint of invariance under arbitrary substitutions
and we require that it holds only with respect to these so-called
atomic substitutions (cf. \cref{sec.1} below).

We proceed in this article as follows. In Sections
\ref{sec.1}-\ref{sec.3.2} we develop at length the key aspects of the
theory of algebraizable weak logics. In \cref{sec.1} we define weak
logics and introduce what we call expanded algebras as their
corresponding algebraic notion. In \cref{sec.2} we introduce the
notions of loose and strict algebraizability for weak logics and we
show that the (loose or strict) equivalent algebraic semantics of a
weak logic is unique, thus mirroring the classical result by Blok and
Pigozzi for standard algebraizable logics. Then, in \cref{sec.3.1}, we
study the relation between a weak logic and its schematic fragment of
consequences invariant under arbitrary substitutions. In particular,
we provide a characterisation of the algebraizability of weak logics
in terms of the algebraizability of their schematic
fragment. Additionally, we introduce the notion of \emph{standard
  companion} and generalise to this setting previous results from
\cite{Quadrellaro.2019B} and \cite{nicolau2024polyatomic}. Finally, in
\cref{sec.3.2} we develop on these results and we prove a version of
Blok and Pigozzi's Isomorphism Theorem for strictly algebraizable weak
logics. These sections thus establish the fundamentals of the theory
of algebraizable weak logics and show that several key results from
the field of abstract algebraic logic carry on to the setting without
uniform substitution. On a partially separate line of inquiry, in
\cref{sec.5} we take a small excursus explaining how to adapt the
usual matrix semantics of propositional logics to the setting of weak
logics. We show that (similarly to the standard setting) every weak
logic admits a matrix semantics, and we describe the class of reduced
models of (strictly) algebraizable weak logics. We put our abstract
framework to the test in \cref{sec.7}, where we apply it to the
specific case of inquisitive and dependence logics. In particular, we
build on previous results from
\cite{ciardelli2020,quadrellaro.2021,Quadrellaro.2019B} to show that
the classical version of inquisitive and dependence logic is strictly
algebraizable, while their intuitionistic versions are only loosely
so. To our eyes, this indicates a significant difference between the
classical version of inquisitive and dependence logic --- which can
essentially be recasted as a theory over their schematic fragment ---
and their intuitionistic counterparts, which do not admit such
reformulation.

Together with these results on inquisitive and dependence logic, we
regard as the main contribution of the present work the fact that it
provides a framework for reasoning about logical systems lacking
uniform substitution. In particular, this work relates several
algebraic studies on inquisitive logic
\cite{grilletti,Quadrellaro.2019B,Puncochar2021-PUNIHA,Puncochar2016-PUNAGO},
dependence logic \cite{quadrellaro.2021}, and polyatomic logics
\cite{nicolau2024polyatomic}, by showcasing all the algebraic
semantics from these works as instances of what we call core semantics
in the present article. Given the increasing popularity of logics
without uniform substitution in the logic literature, we hope that our
approach will be useful to researchers working in these areas also in
the future.

\section{Preliminaries}

We recall in this section some basic facts concerning logics, algebras
and model theory. We also fix some notation that we shall follow
throughout the rest of the article. The following general context sets
the framework of our work.

\begin{context}
  Throughout this article we always let $\atoms$ be a set of variables
  and we let $\cL$ be an algebraic (i.e., purely functional) signature
  (unless we specify otherwise). We denote by $\formulas_\cL$ both the
  set of first-order terms over $\atoms$ in the signature $\cL$ and
  the term algebra in the signature $\cL$ over $\atoms$. We omit the
  index $\cL$ when it is clear from the context. Notice that, since we
  are often dealing with propositional logical systems, we often refer
  to elements of $\formulas_\cL$ also as \emph{(propositional)
    formulas} in the language $\cL$. These should not be confused with
  the first-order formulas in the signature $\cL$.
\end{context}

Given the language $\cL$, we recall the standard abstract Tarskian
definition of (propositional) logic. We refer the reader to
\cite[\S1]{Font.2016} for more on consequence relations and logics,
and for slight variations of these definitions.

\begin{definition}\label{def:consequence_relation}
	A \emph{(finitary) consequence relation} is a relation $\vdash \; \subseteq \: \wp(\formulas)\times \formulas$ such that, for all $\Gamma \subseteq \formulas$:
	\begin{enumerate}[(1)]
		\item $\Gamma \vdash \phi$ for all $\phi \in \Gamma$;
		\item if $\Gamma \vdash \phi$ for all $\phi \in \Delta $, and $\Delta \vdash \psi $, then $\Gamma \vdash \psi$;
		\item if $\Gamma \vdash \phi$ and $\Gamma\subseteq \Delta$, then $\Delta\vdash \phi$;
		\item if $\Gamma \vdash \phi$ then there is some $\Delta \subseteq \Gamma$ such that $|\Delta| < \aleph_0$ and   $\Delta\vdash\phi$.
	\end{enumerate}
\end{definition}

\begin{remark}\label{remark:finitary}
  Condition~\ref{def:consequence_relation}(3) already follows from
  Condition~\ref{def:consequence_relation}(1) together with
  Condition~\ref{def:consequence_relation}(2) (cf. also
  \cite[p.~14]{Font.2016}), but we include it to make explicit the
  fact that we always assume consequence relations (and logics) to be
  monotone. Differently, Condition~\ref{def:consequence_relation}(4)
  is not always required in the definition of consequence relation and
  it essentially restricts the focus to the finitary ones.  This
  reflects our elementary approach to the subject, as
  Condition~\ref{def:consequence_relation}(4) allows us to translate
  propositional systems into first-order logic and to avoid the use of
  infinitary logical systems lacking compactness. We stress that this
  requirement is not necessary, and it is possible to study
  consequence relations and propositional logics where this condition
  fails. We refer the interested reader to \cite{Font.2016} for
  discussions of Condition~\ref{def:consequence_relation}(4) and also
  for an (algebraic) treatment of non-finitary propositional logics.
\end{remark}

\begin{definition}
  A \emph{substitution} is an endomorphism
  $\sigma : \formulas \to \formulas$ of the $\cL$-term algebra. We
  denote by $\subst(\cL)$ the set of all substitutions in the language
  $\cL$. If $x_1,\dots, x_n\in \atoms$ and $\phi_1,\dots,\phi_n$ are
  arbitrary $\cL$-formulas, we denote by
  $\Gamma[\phi_1 \dots \phi_n/x_1 \dots x_n]$ the result of
  simultaneously substituting each $\phi_i$ for all occurrences of
  $x_i$ in the formulas in $\Gamma$.
\end{definition}

\begin{definition}[Logic]
  A consequence relation $\vdash$ is closed under \emph{uniform
    substitution} if $\Gamma \vdash \phi$ entails
  $ \sigma[\Gamma] \vdash \sigma(\phi)$ for all substitutions
  $\sigma\in \subst(\cL)$.  A \emph{(standard) logic} is a consequence
  relation $\vdash$ which is closed under uniform substitution.
\end{definition}

\begin{example}
  Obvious examples of standard logics are the classical propositional
  logic $\cpc$ and the intuitionistic propositional logic
  $\ipc$. Non-examples of logics in this sense are first-order logic,
  as it is not a consequence relation of the propositional term
  algebra, infinitary logics, and higher-order systems.
\end{example}

In the context of abstract algebraic logic, one is often interested in
the \emph{algebraic semantics} of a propositional logic $\vdash$. This
is provided by algebras, i.e., first-order structures in some purely
functional language $\cL$. We first fix some notation concerning
first-order models.

\begin{notation}
  Let $\cL$ be a first-order language, not necessarily functional. We
  use Latin letters $\cA, \cB, \dots$ both to denote first-order
  $\cL$-structures and their underlying domain. When confusion may
  arise, we also write $\dom(\cA)$ to refer to the underlying universe
  of $\cA$. For all function symbols $f \in \cL$ and all relation
  symbols $R \in \cL$, we write $f^{\cA}$ and $R^{\cA}$ for their
  interpretation in $\cA$. We use the same notations for symbols and
  their interpretation when it does not cause confusion. If
  $t(\bar{x})$ is a term in the language $\cL$ (i.e., a propositional
  formula), we usually call its interpretation $t^{\cA}$ a
  polynomial. If $X\subseteq \cA$, we write $\langle X\rangle $ for
  the substructure of $\cA$ generated by $X$. The symbol $\models$
  refers to the standard satisfaction symbol from first-order
  logic. We usually denote classes of structures by boldface font ---
  both for arbitrary collections $(\bQ, \mathbf{K}, \dots)$ and
  designated ones, e.g., the class of all Heyting algebras
  $\mathbf{HA}$ or the class of all Boolean algebras $ \mathbf{BA}$.
\end{notation}

\noindent We assume the reader is familiar with the usual
constructions from model theory and universal algebra, and refer to
\cite{chang1990model,Burris.1981} for background. We recall in
particular the following notions of maps, as we shall need them in the
rest of the article.

\begin{definition}\label{maps:expanded_algebras} Let $h:\cA \to \cB$
  be a function between two $\cL$-structures $\cA$ and $\cB$, for some
  first-order language $\cL$. We define the following notions:
  \begin{enumerate}[(1)]
      \item we say that $h$ is a \emph{homomorphism} if for every function
        symbol $f\in \cL$ we have
	\[ h(f(a_1,\dots,a_n)) = f(h(a_1),\dots, h(a_n))  \]
	and for every relation symbol $R\in \cL$,
	\[ \cA\models R(a_1,\dots,a_n) \; \Longrightarrow\; \cB\models
          R(h(a_1),\dots, h(a_n));\]
      \item we say that $h$ is a \emph{strong homomorphism} if it is a
        homomorphism and, additionally, we have that $R^B= h[R^A]$ for
        every relation symbol $R\in \cL$;
      \item we say that $h$ is a \emph{strict homomorphism} if it is a
        homomorphism and, for every relation symbol $R\in \cL$,
	\[ \cA\models R(a_1,\dots,a_n) \; \Longleftrightarrow\; \cB\models R(h(a_1),\dots, h(a_n));\]
	\item we say that $h$ is an \emph{embedding} if it is an injective strict homomorphism.
\end{enumerate}
We write $\cA\leq \cB$ if $\cA$ is a substructure of $\cB$, i.e., if
the identity map $\id:\cA\to \cB$ is an embedding. We write
$\cA \cong \cB$ if $\cA$ is isomorphic to $\cB$. We say that $B$ is a
homomorphic image of $A$ if there is a surjective homomorphism
$h:A\twoheadrightarrow B$; we say that $B$ is a strong (resp.~strict)
homomorphic image of $A$ if $h$ is a strong (resp.~strict)
homomorphism.
\end{definition}

\begin{remark}\label{remark:def.homomorphism}
  We briefly explain the rationale behind the different versions of
  homomorphism. The notion of homomorphism from
  Definition~\ref{maps:expanded_algebras} is standard from the
  literature in model theory and universal algebra
  (cf. \cite[pp. 70-71]{chang1990model},
  \cite[p. 203]{Burris.1981}). The notion of strong homomorphism comes
  from \cite[p. 321]{chang1990model} and is motivated by the following
  observation. Let $A$ be an $\cL$-structure and let $\cL'$ consists
  of the functional symbols from $\cL$. If $\theta$ is a congruence of
  the algebraic reduct of $A$ we can consider the quotient
  $\cL'$-structure $A/\theta$ and expand it to a $\cL$-structure by
  letting $R^{A/\theta}=R^A/\theta$ for all relational symbols
  $R\in \cL$. Then, the projection map induced by this quotient is a
  strong homomorphism.  Finally, we take the notion of strict
  homomorphism from \cite{dellunde1996some,casanovas1996elementary},
  and we stress that strict homomorphisms correspond to those
  quotients which are additionally compatible with the relational part
  of the vocabulary from $\cL$, in the sense that, if
  $(a_i,b_i)\in \theta$ for all $1\leq i \leq n$, then for every
  relational symbol $R\in \cL$ we have that
  $A\models R(a_1,\dots,a_n)$ if and only if
  $A\models R(b_1,\dots,b_n)$. In this article we will mostly be
  dealing with strong homomorphisms, but we will consider strict
  homomorphisms in \cref{sec.5}.
\end{remark}

\begin{notation}\label{notation:class_operators}
  Let $\bK$ be any class of first-order structures. We denote by
  $\mathbb{I}(\bK)$ its closure under isomorphic copies, by
  $\mathbb{S}(\bK)$ its closure under substructures, by
  $\mathbb{P}(\bK)$ its closure under (direct) products and by
  $\mathbb{P}_{\mrm{U}}(\bK)$ its closure under
  ultraproducts. Finally, we write $\mathbb{H}(\bK)$ for its closure
  under strong homomorphic images, and $\mathbb{H}_{\mrm{s}}(\bK)$ for
  its closure under strict homomorphic images.

\end{notation}

\begin{definition}
  A class of algebras $\bK$ is a \emph{quasivariety} if it is closed
  under the operators $\mathbb{I}$, $\mathbb{S}$, $\mathbb{P}$,
  $\mathbb{P}_{\mrm{U}} $, i.e., if it is closed under isomorphic
  copies, subalgebras, products and ultraproducts. A class of algebras
  $\bK$ is a \emph{variety} if it is closed under $\mathbb{H}$,
  $\mathbb{S}$ and $\mathbb{P}$, i.e., if it is closed under
  homomorphic images, subalgebras and products. We denote by
  $\mathbb{Q}(\bK)$ and $\mathbb{V}(\bK)$ the quasivariety and the
  variety generated by $\bK$, respectively.
\end{definition}

Crucially, the closure of a class of structures under (some of) the
operators $\mathbb{I}$, $\mathbb{S}$, $\mathbb{P}$,
$\mathbb{P}_{\mrm{U}} $, $\mathbb{H}$ is related to conditions
pertaining its axiomatisability. Most famously, it is a fundamental
result by Keisler and Shelah that a class of structures $\bK$ is
elementary (i.e., first-order axiomatisable) if and only if $\bK$ is
closed under ultraproducts. In this work we are concerned with less
general definability conditions. We introduce them by focusing on
three special subclasses of first-order formulas.

\begin{notation}
  We write $\xbar$ as a shorthand for a sequence of variables
  $(x_0,\dots,x_n)$. Also, as we often deal with equational classes of
  structures, we abide to the usual convention from universal algebra
  to distinguish the syntactical equality symbol, written $\approx$,
  from the semantical equality symbol, which we denote by $=$.  To
  help the reader to distinguish when we talk of first-order formulas
  from when we deal with propositional ones, we use lowercase Greek
  symbols $\phi,\psi,\dots$ for the latter and uppercase Greek symbols
  $\Phi,\Psi,\dots$ for the former.
\end{notation}

\begin{definition}
  Let $\cL$ be an arbitrary first-order signature, not necessarily
  functional. Then we define the following types of formulas:
  \begin{enumerate}[(1)]
  \item an \emph{equation} is a formula of the form
    $\varepsilon \approx \delta$, where $\varepsilon$ and $\delta$ are
    two $\cL$-terms;
  \item a \emph{quasiequation} is a formula of the form
    $\bigwedge_{i\leq n} (\varepsilon_i\approx \delta_i) \to
    \varepsilon\approx \delta$, for some $n<\omega$, where all
    $\varepsilon_i,\delta_i$ and $\varepsilon,\delta$ are $\cL$-terms;
  \item a \emph{basic Horn formula} is a formula of the form
    $\bigvee_{1\leq i\leq n}\Psi_i$, where every $\Psi_i$ is either an
    atomic or a negated atomic formula, and at most one $\Psi_i$ is
    atomic;
  \item a basic Horn formula is \emph{strict} if exactly one of its
    disjuncts is atomic;
  \item A \emph{universal Horn formula} is a formula of the form
    \[ \forall x_1 \dots \forall x_m \; \big( \bigwedge_{1\leq i\leq
        \ell} \Phi_i(\xbar) \big) \] for some $m,\ell<\omega$, and
    where each $\Phi_i$ for $1\leq i\leq m$ is a basic Horn formula;
  \item a universal Horn formulas is \emph{strict} if all the basic
    Horn formulas $\Phi_i$ occurring in it are strict.
  \end{enumerate}
\end{definition}

\noindent The following classical results relate the closure under the
operators $\mathbb{I}$, $\mathbb{S}$, $\mathbb{P}$,
$\mathbb{P}_{\mrm{U}} $, $\mathbb{H}$ with different definability
conditions. We refer the reader to \cite[I: Thm. 2.23, Thm. 2.25; II:
Thm. 11.9]{Burris.1981} for their proofs. Facts (1) and (2) are also
known respectively as Birkhoff's Theorem and Maltsev's Theorem.

\begin{fact}\label{birkhoff_maltsev} Let $\cL$ be an algebraic
  signature an $\bK$ a class of $\cL$-algebras, then:
  \begin{enumerate}[(1)]
  \item $\bK$ is a variety if and only if it is axiomatised by
    equations;
  \item $\bK$ is a quasivariety if and only if it is axiomatised by
    quasiequations.
  \end{enumerate}
  Moreover, if $\cL$ is any first-order signature and $\bK$ is a class
  of $\cL$-structures, then:
\begin{enumerate}[(3)]
\item $\bK$ is closed under $\mathbb{I}$, $\mathbb{S}$, $\mathbb{P}$
  and $\mathbb{P}_{\mrm{U}} $ if and only if it is axiomatised by
  universal Horn sentences.
\end{enumerate}
\end{fact}

\noindent In the light of the previous results, it is convenient to
fix some shorthand notation to talk about (fragments of) first-order
theories, and their classes of models.

\begin{notation}\label{maltsev_birkhoff}
  Let $\cL$ be an arbitrary first-order language. If $T$ is a set of
  first-order sentences in $\cL$, we write $\Mod(T)$ for the class of
  all $\cL$-structures $\cA$ such that $\cA \models T$. We write
  $\Mod(\cL)$ for the class of all $\cL$-structures. Notice that, if
  $T$ is a set of formulas (e.g., a set of equations or
  quasiequations), then $\Mod(T)$ is the class of structures which
  models the sentences $\forall x_1\dots \forall x_n (\Phi(\xbar))$,
  for every $\Phi\in T$. On the other hand, let $\bK$ be a class of
  $\cL$-structures. Then we write $\foltheory{\bK}$ for the set of all
  $\cL$-sentences true in $\bK$; $\horntheory{\bK}$ for the set of all
  universal Horn sentences true in $\bK$; $\theory{\bK}$ for the set
  of all quasiequations true in $\bK$; and, finally, $\equations{\bK}$
  for the set of all equations true in $\bK$. We also write
  $\Eq_{\cL}$ for the set of all equations in $\cL$, and $\QEq_{\cL}$
  for the set of all quasiequations in $\cL$ --- we omit the index
  $\cL$ when it is clear from the context.
\end{notation}

Finally, we conclude this preliminary section by recalling what is the
propositional consequence relation induced by a class of algebras.

\begin{notation}
  Recall that $\Fm_\cL$ is the term algebra in the language $\cL$. A
  \emph{(propositional) assignment} is a homomorphism
  $h:\Fm_{\cL}\to A$, where $A$ is an $\cL$-structure. We denote the
  family of all assignments into $A$ as $\homset(\Fm,\cA)$. If $h:\Fm_{\cL}\to A$  is an assignment and $\alpha\approx \beta$ an equation, then we often write simply $h(\alpha\approx \beta)$ for the formula $h(\alpha)\approx h(\beta)$.
\end{notation}

\begin{definition}\label{remark:closure_quasivariety_validities}
  Let $\bK$ be a class of $\cL$-algebras and let
  $\Theta\cup\{\varepsilon\approx\delta \}\subseteq \Eq$ a set of
  equations, then the \emph{equational consequence relative to $\bK$}
  is defined as follows:
  \begin{align*}
    \Theta \models_{\bK} \varepsilon\approx \delta \Longleftrightarrow
    & \text{ for all } \cA \in \bK, \,  h\in \homset(\Fm,\cA), \\
    & \text{ if } h(\varepsilon_i)= h(\delta_i) \text{ for all } \varepsilon_i\approx \delta_i\in \Theta,  \text{ then } h(\varepsilon)= h(\delta).
\end{align*}
And we write
$\models_\bK \bigwedge_{i\leq n} \varepsilon_i\approx \delta_i \to
\varepsilon \approx \delta$ if
$ \bigwedge_{i\leq n} \varepsilon_i\approx \delta_i \models_{\bK}
\varepsilon\approx \delta $.  We often write
$\cA \models \bigwedge_{i\leq n} \varepsilon_i\approx \delta_i \to
\varepsilon \approx \delta$ in place of
$\models_{\{\cA \}} \bigwedge_{i\leq n} \varepsilon_i\approx \delta_i
\to \varepsilon \approx \delta$. The related notions for equations are
defined analogously.
\end{definition}

\begin{remark}\label{propositional_consequence_relation_first_order}
  Notice that, in
  Definition~\ref{remark:closure_quasivariety_validities}, we can
  assume without loss of generality that $\bK$ is a quasivariety,
  since by Fact~\ref{birkhoff_maltsev} the validity of quasiequations
  is preserved under the operations $\mathbb{I}$, $\mathbb{S}$,
  $\mathbb{P}$, $\mathbb{P}_{\mrm{U}}$. Then, if the set of formulas
  $\Theta$ is finite, we have that
  $ \Theta \models_{\bK} \varepsilon\approx \delta $ is equivalent to
  $\bK\models \forall x_1\dots \forall x_n (\bigwedge \Theta \to
  \varepsilon\approx \delta)$. This shows that, when the set of
  premises $\Theta$ is finite, then the consequence relation
  $\models_{\bK}$ can be encoded by the consequence relation from
  first-order logic. Notice that this ultimately justifies the
  notational conventions from
  \ref{remark:closure_quasivariety_validities} above. If the relation
  $\models_{\bK}$ is not finitary this does not need to be the case: one needs to work with so-called \emph{generalised
    quasivarieties} instead, and replace the consequence relation from
  first-order logic by the consequence relation of some suitable
  infinitary logic. Given our present interest in finitary logical
  systems we will not expand on this issue, and we simply refer the
  interested reader to \cite{Font.2016}.
\end{remark}

\section{Weak Logics and Expanded Algebras}\label{sec.1}

In this section we introduce \emph{weak logics} as a generalisation of
propositional logical systems and we provide several examples of
them. Alongside, we define \emph{expanded algebras} and \emph{core
  semantics} to provide an algebraic interpretation to these logics.

\subsection{Weak Logics}

We start by introducing the notion of weak logical systems, which is
the key object of interest of the present work. As we are interested
in logical systems which are not necessarily closed under uniform
substitution, we firstly identify the restricted class of atomic
substitutions.

\begin{definition}\label{def:atomic_substitution}
  An \emph{atomic substitution} is a substitution
  $\sigma\in\subst(\cL)$ such that $\sigma[\atoms] \subseteq
  \atoms$. We denote by $\atomsubst(\cL)$ the set of all \emph{atomic
    substitutions} in $\cL$.
\end{definition}

\begin{definition}[Weak Logic]\label{def:weak_logic}
  A \emph{weak logic} is a (finitary) consequence relation $\Vdash$
  such that, for all atomic substitutions $\sigma\in \atomsubst(\cL)$,
  $\Gamma\Vdash\phi$ entails $\sigma[\Gamma] \Vdash \sigma(\phi)$.
\end{definition}

\begin{remark}
  A weak logic is thus a consequence relation which is closed under
  atomic substitution. Intuitively, this principle reifies the least
  prerequisite a consequence relation must satisfy in order to be
  characterizable as a logic: the validity of the consequences in a
  weak logic can depend on the logical complexity of its formulas, but
  not on the specific variables that occur in them. Philosophically,
  this can be interpreted as a weakening of the Bolzanian-Tarskian
  notion of logicality.
\end{remark}

Obviously, standard logics are weak logics. More poignantly, there are
several examples of weak logics which are not standard logics and that
have been extensively studied in the literature. Their existence and
recognition constitutes the main motivation behind our interest for
this class of consequence relations and for the abstract results of
this article.

\begin{example}
  Public Announcement Logic ($ \PAL $)
  \cite{10.1007/978-3-642-24130-7_6} is an example of a modal logic
  that is not closed under uniform substitution
  \cite{Holliday2013-HOLIDA-5}. However, it can be shown that $\PAL$
  is closed under atomic substitution \cite[\S
  2.1]{Holliday2013-HOLIDA-5} and it is therefore a weak logic.
  Introducing the proper syntax and semantics of $\PAL$ is out of
  scope of this paper, but we mention the following example from
  \cite{Holliday2013-HOLIDA-5} to provide the reader with some
  intuition why uniform substitution fails. Given a set of agents
  $\cA$, the language of $\PAL$ extends the basic modal language with
  operators $K_i$, for all $i \in \cA$, and $\langle \phi \rangle$ for
  any formula $\phi$. The sentence $K_i \phi$ should be read as
  ``agent $i$ knows that $\phi$'' and $\langle \phi \rangle \psi$ as
  ``after the truthful announcement of $\phi$ to all agents, $\psi$
  holds''. Let the atoms of the language stand for \emph{facts} ---
  that is, sentences that can be truly uttered at any time. Consider
  then the principle:
\begin{gather}
  p \to \langle p \rangle p \quad\text{ (if $p$ is true, $p$ remains true after a truthful announcement)} \label{pal.axiom}  \tag{$\star$}
\end{gather}
The schema (\ref{pal.axiom}) is valid for \emph{facts}, but in general
does not hold if we substitute $p$ with a sentence talking about the
epistemic state of an agent. Let $\phi$ be the sentence
\textit{``Ljubljana became the capital of an independent Slovenia in
  1991, and agent $j$ does not know this''}, with translation
$c \land \lnot K_j c$. Now substituting $\phi$ for $p$ in
(\ref{pal.axiom}) gives us a Moorean sentence --- after truthfully
announcing $\phi$, agent $j$ learns that \textit{``Ljubljana became
  the capital of an independent Slovenia in 1991''}, and thus the
conclusion $\langle \phi \rangle \phi$ is no longer truthful.
\end{example}

\begin{example}\label{example.inquisitive}
  Logics based on team semantics, such as inquisitive and dependence
  logics \cite{ciardelli2020,Ciardelli2011-CIAIL,Yang2016-YANPLO},
  offer a rich supply of examples of weak logics.  In \cref{sec.7} we
  will focus particularly on $\inqB$, $\inqB^\otimes$, $\inqI$ and
  $\inqI^\otimes$, namely the classical and the intuitionistic
  versions of inquisitive and dependence logic. However, already now
  we can provide a conceptual motivation why $\inqB$ is not closed
  under uniform substitution. One of the main goals of $\inqB$ is to
  serve as a basis for a uniform treatment of both truth-conditional
  statements and questions in natural language. To that end, the
  intended semantics of $\inqB$ must establish when a piece of
  information \emph{supports} a statement or \emph{settles} a question
  rather than their truth conditions. We call the evidence an
  \emph{information state} and represent it as a set of possible
  worlds.

\smallskip
\noindent Let $p$ be an arbitrary statement without inquisitive
content, e.g., \textit{``It is raining in Glasgow''}. Assume that $p$
holds in the possible worlds $a$ and $b$, i.e., the information state
$\{a,b\}$ supports $p$ (see \cref{fig:inqb.possibleworlds}). We form
the polar question $?p$ --- \textit{``Is it raining in Glasgow?''},
and model it as the set of alternatives $\{a , b\}$ and
$\{c,d\}$. Let's check the validity of
$\textit{Double Negation Elimination}$ (\texttt{DNE}) ---
$\lnot\lnot q \to q$; we interpret negation as the complement of the
union of alternatives. Thus any information state supporting
$\lnot\lnot p$ will support the statement $p$ as well
(\cref{fig:inqb.possibleworlds.c}). However, this is not the case for
questions --- e.g., the state $\{b,d\}$ supports $\lnot\lnot ?p$, but
does not settle $?p$ as the possible worlds $b$ and $d$ do not agree
on a same answer. Hence we can conclude that the schema \texttt{DNE}
is valid only for statements without inquisitive content, i.e., for
propositional atoms.

\begin{figure}[h]
  \centering
  \subfigure[$p$]{\label{fig:inqb.possibleworlds.a}
    \begin{tikzpicture}[>=latex,scale=.7]
      \draw[opaque,rounded corners] (-1.7,1.7) rectangle (1.7, .2);
      \draw[rounded corners,draw opacity=0] (-1.7,-1.7) rectangle (1.7, -.2);
      \draw (-1,1) node[index gray] (yy) {$a$};
      \draw (1,1) node[index gray] (yn) {$b$};
      \draw (-1,-1) node[index gray] (ny) {$c$};
      \draw (1,-1) node[index gray] (nn) {$d$};
    \end{tikzpicture}
  }
  \hspace{.4in}
  \subfigure[$?p$]{\label{fig:inqb.possibleworlds.b}
    \begin{tikzpicture}[>=latex,scale=.7]
      \draw[opaque,rounded corners] (-1.7,1.7) rectangle (1.7, .2);
      \draw[opaque,rounded corners] (-1.7,-1.7) rectangle (1.7, -.2);
      \draw (-1,1) node[index gray] (yy) {$a$};
      \draw (1,1) node[index gray] (yn) {$b$};
      \draw (-1,-1) node[index gray] (ny) {$c$};
      \draw (1,-1) node[index gray] (nn) {$d$};
    \end{tikzpicture}
  }

  \subfigure[$\lnot\lnot p$]{\label{fig:inqb.possibleworlds.c}
     \begin{tikzpicture}[>=latex,scale=.7]
     \draw[opaque,rounded corners] (-1.7,1.7) rectangle (1.7, .2);
     \draw[rounded corners, draw opacity = 0] (-1.7,-1.7) rectangle (1.7, -.2);
     \draw (-1,1) node[index gray] (yy) {$a$};
     \draw (1,1) node[index gray] (yn) {$b$};
     \draw (-1,-1) node[index gray] (ny) {$c$};
     \draw (1,-1) node[index gray] (nn) {$d$};
   \end{tikzpicture}

  }
  \hspace{.4in}
  \subfigure[$\lnot\lnot ?p$]{\label{fig:inqb.possibleworlds.d}
    \begin{tikzpicture}[>=latex,scale=.7]
      \draw[opaque,counter,rounded corners, draw = red] (.2,1.6) rectangle (1.7, -1.6);
      \draw[opaque,rounded corners] (-1.7,1.7) rectangle (1.8, -1.7);
      \draw (-1,1) node[index gray] (yy) {$a$};
      \draw (1,1) node[index gray] (yn) {$b$};
      \draw (-1,-1) node[index gray] (ny) {$c$};
      \draw (1,-1) node[index gray] (nn) {$d$};

    \end{tikzpicture}
  }
  \caption{Double-negation elimination for statements and polar questions.}
  \label{fig:inqb.possibleworlds}
\end{figure}

\smallskip
\noindent Actually, $\inqB$ is a concrete example of a wider class of
weak logics --- a \emph{double negation atoms logic} or
\emph{$\dna$-logic}. A \emph{$\dna$-logic} (also \emph{negative
  variant} of an intermediate logic
\cite{Ciardelli:09thesis,Miglioli1989-PIESRO}) is a set of formulas
$\intlog^\neg=\{ \phi[\neg p_0,\dots,\neg p_n/p_0,\dots,p_n] : \phi
\in \intlog \}$, where $\intlog$ is an intermediate logic, namely a
logic comprised between $\ipc$ and $\cpc$. It can be proved (see
e.g. \cite[Prop. 3.2.15]{Ciardelli:09thesis}) that $\dna$-logics are
closed under atomic substitutions. However, for any $\dna$-logic
$\intlog\neq\cpc$ it is the case that
$\neg\neg p\to p\in \intlog^\neg$, but
$(\neg\neg ( p \lor \lnot p) \to p \lor \lnot p) \notin \intlog^\neg$,
showing that $\dna$-logics are not standard logics. We also notice
that $\dna$-logics can be further generalised to $\chi$-logics,
defined in \cite{quadrellaro.2021}, which provide yet another
non-trivial example of weak logics.
\end{example}

We briefly mention the following natural notions, although we will not
use them in the rest of the paper. If $\Vdash$ is a weak logic we know
that it is at least closed under all atomic substitutions
$\sigma \in \atomsubst(\cL)$, but in general there could be more
substitutions for which the logic $\Vdash$ is closed. We call such
substitutions \emph{admissible} for $\Vdash$.

\begin{definition}[Admissible Substitutions]
  Let $\Vdash$ be a weak logic. The set of \emph{admissible
    substitutions} $\admsubst(\Vdash)$ is the set of all substitutions
  $\sigma$ such that, for all sets of formulas
  $\Gamma \cup \{ \phi \} \subseteq \formulas$, $ \Gamma \Vdash \phi$
  entails $ \sigma[\Gamma] \Vdash \sigma(\phi)$.
\end{definition}

\begin{remark}
  As noticed above, we immediately have that
  $\atomsubst(\cL) \subseteq \admsubst(\Vdash)$. However, in stark
  contrast with the set of atomic substitutions, determining the set
  of admissible substitutions of a weak logic is in principle much
  harder. An example of such a characterization can be given for the
  case of inquisitive logic $\inqB$: one can in fact verify that
  $\sigma\in \admsubst(\inqB) $ if and only if $\sigma $ is a
  classical substitution, namely if for all $p\in\atoms$,
  $\sigma(p)\equiv_{\inqB}\psi$ where $\psi$ is a disjunction-free
  formula.
\end{remark}

\noindent Even if in weak logics we cannot freely substitute formulas
in place of variables, we often want to consider the subset of
formulas for which this is possible. We refer to this subset as the
core of a logic.

\begin{definition}[Core of a Logic]
  The \emph{core} of a weak logic $\Vdash$ is the set $\core(\Vdash)$
  of all formulas $\psi \in \formulas$ such that for all sets of
  formulas $\Gamma \cup \{\phi\}$ we have that:
  \[ \Gamma \Vdash \phi \implies \Gamma[\psi/x] \Vdash
    \phi[\psi/x], \]
  \noindent where $x\in \atoms$ is any atomic variable.
\end{definition}

\begin{remark}
  Equivalently, we can say that $\psi$ is a core formula of $\Vdash$
  if and only if for all $x \in \atoms$ the substitution $\sigma$ such
  that $\sigma {\restriction} \atoms \setminus\{x \} = \id_{\atoms}$
  and $\sigma(x) = \psi$ is admissible. Clearly, we always have that
  $\atoms \subseteq \core(\Vdash)$.
\end{remark}

\subsection{Expanded Algebras}

In order to make sense of weak logics from an algebraic perspective,
we need to refine the usual algebraic semantics from abstract
algebraic logic in order to handle the failure of uniform
substitution. To this end, we introduce \emph{expanded algebras} as
the expansion of standard algebras by an extra predicate symbol.

\begin{definition}[Expanded Algebra]
  Let $\cA$ be an $\cL$-algebra and $P$ a unary predicate, an
  \emph{expanded algebra} is a structure in the language
  $\cL \cup \{P\}$. We denote $P^\cA$ also by $\core(\cA)$ and we
  refer to it as the \emph{core} of the expanded algebra $\cA$.
\end{definition}

\begin{remark}
  Essentially, expanded algebras are first-order structures with
  exactly one predicate symbol of arity 1, and arbitrary many
  functional symbols. Since the relational part of the language
  consists of only one predicate, we always assume without loss of
  generality that it consists of the same symbol, so that we always
  regard any two expanded $\cL$-algebras as structures in the same
  vocabulary. When the signature $\cL$ is clear from the context we
  talk simply of algebras and expanded algebras. Notice that algebras
  augmented by a single predicate symbol often appear in abstract
  algebraic logic as \emph{matrices}
  (cf.~\cite[Ch.~4]{Font.2016}). The key difference between expanded
  algebras and matrices is in the role played by the additional
  predicate symbol. In fact, in the context of expanded algebras the
  predicate identifies the \emph{core} of the algebra, which we use to
  restrict the possible interpretation of the atomic formulas
  (cf.~\ref{def:core_semantics} below). Differently, in the context of
  logical matrices, the predicate rather identifies a designated
  ``truth set''. We will consider in Section~\ref{sec.5} a matrix-like
  approach to weak logics which unifies these two points of views.
\end{remark}

Since expanded algebras are first-order structures, we can apply the
definition of maps from Definition~\ref{maps:expanded_algebras} in
their setting. We recall that a strong homomorphism $h:\cA\to \cB$ of
two expanded algebras is an $\cL$-algebra homomorphism such that
$h[\core(\cA)]=\core(B)$. We recall from
\ref{notation:class_operators} that $\mathbb{H}(\bK)$ indicates the
closure of $\bK$ under images by strong homomorphisms.  We can then
extend the notions of quasivarieties and varieties to the setting of
expanded algebras.

\begin{definition}
  A class of expanded algebras $\bK$ is a \textit{quasivariety} if it
  is closed under $\mathbb{I}$, $\mathbb{S}$, $\mathbb{P}$ and
  $\mathbb{P}_{\mrm{U}}$.  A class of expanded algebras $\bK$ is a
  \textit{variety} if it is closed under $\mathbb{H}$, $\mathbb{S}$
  and $\mathbb{P}$. We denote by $\mathbb{Q}(\bK)$ the quasivariety
  generated by $\bK$ and by $\mathbb{V}(\bK)$ the variety generated by
  $\bK$.
\end{definition}

\begin{notation}
  Let $\bK$ be a class of expanded $\cL$-algebras, then we write
  $\bK\restriction\cL$ for the class of its $\cL$-reducts. If it is
  clear from the context, we also write $\bK$ for
  $\bK\restriction\cL$.
\end{notation}

\begin{remark}\label{remark:expanded_quasivarieties_1}
  We notice the following: if $\bK$ is a quasivariety of expanded
  $\cL$-algebras, then its $\cL$-reducts $\bK\restriction\cL$ form a
  quasivariety of $\cL$-algebras. However, if $\bK$ is a quasivariety
  of $\cL$-algebras, then it is not the case that an arbitrary
  expansion $\bK'$ of the algebras in $\bK$ gives rise to a
  quasivariety of expanded algebras. To obtain a quasivariety of
  expanded algebras we need to consider the generated quasivariety
  $\mathbb{Q}(\bK')$. We shall consider in \ref{sec.2.3} later some
  cases when this additional step is not necessary, namely the case
  when a quasivariety of algebras determines uniquely a quasivariety
  of expanded algebras.
\end{remark}

Crucially, expanded algebras allow us to define a more fine-grained
consequence relation than the one we introduced in
Definition~\ref{remark:closure_quasivariety_validities}. The key idea
is to use the core of the algebra to restrict the possible valuation
of the atomic variables of the language. To our knowledge, the idea of
restricting the possible valuations of the atomic variables to a
specific subset of an algebra first appeared in \cite{grilletti}.

\begin{definition}
  The \textit{expanded term algebra} in the signature $\cL$ is the
  structure $\Fm_\cL$ augmented with
  $\core(\Fm_\cL)=\atoms$. We often write $ (\Fm,\atoms) $ and omit
  the index $\cL$ and its signature operations when the language is
  clear from the context. A homomorphism from $\Fm_\cL$ to an expanded
  $\cL$-algebra $\cA$ is a \emph{core assignment}, i.e., it is a
  homomorphism $h:\Fm_\cL\to \cA$ such that $h(x)\in\core(\cA)$ for
  all $x\in \atoms$. We write $\homset^c(\Fm_\cL,\cA)$ for the set of
  core assignments from $\Fm_\cL$ to $A$.
\end{definition}

\begin{remark}
  We notice that, alternatively, one could also consider the expansion
  of the term algebra $\Fm$ with the core defined by letting
  $\core(\Fm_\cL)=\core(\Vdash)$. This means that the endomorphisms of
  the term algebra are all the admissible substitutions of $\Vdash$,
  and not only the atomic substitutions. As this makes only for a
  minor generalisation of our results, we stick to the former
  definition and always consider the atomic formulas as the underlying
  core of the term algebra.
\end{remark}

\begin{definition}[Core Semantics]\label{def:core_semantics}
  Let $\bK$ be a class of expanded algebras and
  $\Theta \cup \{\varepsilon\approx\delta\}$ a set of equations, then
  the \emph{equational core-consequence relation relative to $\bK$} is defined
  as follows:
  \begin{align*}
    \Theta \corevDash_{\bK} \varepsilon\approx \delta \Longleftrightarrow
    & \text{ for all } \cA \in \bK,\, h\in \homset^c(\Fm,\cA), \\
    & \text{ if } h(\epsilon_i)= h(\delta_i) \text{ for all } \epsilon_i\approx \delta_i\in \Theta,  \text{ then } h(\varepsilon)= h(\delta).
  \end{align*}
  We then write
  $\corevDash_\bK \bigwedge_{i\leq n} \varepsilon_i\approx \delta_i
  \to \varepsilon \approx \delta$ if
  $ \bigwedge_{i\leq n} \varepsilon_i\approx \delta_i \corevDash_\bK
  \varepsilon\approx \delta $.  We usually write
  $\cA \corevDash \bigwedge_{i\leq n} \varepsilon_i\approx \delta_i
  \to \varepsilon \approx \delta$ in place of
  $\corevDash_{\{\cA \}} \bigwedge_{i\leq n} \varepsilon_i\approx
  \delta_i \to \varepsilon \approx \delta$. The related notions for
  equations are defined analogously.
\end{definition}

\begin{remark}\label{remark:first-order-core}
  Crucially, if $\Theta$ is a finite set of equations
  $\{\varepsilon_i\approx \delta_i : i\leq n\}$, then we have that
  $ \bigwedge_{i\leq n} \varepsilon_i\approx \delta_i \corevDash_{\bK}
  \varepsilon \approx \delta $ holds if and only if
  \[\bK \models \forall x_0,\dots ,\forall x_m \Big( \bigwedge_{i\leq
      m} \core(x_i) \land \bigwedge_{i\leq n} \varepsilon_i\approx
    \delta_i \to \varepsilon \approx \delta \Big),  \] 
    where all variables from $\Theta\cup\{\varepsilon \approx \delta\} $ are from $x_0,\dots, x_m$. In other words,
  when $\Theta$ is finite the core semantics over a class of expanded
  algebra can be encoded in terms of the standard first-order
  consequence relation $\models$, exactly as in the case of the
  relation $\models_{\bK}$. As we stressed in
  Remark~\ref{propositional_consequence_relation_first_order}, if the
  relation $\corevDash_{\bK}$ is not finitary, then one could provide
  a similar translation into a suitable infinitary logic and replace
  quasivarieties by generalised quasivarieties.
\end{remark}

By Definition~\ref{def:core_semantics}, in core semantics atomic
variables are always assigned to core elements and, as a result of
this feature, arbitrary formulas are always interpreted inside the
subalgebra generated by the core. This motivates the interest in
core-generated structures and in quasivarieties which are generated by
these structures.

\begin{definition}
  An expanded algebra $\cA$ is \emph{core-generated} if
  $\cA = \langle \core(\cA) \rangle$. If $\bQ$ is a class of algebras,
  we write $\bQ_{\mrm{CG}}$ for its subclass of core-generated
  structures. A quasivariety $\bQ$ of expanded algebras is
  \emph{core-generated} if it is generated by its subclass of
  core-generated expanded algebras, i.e.,
  $\bQ=\mathbb{Q}(\bQ_{\mrm{CG}})$. Similarly, a variety $\bV$ of
  expanded algebras is \emph{core-generated} if
  $\bV = \mathbb{V}(\bV_{\mrm{CG}})$.
\end{definition}

The next proposition shows that unrestricted substitutions interact
nicely with core-generated quasivarieties in core semantics. In a
sense, the next result shows that core semantics in core-generated
classes of structures approximates as much as possible the standard
consequence relations defined over quasivarieties.

\begin{lemma}\label{assignment_vs_core_assignment}
  Let $A$ be a core-generated expanded algebra and let
  $h\in \homset(\Fm,\cA)$, then there are a core assignment
  $g\in \homset^c(\Fm,\cA)$ and a substitution $\sigma$ such that
  $h(\phi)=g(\sigma(\phi))$ for all $\phi\in \Fm$.
\end{lemma}
\begin{proof}
  Let $(x_i)_{i<\omega}$ be an enumeration of $\atoms$. Since $\cA$ is
  core-generated, we have that for all $x_i\in \atoms$ there is a
  polynomial $t_i$ such that $h(x_i) = t_i(a^i_0,\dots,a^i_{n_i})$
  with $a^i_j\in\core(\cA)$ for all $j\leq n_i$. Let now
  $\{ y^i_j : i<\omega, \; j\leq n_i \}$ be another enumeration of
  $\atoms$ and define $g \in \homset^c(\Fm,\cA)$ such that
  $g(y^i_j) = a^i_j$ for all $ i<\omega$, $j\leq n_i$. In particular,
  this means that for all formulas $\psi(x_0, \dots, x_m)\in \Fm$ we
  have
  \[ h(\psi(x_0, \dots x_m)) = g(\psi( t_0(y^0_0,\dots,y^0_{n_0}),
    \dots, t_m(y^m_0,\dots,y^m_{n_m}) ) ).\] By construction we have
  that $g$ is a core assignment. Let $\sigma$ be the substitution
  defined by letting $\sigma(x_i) = t_i(y^i_0,\dots,y^i_{n_i})$ for
  all $x_i\in \atoms$, then from the display above we derive that
  $h(\phi)=g(\sigma(\phi))$ for all $\phi\in \Fm$.
\end{proof}

\begin{proposition}\label{th:subst-core} Let $\bQ$ be a core-generated
  quasivariety of expanded algebras, then
  $\sigma(\Theta) \corevDash_{\bQ} \sigma(\varepsilon \approx \delta)$
  for all $\sigma \in \subst(\cL)$ holds if and only if
  $ \Theta \models_{\bQ} \varepsilon \approx \delta$.
\end{proposition}
\begin{proof}
  We consider the left-to-right direction. Suppose that
  $\Theta \not\models_\bQ \varepsilon \approx \delta$. Since $\bQ$ is
  core-generated, there is a core-generated expanded algebra
  $\cA \in \bQ$ and some $h\in \homset(\Fm,\cA)$ such that
  $\cA\models h(\Theta)$ and
  $\cA\not\models h(\varepsilon) \approx h(\delta)$ (recall
  Definition~\ref{remark:closure_quasivariety_validities}). From
  Lemma~\ref{assignment_vs_core_assignment} it follows that there is a
  core assignment $g\in\homset^c(\Fm,\cA)$ and a substitution $\sigma$
  such that $h(\phi)=g(\sigma(\phi))$ for all $\phi\in \Fm$. Then it
  follows that $\cA \models  g(\sigma[\Theta])$ but
  $\cA \not\models g(\sigma(\varepsilon \approx \delta))$, whence
  $\sigma(\Theta) \not\corevDash_{\bQ} \sigma(\varepsilon \approx
  \delta)$.

  \smallskip
  \noindent We consider the right-to-left direction. Suppose that
  $\sigma(\Theta) \ncorevDash_\bQ \sigma(\varepsilon \approx \delta)$,
  thus we can find an algebra $A\in \bQ$ and an assignment
  $h \in \homset^c(\Fm, \cA)$ such that $\cA \models h(\sigma[\Theta])$
  and $\cA \not\models h(\sigma(\varepsilon\approx \delta))$. Then we
  have that $h \circ \sigma\in \homset(\Fm,\cA)$, 
  whence $ \Theta \not \models_{\bQ} \varepsilon \approx \delta$.
\end{proof}

Since core semantics only looks at the substructure
$\langle \core(\cA) \rangle$ of an expanded algebra $\cA$, it make
sense to consider extensions of $A$ that preserve the core.

\begin{definition}\label{def:core_superstructure}
  We say that $\cB$ is a \emph{core superalgebra} of $\cA$ if
  $ B\in \Mod(\cL)$, $\cA \leq \cB$ (i.e., $\cA$ is a subalgebra of
  $\cB$) and $\core(\cA) = \core(\cB)$. If $\bK$ is class of expanded
  algebras, we write $\mathbb{C}(\bK) $ for the class of all core
  superalgebras of elements of $\bK$.
\end{definition}

\noindent The next lemma and the following proposition establish some
fundamental facts about the relation $\corevDash$. Importantly, they
show that the validity of quasiequations is preserved both under the
operators $\mathbb{I}$, $\mathbb{S}$, $\mathbb{P}$,
$\mathbb{P}_{\mrm{U}} $, and also under the core superstructure
operator $\mathbb{C}$ from
Definition~\ref{def:core_superstructure}. Additionally, we also show
that equations are preserved under strong homomorphisms. The following
lemma is essentially a rephrasing of
Remark~\ref{remark:first-order-core}.

\begin{lemma}\label{semantic:translation}
  Let $\bK$ be a class of expanded algebras, let
  $\bigwedge_{i\leq n} \varepsilon_i \approx \delta_i\to \varepsilon
  \approx \delta $ be a quasiequation and let $V$ be all the variables
  occurring in it, then
  \begin{align*}
    \{\varepsilon_i \approx \delta_i: i \leq n \} \corevDash_{\bK} \varepsilon \approx \delta \; \Longleftrightarrow \; \{\varepsilon_i \approx \delta_i : i \leq n \} \cup \bigcup_{x\in V}\core(x)	\models_{\bK} \varepsilon \approx \delta.
  \end{align*}
\end{lemma}
\begin{proof}
  This follows immediately from the fact that an assignment
  $h\in \homset(\Fm, \cA)$ is a core assignment if and only if
  $h(x)\in \core(\cA)$ for all $x\in \atoms$.
\end{proof}

\begin{proposition}\label{preservation.thm}
  Let $\bK$ be a class of expanded algebras, then the following hold:
  \begin{enumerate}[(1)]
  \item\label{preservation.thm.pt2} let
    $\bigwedge_{i\leq n} \varepsilon_i \approx \delta_i \to
    \varepsilon \approx \delta $ be a quasi-equation, then for all
    $\mathbb{O}\in \{ \mathbb{I},
    \mathbb{S},\mathbb{P},\mathbb{P}_{\mrm{U}}, \mathbb{C} \}$ we have
    that
    \begin{align*}
      \{\varepsilon_i \approx \delta_i : i \leq n \}\corevDash_{\bK} \varepsilon \approx \delta  \; \Longrightarrow \;  \{\varepsilon_i \approx \delta_i : i \leq n \} \corevDash_{\mathbb{O}(\bK)} \varepsilon \approx \delta;
    \end{align*}
  \item\label{preservation.thm.pt3} if
    $\corevDash_{\bK} \varepsilon \approx \delta$, then
    $ \corevDash_{\mathbb{H}(\bK)} \varepsilon \approx \delta $.
  \end{enumerate}
\end{proposition}
\begin{proof}
  For Claim (1) recall that the validity of universal Horn formulas is
  preserved by the operators
  $\mathbb{O}\in \{ \mathbb{I},
  \mathbb{S},\mathbb{P},\mathbb{P}_{\mrm{U}} \}$. Then
  Lemma~\ref{semantic:translation} implies that the validity of
  formulas in core semantics is preserved by
  $\mathbb{O}\in \{ \mathbb{I},
  \mathbb{S},\mathbb{P},\mathbb{P}_{\mrm{U}} \}$. Preservation of
  validity by $\mathbb{C}$ follows immediately by the definitions of
  core semantics and core superstructure. Claim (2) is immediate by
  the definition of core semantics and the fact that, if
  $B\in \mathbb{H}(\bK)$, then there is $A\in\bK$ and a surjective
  homomorphism $h:A\to B$ such that $\core(B)=h[\core(A)]$.
\end{proof}

We conclude this section by showing a version of Maltsev's Theorem for
the setting of core-generated quasivarieties, i.e., we prove that
every core-generated quasivariety is axiomatised by its validities
under core semantics.

\begin{definition} We define the following notions.
  \begin{enumerate}[(1)]
  \item For any set $T$ of quasiequations (or equations), we let
    $\coremodels{T}$ be the class of expanded algebras $\cA$ such that
    $\cA \corevDash T$ and $\coregenmodels{T}$ for its subclass of
    core-generated models.
  \item For any class $\bK$ of expanded algebras, we denote by
    $\coretheories{\bK}$ the set of all quasiequations true in $\bK$
    under core semantics, and we denote by $\coreequations{\bK}$ the
    set of all equations true in $\bK$ under core semantics.
  \end{enumerate}
\end{definition}

\noindent The following proposition is an immediate corollary of
Maltsev's Theorem (cf. Fact~\ref{birkhoff_maltsev}) and
Proposition~\ref{th:subst-core}.

\begin{proposition}\label{thm:maltsev}
  Let $\bQ$ be a quasivariety of expanded algebras and let $\cA$ be a
  core-generated expanded algebra. Then $\cA \in \coregen{\bQ}$ if and
  only if $\cA \corevDash \coretheories{\bQ}$.
\end{proposition}
\begin{proof}
  The direction from left-to-right follows immediately by the
  definition of $\coretheories{\bQ}$. We consider the direction from
  right-to-left. Let $\cA$ be core-generated and suppose
  $A\notin \coregen{\bQ}$. It follows that $A\notin \bQ$ and so by
  Maltsev's Theorem $\cA \not\models \theory{\bQ}$. Let
  $\bigwedge_{i\leq n} \varepsilon_i \approx \delta_i \to \alpha
  \approx \beta \in \theory{\bQ}$ be such that
  $\cA\not \models \bigwedge_{i\leq n} \varepsilon_i \approx \delta_i
  \to \alpha \approx \beta$. By Proposition~\ref{th:subst-core} it
  follows that there is a substitution $\sigma$ such that
  $\cA\ncorevDash \sigma(\bigwedge_{i\leq n} \varepsilon_i \approx
  \delta_i \to \alpha \approx \beta)$. Now, since the relation
  $\models_{\bQ}$ is closed under uniform substitution, it follows in
  particular that
  $\sigma(\bigwedge_{i\leq n} \varepsilon_i \approx \delta_i \to
  \alpha \approx \beta) \in \theory{\bQ}$. Moreover, since obviously
  $\theory{\bQ}\subseteq\coretheories{\bQ} $, we obtain that
  $\sigma(\bigwedge_{i\leq n} \varepsilon_i \approx \delta_i \to
  \alpha \approx \beta) \in \coretheories{\bQ}$. We conclude that
  $\cA \ncorevDash \coretheories{\bQ}$.
\end{proof}

\section{Algebraizability of Weak Logics}\label{sec.2}

In this section we use quasivarieties of expanded algebras to provide
three different notions of algebraizability for the setting of weak
logics, which we respectively call \emph{loose algebraizability},
\emph{strict algebraizability} and \emph{fixed-point
  algebraizability}. We proceed as follows: firstly, we review the
notion of algebraizability for the setting of standard logics in
\cref{sec.2.1}, then in \cref{sec.2.2} we introduce the loose notion
of algebraizability for weak logics, in \cref{sec.2.3} we consider the
strict version of algebraizability for weak logics, and finally in
\cref{sec.2.4} we focus on the very specific notion of fixed-point
algebraizability.

\subsection{Algebraizability of Standard Logics}\label{sec.2.1}

We first review in this section the notion of algebraizability for the
setting of standard logics. Recall that we denote by $\Fm$ and $\Eq$
respectively the set of formulas and equations in the signature
$\cL$. We shall refer to two maps, called transformers
\cite{Font.2016}, that allow us to translate formulas into equations
and \textit{vice versa}.

\begin{definition}
  A pair of \emph{transformers} in the language $\cL$ is a pair of
  functions $\tau : \Fm \to \wp(\Eq)$ and $\Delta: \Eq \to
  \wp(\Fm)$. We say that $\tau$ and $\Delta$ are \emph{structural} if
  for all substitutions $\sigma \in \subst(\cL)$,
  $\tau(\sigma(\phi)) = \sigma(\tau(\phi))$ and
  $\sigma(\Delta(\varepsilon,\delta)) =
  \Delta(\sigma(\varepsilon,\delta))$, where we let
  $\sigma(\varepsilon,\delta) =
  (\sigma(\varepsilon),\sigma(\delta))$. We say that the transformers
  $\tau, \Delta$ are \emph{finitary} if for all $\phi\in\Fm$ and
  $\varepsilon\approx \delta \in \Eq$, $|\tau(\phi)|<\aleph_0$ and
  $|\Delta(\varepsilon,\delta)|<\aleph_0$. For any set of formulas
  $\Gamma\subseteq \Fm$, we let
  $\tau(\Gamma) = \bigcup_{\phi\in\Gamma} \tau(\phi)$ and for all set
  of equations $\Theta\subseteq \Eq$ we let
  $\Delta(\Theta) = \bigcup_{\varepsilon \approx \delta
    \in\Theta}\Delta(\varepsilon,\delta)$.
\end{definition}

\noindent The following notion of algebraizability was introduced by
Blok and Pigozzi in their seminal article \cite{Blok1989-BLOAL}. We
write $\Gamma_0 \equiv_{\bK} \Gamma_1$ as a shorthand for
$\Gamma_0 \vDash_{\bK} \Gamma_1$ and $\Gamma_1 \vDash_{\bK} \Gamma_0$.

\begin{definition}[Algebraizability]\label{standard.Alg}
  A (standard) logic $\vdash$ is \emph{algebraizable} if there are a
  quasivariety of (standard) algebras $\bQ$ and structural
  transformers $\tau : \Fm \to \wp(\Eq)$ and
  $\Delta: \Eq \to \wp(\Fm)$ such that:
  \begin{align}
    \Gamma \vdash \phi &\Longleftrightarrow \tau[\Gamma]\models_{\bQ} \tau(\phi) \tag{A1} \label{s.alg1} \\
    \Delta[\Theta] \vdash \Delta(\varepsilon,\delta) &\Longleftrightarrow \Theta\models_{\bQ} \varepsilon\approx\delta \tag{A2} \label{s.alg2} \\
    \phi &\dashv\vdash \Delta[\tau(\phi)] \tag{A3} \label{s.alg3} \\
    \varepsilon \approx \delta  &\equiv_{\bQ}  \tau[ \Delta(\varepsilon,\delta)]. \tag{A4} \label{s.alg4}
  \end{align}
  \noindent If these conditions are met, we then say that $\bQ$ is an
  \emph{equivalent algebraic semantics} of $\vdash$, and that
  $(\bQ,\tau,\Delta)$ witnesses the algebraizability of
  $\vdash$.
\end{definition}

\begin{remark}
  The restriction to quasivarieties of algebras stems from the fact
  that we are exclusively considering finitary consequence relations
  $\vdash$. See \cite[\S 3]{Font.2016} for an extension of this
  definition to generalised quasivarieties and non-finitary logics.
\end{remark}

\noindent We recall the following two important facts about
algebraizability. We refer the reader to \cite[Prop. 3.12,
Thm. 3.37]{Font.2016} for a proof of these results.

\begin{fact}\label{facts.algebraizability} Let $\vdash$ be a standard logic, then:
  \begin{enumerate}[(1)]
  \item if Conditions \ref{standard.Alg}(\ref{s.alg1}) and
    \ref{standard.Alg}(\ref{s.alg4}), or
    \ref{standard.Alg}(\ref{s.alg2}) and
    \ref{standard.Alg}(\ref{s.alg3}) hold, then $\vdash$ is
    algebraized by $(\bQ, \tau, \Delta)$;
  \item if $\vdash$ is algebraizable, then there are two finitary
    transformers $\tau, \Delta$ witnessing this fact.
  \end{enumerate}
\end{fact}

\noindent It is a key property of algebraizability that the equivalent
algebraic semantics of a standard logic is unique. See
\cite{Font.2016} for a proof of this result.

\begin{fact}\label{uniqueness.standard.alg}
  If the tuples $(\bQ_0, \tau_0, \Delta_0)$ and
  $(\bQ_1, \tau_1, \Delta_1)$ both witness the algebraizability of a
  standard logic $\vdash$, then:
  \begin{enumerate}[(1)]
  \item $\bQ_0= \bQ_1$;
  \item
    $\Delta_0(\varepsilon,\delta) \dashv\vdash
    \Delta_1(\varepsilon,\delta)$, for all
    $\varepsilon\approx \delta\in \Eq$;
  \item $\tau_0(\phi) \equiv_{\bQ_i} \tau_1(\phi)$ with
    $i\in \{0,1\}$, for all $\phi\in \Fm$.
  \end{enumerate}
\end{fact}

\subsection{Loose Algebraizability of Weak Logics}\label{sec.2.2}

By using the core consequence relation $\corevDash$ in place of the
standard one $\models$, we can provide a first version of the notion
of algebraizability in the setting of weak logics. Obviously, we write
$\Gamma_0 \equiv^c_{\bK} \Gamma_1$ as a shorthand for
$\Gamma_0 \corevDash_{\bK} \Gamma_1$ and
$\Gamma_1 \corevDash_{\bK} \Gamma_0$.

\begin{definition}\label{loose_Alg}
  A weak logic $\Vdash$ is \emph{loosely algebraizable} if there is a
  core-generated quasivariety of expanded algebras $\bQ$ and two
  structural transformers $\tau : \Fm \to \wp(\Eq)$ and
  $\Delta: \Eq \to \wp(\Fm)$ such that:
  \begin{align}
    \Gamma \Vdash \phi &\Longleftrightarrow \tau[\Gamma]\corevDash_{\bQ} \tau(\phi) \tag{W1} \label{l.alg1} \\
    \Delta[\Theta] \Vdash \Delta(\varepsilon,\delta) &\Longleftrightarrow \Theta\corevDash_{\bQ} \varepsilon\approx\delta \tag{W2} \label{l.alg2} \\
    \phi &\VdashR\,\Vdash \Delta[\tau(\phi)] \tag{W3} \label{l.alg3} \\
    \varepsilon \approx \delta  &\equiv^c_{\bQ}  \tau[ \Delta(\varepsilon,\delta)]. \tag{W4} \label{l.alg4}
  \end{align}
  \noindent If these conditions are met, we then say that $\bQ$ is a 
  \emph{loose algebraic semantics} of $\Vdash$, and that the tuple  $(\bQ,\tau,\Delta)$ witnesses the loose algebraizability of
  $\Vdash$.
\end{definition}

To show this definition is indeed robust, we adapt the proof of the
uniqueness of the equivalent algebraic semantic of standard logics
\cite[Thm. 3.17]{Font.2016} to the case of loosely algebraizable weak
logics. In the light of the following theorem, we often refer to the tuple $(\bQ, \tau, \Delta)$ witnessing the loose algebraizability of the weak logic $\Vdash$ as \emph{the equivalent loose algebraic semantics} of $\Vdash$.

\begin{theorem}\label{uniquenessalg_1}
  If both the tuples $(\bQ_0, \tau_0, \Delta_0)$ and
  $(\bQ_1, \tau_1, \Delta_1)$ witness the algebraizability of
  $\Vdash$, then:
  \begin{enumerate}[(1)]
  \item $\bQ_0= \bQ_1$;
  \item
    $\Delta_0(\varepsilon,\delta) \VdashR\,\Vdash
    \Delta_1(\varepsilon,\delta) $, for all
    $\varepsilon\approx \delta\in \Eq$;
  \item $\tau_0(\phi)\equiv^c_{\bQ_i}\tau_1(\phi) $ with
    $i\in \{0,1\}$, for all $\phi\in\Fm$.
  \end{enumerate}
\end{theorem}
\begin{proof}
  Notice that the two witnesses of algebraizability give rise to two
  different consequence relations, which we shall denote by
  $\models^0_{\bQ_0}$ and $\models^1_{\bQ_1}$.

  \medskip
  \noindent Clause (2):
  $\Delta_0(\varepsilon,\delta) \VdashR\,\Vdash
  \Delta_1(\varepsilon,\delta) $, for all
  $\varepsilon\approx \delta\in \Eq$.  \newline We prove
  $\Delta_0(\varepsilon,\delta) \Vdash \Delta_1(\varepsilon,\delta) $.
  Let $\phi \in \Delta_1(\varepsilon,\delta) $, then we clearly have
  that:
  \[\tau_0(\phi(\varepsilon,\varepsilon)),
    \phi(\varepsilon,\varepsilon)\approx \phi(\varepsilon,\delta)
    \models^0_{\bQ_0} \tau_0(\phi(\varepsilon,\delta)).\]
  \noindent By \ref{loose_Alg}(\ref{l.alg2}) it follows:
  \[\Delta_0(\tau_0(\phi(\varepsilon,\varepsilon))),
    \Delta_0(\phi(\varepsilon,\varepsilon), \phi(\varepsilon,\delta))
    \Vdash \Delta_0(\tau_0(\phi(\varepsilon,\delta))),\]
  \noindent hence, by \ref{loose_Alg}(\ref{l.alg3}):
  \begin{equation}\tag{a}
    \phi(\varepsilon,\varepsilon), \Delta_0(\phi(\varepsilon,\varepsilon), \phi(\varepsilon,\delta)) \Vdash \phi(\varepsilon,\delta).
  \end{equation}
  Now, we have that
  $\varnothing \models^1_{\bQ_1} \varepsilon\approx \varepsilon $,
  hence by $\phi\in \Delta_1(\varepsilon,\varepsilon)$ and
  \ref{loose_Alg}(\ref{l.alg2}), we obtain:
  \begin{equation}\tag{b}
    \varnothing \Vdash \phi(\varepsilon,\varepsilon).
  \end{equation}
  Moreover, it also follows that
  $\varepsilon\approx \delta \models^0_{\bQ_0}
  \phi(\varepsilon,\varepsilon)\approx \phi(\varepsilon,\delta)$,
  hence by \ref{loose_Alg}(\ref{l.alg2}):
  \begin{equation}\tag{c}
    \Delta_0(\varepsilon,\delta)\Vdash \Delta_0(\phi(\varepsilon,\varepsilon),\phi(\varepsilon,\delta)).
  \end{equation}

  \noindent Finally, by (a), (b) and (c), it follows that
  $\Delta_0(\varepsilon,\delta) \Vdash \phi(\varepsilon,\delta) $,
  hence
  $\Delta_0(\varepsilon,\delta) \Vdash \Delta_1(\varepsilon,\delta)
  $. The converse direction is proven analogously.

  \medskip
  \noindent Clause (1): $\bQ_0= \bQ_1$.  \newline We first prove that
  $\bQ_0$ and $\bQ_1$ satisfy the same quasiequations under core
  semantics.  We show only that
  $\coretheories{\bQ_0} \subseteq \coretheories{\bQ_1}$, as the other
  direction follows analogously.  Let
  $\bigwedge_{i\leq n} \varepsilon_i \approx \delta_i \to
  \varepsilon\approx \delta \in \coretheories{\bQ_0} $, then it
  follows that
  $\bigwedge_{i\leq n} \varepsilon_i \approx \delta_i
  \models^0_{\bQ_0} \varepsilon \approx \delta$ and this yields that
  $ \bigcup_{i\leq n} \Delta_0(\varepsilon_i,\delta_i) \Vdash
  \Delta_0(\varepsilon, \delta)$ by \ref{loose_Alg}(\ref{l.alg2}).  By
  point (3) above, it follows that
  $\bigcup_{i\leq n} \Delta_1(\varepsilon_i,\delta_i) \Vdash \Delta_1
  (\varepsilon, \delta) $, hence by \ref{loose_Alg}(\ref{l.alg2}) we
  get
  $\bigwedge_{i\leq n} \varepsilon_i \approx \delta_i
  \models^1_{\bQ_1} \varepsilon \approx \delta$.  The latter   entails
  $\bigwedge_{i\leq n} \varepsilon_i \approx \delta_i \to \varepsilon
  \approx \delta \in \coretheories{\bQ_1} $ and thus
  $\coretheories{\bQ_0}\subseteq \coretheories{\bQ_1}$.  By reasoning
  analogously we obtain that
  $\coretheories{\bQ_1}\subseteq \coretheories{\bQ_0} $, hence
  $\coretheories{\bQ_0}= \coretheories{\bQ_1} $.  It then follows by
  Proposition~\ref{thm:maltsev} above that
  $(\bQ_0)_{\mrm{CG}}=(\bQ_1)_{\mrm{CG}}$ and since both $\bQ_0 $ and
  $\bQ_1 $ are core-generated we conclude that $\bQ_0=\bQ_1$.

  \medskip
  \noindent Clause (3): $\tau_0(\phi)\equiv^c_{\bQ_i}\tau_1(\phi) $
  with $i\in \{0,1\}$, for all $\phi\in\Fm$.  \newline By Clause (1),
  it suffices to prove that
  $\tau_0(\phi) \equiv_{\bQ_0}^0 \tau_1(\phi)$.  By
  \ref{loose_Alg}(\ref{l.alg3}), we have
  $\Delta_0(\tau_0(\phi)) \VdashR\,\Vdash \Delta_1(\tau_1(\phi))$ and
  by Clause (2) this is equivalent to
  $\Delta_0(\tau_0(\phi)) \VdashR\,\Vdash \Delta_0(\tau_1(\phi))$.  It
  then follows by \ref{loose_Alg}(\ref{l.alg2}) that
  $\tau_0(\phi) \equiv_{\bQ_0}^0 \tau_1(\phi)$.
\end{proof}

We have thus established that every algebraizable weak logic has a
unique equivalent algebraic semantics, up to equivalence under the
core consequence relation. We conclude this section by proving an
analogue of Fact~\ref{facts.algebraizability} for weak logics.

\begin{proposition}\label{facts.weak.algebraizability} Let $\Vdash$ be
  a weak logic, then:
  \begin{enumerate}[(1)]
  \item if Conditions \ref{loose_Alg}(\ref{l.alg1}) and
    \ref{loose_Alg}(\ref{l.alg4}), or Conditions
    \ref{loose_Alg}(\ref{l.alg2}) and \ref{loose_Alg}(\ref{l.alg3})
    hold, then $\Vdash$ is loosely algebraized by
    $(\bQ, \tau, \Delta)$;
  \item if $\Vdash$ is loosely algebraizable, then there are two
    finitary transformers $\tau, \Delta$ witnessing this fact.
  \end{enumerate}
\end{proposition}
\begin{proof} We prove (1). Suppose $\Vdash$ is a weak logic and let
  $(\bQ, \tau, \Delta)$ satisfy \ref{loose_Alg}(\ref{l.alg1}) and
  \ref{loose_Alg}(\ref{l.alg4}). We verify that
  \ref{loose_Alg}(\ref{l.alg2}) and \ref{loose_Alg}(\ref{l.alg3}) hold
  as well.  By \ref{loose_Alg}(\ref{l.alg4})
  $\Theta\corevDash_{\bQ} \varepsilon\approx\delta$ is equivalent to
  $\tau[\Delta(\Theta)]\corevDash_{\bQ}\tau(\Delta(\varepsilon,\delta))
  $, which by \ref{loose_Alg}(\ref{l.alg1}) is equivalent to
  $\Delta(\Theta)\Vdash\Delta(\varepsilon,\delta) $, proving
  \ref{loose_Alg}(\ref{l.alg2}). Also, for any formula $\phi$, we have
  that $ \tau(\phi)\equiv^c_\bQ \tau(\phi)$, hence by
  \ref{loose_Alg}(\ref{l.alg4})
  $\tau[\Delta( \tau(\phi))]\equiv^c_\bQ \tau(\phi)$ and by
  \ref{loose_Alg}(\ref{l.alg1})
  $\Delta( \tau(\phi))\VdashR\,\Vdash \phi$, proving
  \ref{loose_Alg}(\ref{l.alg3}). If $(\bQ, \tau, \Delta)$ satisfies
  \ref{loose_Alg}(\ref{l.alg2}) and \ref{loose_Alg}(\ref{l.alg3}),
  then we proceed analogously.

  \smallskip
  \noindent We prove (2). For any $\phi\in \Fm$, we have by
  \ref{loose_Alg}(\ref{l.alg3}) that
  $\phi \VdashR\,\Vdash \Delta[\tau(\phi)] $, thus by $\Vdash$ being
  finitary there is some $\tau_0(\phi)\subseteq \tau(\phi) $ such that
  $|\tau_0(\phi)|<\aleph_0$ and
  $ \phi \VdashR\,\Vdash \Delta[\tau_0(\phi)] $. Moreover, for any
  equation $\varepsilon \approx \delta\in \Eq$, we have by
  \ref{loose_Alg}(\ref{l.alg4}) that
  $\varepsilon \approx \delta \equiv^c_{\bQ} \tau[
  \Delta(\varepsilon,\delta)]$, so we obtain by finitarity a finite
  subset
  $\Delta_0(\varepsilon,\delta) \subseteq \Delta(\varepsilon,\delta)$
  such that
  $\varepsilon \approx \delta \equiv^c_{\bQ} \tau[
  \Delta_0(\varepsilon,\delta)]$. Finally, it follows by the choice of
  $\tau_0,\Delta_0$ that
  $ \Delta[\tau(\phi)] \VdashR\,\Vdash \Delta(\tau_0(\phi))$, hence
  $\tau(\phi)\equiv^c_{\bQ} \tau_0(\phi)$. Similarly, from
  $\tau[ \Delta(\varepsilon,\delta)] \equiv^c_{\bQ} \tau[
  \Delta_0(\varepsilon,\delta)]$ we obtain
  $\Delta_0(\varepsilon,\delta)\VdashR\,\Vdash
  \Delta(\varepsilon,\delta)$. Thus $\tau_0,\Delta_0$ and $\bQ$
  witness the algebraizability of $\Vdash$.
\end{proof}

\begin{remark}\label{remark:problem_loose_alg}
  As witnessed by the previous results, loose algebraizability
  satisfies the same uniqueness property of standard
  algebraizability. However, we believe this is too weak of a notion,
  as it does not really meet the fundamental intuition behind
  algebraizability. In fact, in contrast to the matrix semantics of
  logics (which we shall explore later in \cref{sec.5}), the
  fundamental aspect of algebraizability is that it allows us to
  translate logical systems into equational consequences defined by
  classes of algebras, and not by arbitrary first-order structures
  containing relations. We will achieve this by considering strict
  algebraizability in the next section.
\end{remark}

\subsection{Strict Algebraizability of Weak Logics}\label{sec.2.3}

We introduce strict algebraizability as a refined notion of
algebraizability for weak logical systems. As we mentioned in
Remark~\ref{remark:problem_loose_alg}, the problem with loose
algebraizability is that it relates weak logics to (universal Horn)
classes of first-order structures, and not really to algebras. To
overcome this issue, we need to look for ways in which one can, so to
speak, \emph{eliminate} the core predicate $\core(\cA)$. The key idea
is to restrict attention to classes of expanded algebras in which the
core is already definable in the functional part of the
signature. Furthermore, since we are dealing with quasivarieties and
we want the validity of formulas to be preserved by the operators
$\mathbb{I}$, $\mathbb{S}$, $\mathbb{P}$ and $\mathbb{P}_{\mrm{U}}$,
we exclusively consider definability by means of equations. We make
this explicit in the following definitions.

\begin{notation}
  If $\cA$ is an $\cL$-algebra and $\Sigma(x)$ a set of equations in
  the variable $x$, we let
  $\Sigma(\cA)=\{a\in \cA : \cA \models \varepsilon(a) \approx
  \delta(a) \text{ for all } \varepsilon\approx\delta\in\Sigma \}$.
\end{notation}

\begin{definition}
  An expanded algebra $\cA$ is said to have an \emph{equationally
    definable core} if there is a finite set of equations $\Sigma$ in
  the variable $x$ such that $\core(\cA)=\Sigma(\cA)$.  A class of
  expanded algebras $\bK$ is said to have a \emph{(uniformly)
    equationally definable core} if there is a finite set of equations
  $\Sigma$ such that $\core(\cA)=\Sigma(\cA)$ for all $\cA \in \bK$.
\end{definition}

The key reason that explains our interest in quasivarieties of
expanded algebras with an equationally definable core is that, given
an algebra $\cA$ and a finite set of equations $\Sigma$, there is a
\emph{unique} way to expand $\cA$ into an expanded algebra with core
defined by $\Sigma$. As the following proposition shows, this provides
us with a canonical expansion of $\bQ$ into a quasivariety of expanded
algebras with core defined by $\Sigma$.

\begin{proposition}\label{projection_commutes}
  Let $\bQ$ be a quasivariety of $\cL$-algebras and $\Sigma(x)$ a
  finite set of equations, then the class of structures
  $(A,\core(\cA))$ with $A\in \bQ$ and $ \core(\cA)=\Sigma(\cA) $ is a
  quasivariety of expanded algebras.
\end{proposition}
\begin{proof}
  Let $\bQ$ be a quasivariety of $\cL$-algebras and let $\bK$ be the
  class of structures $(A,\core(\cA))$ with $A\in \bQ$ and
  $ \core(\cA)=\Sigma(\cA) $. Let
  $\{\alpha_i\approx \beta_i : i\leq n\}$ enumerate all equations in
  $\Sigma$. Then consider the formulas
  \begin{align*}
    \Phi_i\coloneqq \; &\forall x(\core(x)\to \alpha_i(x)\approx \beta_i(x) )\\
    \Psi\coloneqq \; &\forall x (\bigwedge_{i\leq n}\alpha_i(x)\approx \beta_i(x) \to \core(x)).
  \end{align*}
  Clearly, both $\Psi$ and each $\Phi_i$ are universal Horn sentences,
  thus by Fact~\ref{birkhoff_maltsev} they are preserved by the
  quasivariety operators $\mathbb{I}$, $\mathbb{S}$, $\mathbb{P}$ and
  $\mathbb{P}_{\mrm{U}}$. It follows that $\core(\cA)=\Sigma(\cA)$ for
  all $A\in \mathbb{Q}(\bK)$, and thus that $\bK$ is already a
  quasivariety.
\end{proof}

\begin{remark}
  Crucially, the previous proposition cannot be extended to varieties
  of expanded algebras. In particular, if $\bK$ is a class of expanded
  algebras with $\core(A)=\Sigma(A)$ for all $A\in \bK$, it is not
  necessarily the case that $\core(A)=\Sigma(A)$ for all
  $A\in \mathbb{H}(\bK)$. For example, consider the equation
  $x\approx \neg x$, let $A$ be any Boolean algebra with size
  $|A|\geq 2$ and let $B$ be the trivial one-element Boolean
  algebra. Clearly, there is a surjective homomorphism
  $h:A\twoheadrightarrow B$ and, at the same time, we also have that
  $\Sigma(\cB)=B$ and $\Sigma(\cA)=\emptyset$, showing that
  $\core(\cB)\neq h(\core(A))$. This shows that we cannot replace
  quasivarieties by varieties in the statement of
  Proposition~\ref{projection_commutes}. We shall see later that this
  is possible in the restricted setting of fixed-point
  algebraizability.
\end{remark}

It is easy to find concrete examples of quasivarieties of expanded
algebras whose core is equationally definable. We list here some
examples which are determined by some very basic equations.

\begin{example}
  Let $M$ be a monoid in $\cL=(\cdot,e)$ and define
  $\core(\cM) = \{a \in \cM : a^n = e^\cM \}$, then $\cM$ is an
  expanded algebra with core defined by $\Sigma=\{ x^n\approx e \}$.
\end{example}

\begin{example}
  Let $\HA$ be the variety of Heyting algebras and for all
  $\cA\in \HA$ let
  $\core(\cA) = \cA_\neg =\{ x\in \cA : \cA \models x \approx \neg\neg
  x \}$, i.e., the core of $\cA$ is its subset of regular
  elements. Let $\mathbf{ML}$ be the variety of all Medvedev algebras,
  then $\mathbf{ML}$ is generated by its subclass of core-generated
  Heyting algebras $\cA$ with core $\cA_\neg$ (cf. \cite{grilletti,
    Quadrellaro.2019B} and Fact~\cref{well-known-facts} below). As we
  shall see in \cref{sec.7}, this class plays an important role in the
  semantics of classical inquisitive propositional logic.
\end{example}

We define the strict version of algebraizability for weak logics by
requiring that the core-generated algebras corresponding to a weak
logic $\Vdash$ have the core defined by a finite set of equations
$\Sigma$. This definition is significantly stronger than the one we
introduced in the previous section, and it essentially rephrase the
weak logic $\Vdash$ in terms of the relative consequence relation
$\models_{\bQ}$ of the corresponding quasivariety of algebras.

\begin{definition}\label{strict.Alg}
  A weak logic $\Vdash$ is \emph{strictly algebraizable} if it is
  loosely algebraized (in the sense of Definition~\ref{loose_Alg}) by
  $(\bQ,\tau,\Delta)$ and, additionally, there is a finite set of
  equations $\Sigma$ defining the core of the expanded algebras in
  $\bQ$. We then say that $\Vdash$ is strictly algebraized by
  $(\bQ,\Sigma,\tau,\Delta)$ and call $(\bQ,\Sigma,\tau,\Delta)$ its
  \emph{(equivalent) strict algebraic semantics}.
\end{definition}

\noindent We can then strengthen \cref{uniquenessalg_1} to the present
context of core-generated quasivarieties with a definable core.

\begin{theorem}\label{uniquenessalg}
  If both the tuples $(\bQ_0, \Sigma_0, \tau_0, \Delta_0)$ and
  $(\bQ_1, \Sigma_1, \tau_1, \Delta_1)$ witness the strict
  algebraizability of $\Vdash$, then:
  \begin{enumerate}[(1)]
  \item $\bQ_0= \bQ_1$;
  \item $\Sigma_0 \equiv_{\bQ_i} \Sigma_1 $ with $i\in \{0,1\}$;
  \item
    $\Delta_0(\varepsilon,\delta) \dashv\Vdash
    \Delta_1(\varepsilon,\delta) $, for all
    $\varepsilon\approx \delta\in \Eq$;
  \item $\tau_0(\phi)\equiv^c_{\bQ_i}\tau_1(\phi) $ with
    $i\in \{0,1\}$, for all $\phi\in\Fm$.
  \end{enumerate}
\end{theorem}
\begin{proof}
  Clauses (1), (3) and (4) follow immediately from
  \cref{uniquenessalg_1}, so we prove (2). We write $\core_0(\cA)$ and $\core_1(\cA)$ for the core predicates induced by $\Sigma_0$ and $\Sigma_1$ on $A$, respectively. First, notice that by (1)
  we have $\bQ_0 = \bQ_1$, so let $\bQ = \bQ_i$, $i \in \{0,1\}$. Let
  $\alpha_0\approx\beta_0\in \Sigma_0$, then
  $\models^{0}_{\bQ} \alpha_0\approx\beta_0$ hence by
  \ref{loose_Alg}(\ref{l.alg2}) we obtain
  $\varnothing \Vdash \Delta_0(\alpha_0,\beta_0)$ and thus
  $\varnothing \Vdash \Delta_1(\alpha_0,\beta_0)$. It follows
  that $\models_{\bQ}^{1}\alpha_0\approx\beta_0$, meaning that
  $(\cA,\core_1(\cA)) \models \alpha_0\approx\beta_0$. Since this
  holds for all equations $\alpha_0\approx\beta_0\in \Sigma_0$, it
  follows that $\core_1(\cA) \subseteq \core_0(\cA)$. The other
  direction is proven analogously.
\end{proof}

\noindent The following corollary follows immediately from
\cref{facts.weak.algebraizability}.

\begin{corollary}\label{facts.strict.algebraizability} Let $\Vdash$ be
  a standard logic, then:
  \begin{enumerate}[(1)]
  \item if Conditions~\ref{loose_Alg}(\ref{l.alg1}) and
    \ref{loose_Alg}(\ref{l.alg4}), or
    Conditions~\ref{loose_Alg}(\ref{l.alg2}) and
    \ref{loose_Alg}(\ref{l.alg3}) hold, then $\Vdash$ is algebraized
    by $(\bQ, \Sigma, \tau, \Delta)$;
  \item if $\Vdash$ is algebraizable, then there are two finitary
    transformers $\tau$ and $ \Delta$ witnessing this fact.
  \end{enumerate}
\end{corollary}

\subsection{Fixed-Point Algebraizability of Weak Logics}\label{sec.2.4}

We introduce a third alternative definition of algebraizability for
weak logics, which can be seen as a refinement of strict
algebraizability. We first define the following notion of selector
term, essentially from \cite{nicolau2024polyatomic}.

\begin{notation}
  Let $t(x)$ be unary $\cL$-term, then we define recursively
  $t^1(x)=t(x)$ and $t^{n+1}(x)=t(t^n(x))$ for all $n<\omega$.
\end{notation}

\begin{definition}
  Let $\bQ$ be a quasivariety of $\cL$-algebras and $\delta(x)$ be a
  unary term in $\cL$, then we say that $\delta(x)$ is a
  \emph{selector term for $\bQ$} if
  $\bQ\models \delta(x) \approx \delta^2(x)$.
\end{definition}

\begin{remark}
  Selector terms essentially identifies the fixed points of
  polynomials in $\bQ$. If there is some $n<\omega$ such that for all
  $m<\omega$ we have $\delta^{n+m}(x)=\delta^{n}(x)$, then in
  particular we obtain
  \begin{align*}
    \delta^{n}(\delta^{n}(a)) =  \delta^{n+n}(a)= \delta^{n}(a)
  \end{align*}
  for all $A\in \bQ$. Thus $\delta^{n}$ is a selector term for $\bQ$.
\end{remark}

\begin{definition}\label{inv_Alg}
  A weak logic $\Vdash$ is \emph{fixed-point algebraizable} if it is
  strictly algebraized by a tuple $(\bQ,\Sigma,\tau,\Delta)$ (as in
  \ref{strict.Alg}) and, additionally,
  $\Sigma=\{ x\approx \delta(x) \}$ for some selector term
  $\delta(x)$.
\end{definition}

The key idea in the previous definition is that the core of the
expanded algebras corresponding to the logic $\Vdash$ is not simply
defined by a finite sets of equations, but by one single equation
characterising it as the set of fixed points of a certain
polynomial. The main motivation lies in the fact that fixed-point
algebraizability behaves very well in the context of varieties, as the
following result shows (this is essentially \cite[Thm.~19]{grillq} and
\cite[Thm.~3.14]{nicolau2024polyatomic}).

\begin{proposition}\label{projection_commutes_2}
  Let $\bV$ be a variety of $\cL$-algebras and
  $\Sigma(x)=\{\delta(x)\approx x\}$ for some selector term $\delta$,
  then the class of structures $(A,\core(\cA))$ with $A\in \bV$ and
  $ \core(\cA)=\Sigma(\cA) $ is a variety of expanded algebras.
\end{proposition}
\begin{proof}
  Let $\bV$ be a variety of $\cL$-algebras and let $\bK$ be the class
  of structures $(A,\core(\cA))$ with $A\in \bV$ and
  $ \core(\cA)=\Sigma(\cA) $. By Proposition~\ref{projection_commutes}
  it follows that $\bK$ is a quasivariety. We show that $\bK$ is also
  closed under homomorphic images. Let $A\in \bV$ and consider a
  surjective homomorphism $h:A\to B$, it suffices to show that
  $\core(\cB)=h[\core(\cA)]$. Since $\core(\cA)=\Sigma(\cA)$ and
  $\core(\cB)=\Sigma(\cB)$, it follows from the fact that homomorphic
  images preserve the validity of equations that
  $h[\core(\cA)]\subseteq \core(\cB)$. Consider now some
  $b\in \core(B)$ and notice that since $h$ is surjective there is
  some $a\in A$ with $h(a)=b$. Then we have that
  $h(\delta(a))=\delta(h(a))=h(a)$, given that
  $h(a)=b\in\core(\cB)=\Sigma(\cB)$. Since $\delta$ is a selector
  term, we have $A\models\delta(a)\approx \delta^2(a)$ and so
  $\delta(a) \in \core(\cA)$ and $b \in h[\core(\cA)]$.
\end{proof}

\begin{remark}
  We notice that most of our examples of algebraizable weak logics are
  actually fixed-point algebraizable. However, in the present work we
  shall not focus on the special features of fixed-point
  algebraizability, as our interest is rather in the more general
  properties of loose and strict algebraizability. We refer the reader
  to \cite{nicolau2024polyatomic} for an in-depth study of logics
  arising from selector terms.
\end{remark}

\section{Schematic Fragment and Standard Companions}\label{sec.3.1}

We start in this section to investigate the relation between the loose
and the strict version of algebraizability. In particular, we
introduce the notion of schematic fragment of a weak logic and relate
the algebraizability of a weak logic to the standard algebraizability
of its schematic fragment.

\subsection{The Schematic Fragment of a Weak Logic}

Given a weak logic $\Vdash$, it is natural to ask if we can associate
to it to some specific standard logical system. The study of negative
variants of intermediate logics led to the notion of schematic
fragment, which was originally introduced in
\cite[p.~545]{Miglioli1989-PIESRO} (under the name of
\emph{standardization}) and further investigated in
\cite{Ciardelli2011-CIAIL, Quadrellaro.2019B}. Here we generalise it
to arbitrary weak logics.

\begin{definition}[Schematic Fragment]\label{def:schematic_fragment}
  Let $\Vdash$ be a weak logic, we define its \emph{schematic
    fragment} $\Schm(\Vdash)$ as follows:
  \[\Schm(\Vdash) := \{ (\Gamma,\phi) \,: \,\forall\sigma\in
    \subst(\cL),\, \sigma[\Gamma]\Vdash\sigma(\phi) \}\] and we then
  also write $\Gamma\Vdash_{\mrm{s}}\phi$ if
  $(\Gamma,\phi)\in \Schm(\Vdash)$.
\end{definition}

\noindent It is clear from the definition that $\Schm(\Vdash)$ is the
largest standard logic contained in $\Vdash$. We use the schematic
fragment of a weak logic to relate loose and strict algebraizability
to standard algebraizability. To this end, we first introduce the
following notion of (finite) representability of a weak logic.

\begin{notation}
  Let $\Gamma\subseteq \Fm$ be a set of formulas, we let
  $\atomsubst[\Gamma]$ be the closure of $\Gamma $ under all atomic
  substitutions $\sigma\in\atomsubst(\cL)$. If $\Theta$ is a set of
  equations we similarly denote by $\atomsubst[\Theta]$ the closure of
  $\Theta $ under atomic substitutions.
\end{notation}

\begin{definition}\label{def:representability}
  We say that a weak logic $\Vdash$ is \emph{representable} if there
  is a set of formulas $\Lambda$ in one variable such that for all
  $\Gamma\cup \{\phi\} \subseteq \Fm$:
  \[ \Gamma \Vdash\phi \; \Longleftrightarrow \; \Gamma \cup
    \atomsubst[\Lambda] \Vdash_{\mrm{s}} \phi. \]
  \noindent We say that $\Vdash$ is \emph{finitely representable} if
  the condition above holds for some finite set
  $\Lambda\subseteq \Fm$.
\end{definition}

\begin{proposition}\label{characterisationalg_loose_1}
  If a weak logic $\Vdash$ is loosely algebraized by
  $(\bQ,\tau,\Delta)$, then its schematic fragment $\Schm(\Vdash)$ is
  algebraized by $(\bQ,\tau,\Delta)$ as well.
\end{proposition}
\begin{proof}
  Let $(\bQ,\tau,\Delta)$ witness the loose algebraizability of
  $\Vdash$. In particular $\bQ$ is a core-generated quasivariety of
  expanded algebra, and thus the collection of its algebraic reducts
  $\bQ\restriction\cL$ is a quasivariety of algebras. We claim that
  $(\bQ\restriction\cL,\tau,\Delta)$ witnesses the algebraizability of
  $\Schm(\Vdash)$. Notice that by Fact~\ref{facts.algebraizability} it
  suffices to verify that $(\bQ\restriction\cL, \tau, \Delta)$
  satisfies \ref{standard.Alg}(\ref{s.alg1}) and
  \ref{standard.Alg}(\ref{s.alg4}).  \smallskip \newline Firstly, we
  consider \ref{standard.Alg}(\ref{s.alg1}):
  \begin{flalign*}
	\Gamma \Vdash_{\mrm{s}} \phi &\Longleftrightarrow \forall \sigma \in \subst(\cL),\; \sigma[\Gamma] \Vdash \sigma(\phi) & \text{(by definition)} \\
	&\Longleftrightarrow  \forall \sigma \in \subst(\cL),\;\tau(\sigma[\Gamma]) \corevDash_{\bQ\restriction\cL} \tau(\sigma(\phi)) & \text{(by \ref{loose_Alg}(\ref{l.alg2}))} \\
	&\Longleftrightarrow  \forall \sigma \in \subst(\cL),\;\sigma(\tau[\Gamma]) \corevDash_{\bQ\restriction\cL} \sigma(\tau(\phi)) &  \text{(by structurality of $\tau$)} \\
	&\Longleftrightarrow \tau[\Gamma] \models_{\bQ\restriction\cL} \tau(\phi) &   \text{(by Proposition~\ref{th:subst-core}).}
  \end{flalign*}
  \smallskip \noindent Next, consider
  \ref{standard.Alg}(\ref{s.alg4}). To this end, suppose that
  $\tau(\Delta(x,y)) \not \equiv_{\bQ\restriction\cL} x \approx y$,
  then by Proposition~\ref{th:subst-core} and the fact that $\bQ$ is
  core-generated, we get that there exists a substitution $\sigma$
  such that
  $\sigma(\tau(\Delta(x,y))) \not\equiv^c_{\bQ} \sigma(x \approx
  y)$. Then, we obtain by the structurality of $\tau$ and $\Delta$
  that
  $\tau(\Delta(\sigma(x), \sigma(y))) \not\equiv^c_{\bQ} \sigma(x)
  \approx \sigma(y)$, contradicting the algebraizability of
  $\Vdash$. This completes the proof.
\end{proof}

The following theorem allows us to relate loose and strict
algebraizability of weak logics to the standard algebraizability of
their schematic fragment.

\begin{theorem}\label{characterisationalg_loose}
  For a weak logic $\Vdash$, the following are equivalent:
  \begin{enumerate}[(1)]
  \item $\Vdash$ is strictly algebraizable;
  \item $\Vdash$ is loosely algebraizable and $\Vdash$ is finitely
    representable;
  \item $\Schm(\Vdash)$ is algebraizable and $\Vdash$ is finitely
    representable.
  \end{enumerate}
\end{theorem}
\begin{proof}
  We first show that (1) entails (2). By definition, if $\Vdash$ is
  strictly algebraizable, then it is also loosely algebraizable. We
  show that it is also finitely representable. Let
  $(\bQ,\Sigma,\tau,\Delta)$ be the strict algebraic semantics of
  $\Vdash$. In particular we have that $\core(A)=\Sigma(\cA)$ for all
  $A\in \bQ$. We then have the following equivalences:
  \begin{flalign*}
	\Gamma \Vdash\phi & \Longleftrightarrow \tau[\Gamma] \models^c_{\bQ} \tau(\phi) & \text{(by \ref{loose_Alg}(\ref{l.alg1}))} \\
	& \Longleftrightarrow \tau[\Gamma] \cup \atomsubst[\Sigma] \models_{\bQ} \tau(\phi) & \text{(by Lemma~\ref{semantic:translation})} \\
	& \Longleftrightarrow \forall \sigma \in \subst(\cL),\; \sigma(\tau(\Gamma)) \cup \sigma(\atomsubst[\Sigma]) \models^c_{\bQ} \sigma(\tau(\phi))& \text{(by Proposition~\ref{th:subst-core})} \\
	& \Longleftrightarrow  \forall \sigma \in \subst(\cL),\; \tau(\sigma(\Gamma)) \cup \sigma(\atomsubst(\Sigma)) \models^c_{\bQ} \tau(\sigma(\phi))& \text{(by structurality of $\tau$)} \\
	& \Longleftrightarrow  \forall \sigma \in \subst(\cL),\; \sigma[\Gamma] \cup \Delta(\sigma(\atomsubst(\Sigma))) \Vdash \sigma(\phi)  & \text{(by \ref{loose_Alg}(\ref{l.alg2}),  \ref{loose_Alg}(\ref{l.alg3}))}	\\
	& \Longleftrightarrow  \forall \sigma \in \subst(\cL),\; \sigma[\Gamma] \cup \sigma(\atomsubst(\Delta(\Sigma))) \Vdash \sigma(\phi)  & \text{(by structurality of $\Delta$)}	\\
	& \Longleftrightarrow \Gamma \cup \atomsubst(\Delta(\Sigma)) \Vdash_{\mrm{s}} \phi,  & \text{}
  \end{flalign*}
  and thus $\Lambda=\Delta(\Sigma)$ witnesses the fact that $\Vdash$
  is finitely representable.

\smallskip
\noindent The direction from (2) to (3) follows immediately by
Proposition~\ref{characterisationalg_loose_1}.

\smallskip
\noindent Finally, we show that (3) entails (1). Suppose that
$\Schm(\Vdash)$ is algebraized by $(\bQ, \tau, \Delta)$ and that
$\Vdash$ is finitely represented via a set of formulas $\Lambda$. Let
$\core(\cA)=\tau[\Lambda](\cA)$ for all $\cA\in \bQ$ and consider
$\mathbf{Q'} := \mathbb{Q}(\{\langle \core(\cA) \rangle : \cA \in
\bQ\})$, i.e., $\mathbf{Q'}$ is the quasivariety of expanded algebras
generated by the algebras of the form $\langle \core(\cA) \rangle$,
for $\cA \in \bQ$. We then derive the following equivalences:
\begin{flalign*}
\Gamma \Vdash \phi &\iff \Gamma \cup \atomsubst[\Lambda] \Vdash_{\mrm{s}} \phi & \text{(by assumption)}  \\
&\iff  \tau[\Gamma] \cup \tau[\atomsubst[\Lambda]] \models_\bQ \tau(\phi)& \text{(by \ref{standard.Alg}(\ref{s.alg1}))}  \\
&\iff  \tau[\Gamma] \models_\bQ^c \tau(\phi). & \text{(by Lemma~\ref{semantic:translation})} \\
&\iff  \tau[\Gamma] \models_{\bQ'}^c \tau(\phi) & \text{(by Proposition~\ref{preservation.thm})}
\end{flalign*}
which establish \ref{loose_Alg}(\ref{l.alg1}).  We next verify that
\ref{loose_Alg}(\ref{l.alg4}) holds as well. Suppose that
$\tau(\Delta(x,y)) \not \equiv^c_\mathbf{Q'} x \approx y$, then since
$\core(\cA)=\tau[\Lambda](\cA)$ we obtain that
\[ \tau(\Delta(x,y) ) \cup \tau[\Lambda] \not \equiv_\mathbf{Q'} \{ x \approx y\} \cup \tau[\Lambda]  \]
which contradicts the fact that $\Schm(\Vdash)$ is algebraized by $(\bQ, \tau, \Delta)$.
\end{proof}

\noindent The following corollary is an immediate consequence of the
proofs of the previous theorem and lemma.

\begin{corollary}\label{characterisationalg_loose_5}
  A weak logic $\Vdash$ is strictly algebraized by
  $(\bQ,\Sigma,\tau,\Delta)$ if and only if it is represented by
  $\Delta[\Sigma]$ and $\Schm(\Vdash)$ is algebraized by
  $(\bQ,\tau,\Delta)$.
\end{corollary}

\subsection{The Lattice of Standard Companions}

The previous results established that, if $\Vdash$ is strictly
algebraizable, then there is a finite set $\Lambda$ witnessing its
finite representability. Essentially, this means that the consequences
in $\Vdash$ can be encoded as logical consequences in its schematic
fragment $\Vdash_{\mrm{s}}$, modulo the set of formulas
$\Lambda$. However, the schematic fragment of $\Vdash$ is not
necessarily the only standard logic which bears this property. We
study what are the standard logics that, up to $\Lambda$, share the
same consequences.

\begin{definition}\label{def:standard_companion_logic}
  Let $\Vdash$ be a strictly algebraizable weak logic which is
  finitely representable by $\Lambda$ and with strict algebraic semantics
  $(\bQ,\Sigma,\tau,\Delta)$. A standard logic $\vdash$ is a
  \emph{standard companion} of $\Vdash$ if the following conditions
  hold:
\begin{enumerate}[(1)]
\item
  $\Gamma \Vdash\phi \; \Longleftrightarrow \; \Gamma \cup
  \atomsubst[\Lambda] \vdash \phi $;
\item $\vdash$ is algebraized by a quasivariety $\bK$
  together with the transformers $\tau$ and $\Delta$.
\end{enumerate}
We denote by $\mrm{St}(\Vdash)$ the family of all standard
companions of $\Vdash$.
\end{definition}

\begin{remark}
  We provide some explanations of the previous
  definition. Condition~\ref{def:standard_companion_logic}(1) is
  essentially the same condition of representability from
  Definition~\ref{def:representability}, but with respect to an
  arbitrary standard logic. Thus
  Condition~\ref{def:standard_companion_logic}(1) identifies those
  standard logics that, up to $\Lambda$, deliver the same weak logic
  $\Vdash$. However, this condition by itself is quite weak, and thus
  we focus on logics that satisfy also condition
  \ref{def:standard_companion_logic}(2), i.e., that can be algebraized
  \emph{via} the same transformers $\tau$ and $\Delta$ as
  $\Vdash$. Notice that, in the specific cases of negative variants
  and polyatomic logics, these families had already been identified
  and studied in \cite{Quadrellaro.2019B} and
  \cite{nicolau2024polyatomic}, respectively.
\end{remark}

Our underlying intuition is that the family of standard companions
defined above must give rise to a corresponding notion on the side of
quasivarieties of algebras. 

\begin{definition}\label{def:standard_companion_algebra}
  Let $\bQ$ be a core-generated quasivariety of expanded algebras with
  core defined by a finite set of equations $\Sigma(x)$. Let $\bK$ be
  an arbitrary quasivariety of algebras (in the same signature as
  $\bQ$), then we say that $\bK$ is a \emph{standard companion} of
  $\bQ$ if, for all $A\in \bK$, we have
  $ \langle \Sigma(A) \rangle\in \bQ$.
\end{definition}


\begin{lemma}\label{companion_lemma}
  Let $\Vdash$ be a weak logic strictly algebraized by
  $(\bQ,\Sigma,\tau,\Delta)$, and let $\vdash$ be a standard
  logic. Then $\vdash\in \mrm{St}(\Vdash)$ if and only if $\vdash$ is
  algebraized by $(\bK,\tau,\Delta)$ for some $\bK$ standard companion
  of $\bQ$.
\end{lemma}
\begin{proof}
  We first prove the left-to-right direction.  Since
  $\vdash\in \mrm{St}(\Vdash)$ we have by definition that $\vdash$ is
  algebraized by some $(\bK,\tau,\Delta)$, where $\tau$ and $\Delta$
  also witness the strict algebraizability of $\Vdash$. In particular,
  notice that by the strict algebraizability of $\Vdash$ and
  Corollary~\ref{characterisationalg_loose_5} we have that $\Vdash$ is
  represented by $\Delta[\Sigma]$ over its schematic fragment, i.e.,
  that $\Lambda=\Delta[\Sigma]$. Now, we consider $\bK$ as a class of
  expanded algebras with core $\Sigma(\cA)$ for all $\cA\in \bK$ and
  we let $\bC=\mathbb{Q}(\{ \langle \Sigma(A)\rangle : A\in \bK
  \})$. Clearly, we have by Proposition~\ref{thm:maltsev} that
  $\Theta \models^c_{\bC}\varepsilon\approx \delta $ holds if and only
  if $\Theta \models^c_{\bK}\varepsilon\approx \delta$ does. Also, we
  have by strict algebraizability of $\Vdash$ and $\vdash$ that
  $\Theta \models^c_{\bK}\varepsilon\approx \delta$ is equivalent to
  $\Theta \models^c_{\bQ}\varepsilon\approx \delta$, since
  $\vdash\in \mrm{St}(\Vdash)$ and $\Lambda=\Delta[\Sigma]$. Thus we
  obtain that $\coretheories{\bC}=\coretheories{\bQ}$ and so by
  Proposition~\ref{thm:maltsev} $\bC=\bQ$, proving that $\bK$ is a
  standard companion of $\bQ$.

\smallskip
\noindent We consider the right-to-left direction. Suppose $\vdash$ is
algebraized by $(\bK,\tau,\Delta)$ and let $\bK$ be a standard
companion of $\bQ$, we show that $\vdash$ is a standard companion of
$\Vdash$. We view $\bK$ also as a quasivariety of expanded algebras by letting,  $\core(A)=\Sigma(A)$ for each $A\in \bK$. We have the following equivalences:
\begin{flalign*}
\Gamma \Vdash \phi & \Longleftrightarrow  \tau(\Gamma) \models^c_{\bQ} \tau(\phi) & \text{(by \ref{loose_Alg}(\ref{l.alg1}))} \\
& \Longleftrightarrow  \tau(\Gamma) \models^c_{\bK} \tau(\phi) & \text{(by Definition~\ref{def:standard_companion_algebra} and Proposition~\ref{thm:maltsev})} \\
& \Longleftrightarrow  \tau(\Gamma) \cup \atomsubst[\Sigma] \models_{\bK} \tau(\phi) & \text{(by Lemma~\ref{semantic:translation})} \\
& \Longleftrightarrow  \Gamma \cup \Delta(\atomsubst[\Sigma]) \vdash \phi & \text{(by \ref{standard.Alg}(\ref{s.alg2}) and \ref{standard.Alg}(\ref{s.alg3}))} \\
& \Longleftrightarrow  \Gamma \cup \atomsubst(\Delta[\Sigma]) \vdash \phi & \text{(by structurality of $\Delta$)}
\end{flalign*}
and thus, since $\Lambda=\Delta[\Sigma]$ (as we argued in the previous
direction), it follows that $\vdash$ is a standard companion of
$\Vdash$. This completes our proof.
\end{proof}

The correspondence provided by the previous lemma motivates the
following definition. Intuitively, while the schematic fragment
identifies the greatest standard logic contained in a weak logic, the
following notion characterises the greatest quasivariety with core
defined by $\Sigma$ and which has the same core-generated
substructures of a quasivariety $\bQ$. The following definition
essentially refines the notion of core superalgebra
(Definition~\ref{def:core_superstructure}) in the setting with equationally
definable core.

\begin{definition}\label{def:varopcore}
  We say that $\cB$ is a \emph{$\Sigma$-superalgebra} of $\cA$ if
  $\cA \leq \cB$ and $\Sigma(\cA) = \Sigma(\cB)$. If $\bQ$ is a class
  of algebras, then we write $\mathbb{C}_{\Sigma}(\bQ) $ for the class
  of all $\Sigma$-superalgebras of elements of $\bQ$.
\end{definition}

\begin{lemma}\label{axiomatisation_varopcore_1}
  Suppose $\bK$ is a quasivariety of expanded algebras with core
  defined by $\Sigma$, then $\varopCore(\bK)$ is also a quasivariety
  of expanded algebras.
\end{lemma}
\begin{proof}
  We notice that $A\in \varopCore(\bK)$ if and only if
  $A\models \forall x(\core(x)\leftrightarrow\bigwedge\Sigma(x))$,
  which can be rephrased by means of Horn formulas as in
  \ref{projection_commutes}. Also, for all quasiequations
  $\bigwedge_{i\leq n} \varepsilon_i \approx \delta_i \to
  \varepsilon\approx \delta \in \coretheories{\bK}$ we have that this
  is equivalent to saying
  \[ A\models \forall x_1,\dots \forall x_n ( \bigwedge_{1\leq i\leq
      n}\core(x_i) \land \bigwedge_{i\leq n} \varepsilon_i \approx
    \delta_i \to \varepsilon\approx \delta ), \] where these are also
  Horn formulas. It follows immediately by Fact~\ref{birkhoff_maltsev}
  that $\varopCore(\bK)$ is closed under $ \mathbb{I} $,
  $ \mathbb{S} $, $ \mathbb{P} $, $ \mathbb{P}_{\mrm{U}} $.
\end{proof}

We conclude this section with the following characterisation of the
family of standard companions of a strictly algebraizable weak logic
$\Vdash$. This generalises the previous results from
\cite{Quadrellaro.2019B,nicolau2024polyatomic} to our abstract
setting.

\begin{notation}
  Let $\tau:\Fm\to\Eq$ and let $\bQ$ be a quasivariety, then we write
  $\logic_\tau(\bQ)$ for the set of all pairs $(\Gamma,\phi)$ such
  that $\tau(\Gamma)\models_{\bQ} \tau(\phi)$.
\end{notation}

\begin{theorem}
  Let $\Vdash$ be a weak logic algebraized by
  $(\bQ,\Sigma,\tau,\Delta)$, then $\mrm{St}(\Vdash)$ with the subset
  ordering $\subseteq$ forms a bounded lattice with maximum element
  $\Schm(\Vdash)$ and minimum element $\logic_\tau(\varopCore(\bQ))$.
\end{theorem}
\begin{proof}
  Let $\vdash_0,\vdash_1 \in \mrm{St}(\Vdash)$, then
  $\vdash_0\cap \vdash_1$ is a standard logic and it is
  straightforward to verify that
  $\vdash_0\cap \vdash_1\in \mrm{St}(\Vdash)$. Similarly, if
  $\vdash_0,\vdash_1 \in \mrm{St}(\Vdash)$ and they are algebraized
  respectively by $(\bK_0,\tau,\Delta)$ and $(\bK_1,\tau,\Delta)$,
  then $\bK_0\cap \bK_1$ is a quasivariety, and the logic $\vdash_2$
  defined by letting
  \[ \Gamma\vdash_2 \phi \;\Longleftrightarrow\;
    \tau(\Gamma)\models^c_{\bK_0\cap \bK_1} \tau(\phi) \] is the
  supremum of $\vdash_0$ and $\vdash_1$ in $\mrm{St}(\Vdash)$. This
  shows that $(\mrm{St}(\Vdash), \subseteq)$ is a lattice.

  \medskip
  \noindent We show that $\Schm(\Vdash)$ is maximal in
  $\mrm{St}(\Vdash)$. Recall that, by
  Corollary~\ref{characterisationalg_loose_5}, we have that
  $\Schm(\Vdash)=\logic_\tau(\bQ)$. Suppose now that
  $\vdash\;\in \mrm{St}(\Vdash)$ then by monotonicity
  $\Gamma\vdash\phi$ entails $\Gamma\cup \Lambda\vdash \phi$, and thus
  we have $\Gamma\Vdash \phi$. This shows that
  $\vdash\;\subseteq \; \Vdash$. Since $\vdash$ and $\Schm(\Vdash)$
  are both closed under uniform substitution, it follows that
  $\vdash\subseteq \Schm(\Vdash)$.

  \medskip
  \noindent We show that $\logic_\tau(\varopCore(\bQ))$ is minimal in
  $\mrm{St}(\Vdash)$. If $\vdash\; \in \mrm{St}(\Vdash)$, then by
  Lemma~\ref{companion_lemma} $\vdash$ is algebraized by
  $(\bK,\tau,\Delta)$, for some standard companion $\bK$ of $\bQ$. By
  Definition~\ref{def:varopcore}, it follows that
  $\bK\subseteq \mathbb{C}_{\Sigma}(\bQ)$ and so
  $\logic_\tau(\varopCore(\bQ))\subseteq \logic_\tau(\bK)=\; \vdash$.
\end{proof}

\begin{remark}
  We notice that, if $\Vdash$ is fixed-point algebraized by
  $(\bQ,\Sigma,\tau,\Delta)$ with $\Sigma=\{x\approx \delta(x)\}$ as in \ref{inv_Alg}, then
  $\Theta\models^c_{\bQ} \alpha\approx \beta$ if and only if
  $\sigma_{\delta}(\Theta)\models \sigma_{\delta}(\alpha)\approx
  \sigma_{\delta}(\beta)$, where $\sigma_{\delta}$ is the substitution
  $\sigma_{\delta}(x)=\delta(x)$ for all $x\in \atoms$. Then, since
  $\Vdash_{\mrm{s}}$ is the logic $\logic_\tau(\bQ)$, we obtain that
  $\Gamma \Vdash \phi $ if and only if
  $\sigma_{\delta}(\Gamma)\Vdash_{\mrm{s}}
  \sigma_{\delta}(\phi)$. This correspondence showcases a special case
  where the equations identifying the core of the expanded algebra
  induce a translation between a logic and its schematic fragment. We
  refer the reader to \cite{nicolau2024polyatomic} for an in-depth
  study of such translations and to
  \cite{nicolau2024polyatomic,moraschini2018logical} for the relation
  between such translations and corresponding adjunctions between
  quasivarieties of algebras.
\end{remark}

\subsection{Example of Bridge Theorem: Deduction Theorem}

As stressed by Font in \cite[p. 160]{Font.2016}, ``Bridge theorems and
transfer theorems are the ultimate justification of abstract algebraic
logic''. We show that this motivation remains applicable in our
modified setting of weak logics. The key observation is that, in
virtue of \cref{characterisationalg_loose}, we do not have to come up
with novel algebraic descriptions of logical properties, but we can
rely on the established characterisation of such properties for
standard logics.

\smallskip
\noindent To showcase a concrete example, we look at one bridge
theorem for strictly algebraizable weak logics, i.e., we consider the
case of strictly algebraizable logics with a deduction-detachment
theorem. We start by recalling the definition of this property and the
standard bridge theorem relating logics with the deduction-detachment
theorem and quasivarieties with so-called EDPRC property.

\begin{definition}
  A weak (or standard) logic $\Vdash$ has the \emph{deduction-detachment theorem
    (DDT)} if there is a finite set of formulas $I(x,y)$ such that,
  for every set of formulas $\Gamma$,
  \[ \Gamma \cup \{ \phi \}\Vdash \psi \Longleftrightarrow \Gamma
    \Vdash I(\phi,\psi). \]
\end{definition}

\begin{definition}
  A quasivariety $\bK$ has \emph{equationally definable principal
    relative congruences} (EDPRC) if there is a finite set of
  equations $\Theta(x, y, z, v)$ such that, for every $\cA\in \bK$ and
  $a,b,c,d\in \cA$:
  \begin{align*}
    \langle a,b\rangle \in  \mathrm{Cn}_\bK(c,d)  \Longleftrightarrow \cA\models \Theta(c, d, a, b),
  \end{align*}
  where $\mathrm{Cn}_\bK(c,d) $ is the $\bK$-congruence of $\cA$
  induced by the equation $c\approx d$, i.e., the smallest congruence $\theta$ of $A$ such that $(c,d)\in \theta$ and $A/\theta\in \bK$ (cf.~\ref{Q-congruence}).
\end{definition}

\noindent The following result, due to Blok and Pigozzi, is one
classical bridge theorem from abstract algebraic logic. We refer the
reader to \cite[3.85]{Font.2016} for a proof of this result.

\begin{fact}[Blok-Pigozzi]\label{thm:blok_pigozzi_DDT}
  Let $\vdash$ be a standard algebraizable logic with equivalent
  algebraic semantics $\bQ$, then $\vdash$ has a deduction-detachment
  theorem if and only if $\bQ$ has EDPRC.
\end{fact}

In the case of a strictly algebraizable weak logic $\Vdash$, finite
representability allows us to transfer the property of having the DDT
from $\Vdash$ to $\Vdash_{\mrm{s}}$. More precisely, we can prove the
following proposition.

\begin{proposition}\label{prop:strict_DDT}
  Let $\Vdash$ be a strictly algebraizable weak logic, then $\Vdash$
  has DDT if and only of $\Vdash_{\mrm{s}}$ has DDT, witnessed by the same
  finite set of formulas $I(x,y)$.
\end{proposition}
\begin{proof}
  We first prove the right-to-left direction. Suppose
  $\Vdash_{\mrm{s}}$ has DDT, and notice that since $\Vdash$ is
  strictly algebraizable then it is finitely represented by a set of
  formulas $\Lambda$ (by \ref{characterisationalg_loose}). We have the
  following equivalences:
  \begin{flalign*}
	\Gamma\cup \{\phi\}\Vdash \psi &\Longleftrightarrow \Gamma\cup \{\phi\}\cup \atomsubst[\Lambda]\Vdash_{\mrm{s}} \psi & \text{(by representability)} \\
	&\Longleftrightarrow \Gamma\cup \atomsubst[\Lambda]\Vdash_{\mrm{s}} I(\phi,\psi) & \text{(by assumption)} \\
	&\Longleftrightarrow \Gamma \Vdash I(\phi,\psi) &  \text{(by representability)}
  \end{flalign*}
  which establish that $\Vdash$ has DDT.

  \medskip
  \noindent Consider now the left-to-right direction. We assume that
  $\Vdash$ has DDT and that this is witnessed by the finite set of
  formulas $I(x,y)$, then:
  \begin{flalign*}
\Gamma\cup \{\phi\}\Vdash_{\mrm{s}} \psi &\Longleftrightarrow \forall\sigma\in\subst(\cL), \; \sigma(\Gamma\cup \{\phi\})\Vdash \sigma(\psi)  & \text{(by Definition~\ref{def:schematic_fragment})} \\
&\Longleftrightarrow  \forall\sigma\in\subst(\cL), \; \sigma(\Gamma)\Vdash I(\sigma(\phi),\sigma(\psi))  & \text{(by assumption)} \\
&\Longleftrightarrow   \Gamma \Vdash_{\mrm{s}} I(\phi,\psi) & \text{(by Definition~\ref{def:schematic_fragment})}
  \end{flalign*}
  which proves our claim.
\end{proof}

\begin{corollary}
  Let $\Vdash$ be a strictly algebraizable weak logic with equivalent
  algebraic semantics $\bQ$, then $\Vdash$ has DDT if and only of
  $\bQ$ has EDPRC.
\end{corollary}
\begin{proof}
  Immediate from Fact~\ref{thm:blok_pigozzi_DDT} and
  Proposition~\ref{prop:strict_DDT}.
\end{proof}

\begin{remark}
  The previous result provides one example where one can obtain
  results about strictly algebraizable weak logics by investigating
  their schematic fragments. In particular, we did not need to come up
  with a novel algebraic characterisation for the presence of a
  deduction-detachment theorem in a strictly algebraizable weak logic
  $\Vdash$, and we could simply rely on
  \cref{characterisationalg_loose} and the classical bridge theorem by
  Blok and Pigozzi. This shows that strict algebraizability provides
  an especially strong and well-behaved notion in the setting of weak
  logics. We point out that we ascribe the underlying reason of this
  phenomena to the validity of a strong version of the Isomorphism
  Theorem in the case of strict algebraizability (in contrast to loose
  algebraizability). We shall investigate this in \cref{sec.3.2.3}.
\end{remark}

\section{The Isomorphism Theorem}\label{sec.3.2}

Blok and Pigozzi's isomorphism theorem \cite[\S 3.5]{Font.2016} is
possibly the most important single result on algebraizability, as it
relates the algebraizability of a logic to the existence, for every
algebra, of an isomorphism between the lattice of deductive filters of
the logic and the lattice of its relative congruences. We consider in
this section analogues of this result for the setting of algebraizable
weak logics. Firstly, in \cref{sec.3.2.1} we recall the isomorphism
theorem for standard algebraizable logic. In \cref{sec.3.2.2} we
provide a partial analogue of this result for loosely algebraizable
weak logics, and then in \cref{sec.3.2.3} we prove a stronger result
for the case of strictly algebraizable weak logics.

\subsection{Standard Algebraizable Logics}\label{sec.3.2.1}

We start by providing some preliminary definitions needed to state
Blok and Pigozzi's result. We first recall the notion of deductive
filter and fix some notation for algebraic congruences.

\begin{definition}
  For any algebra $\cA$ and standard logic $\vdash$, we say that
  $F\subseteq \cA$ is a \emph{deductive filter} of $\vdash$ on $\cA$
  if:
  \[ \Gamma\vdash \phi \implies \forall h\in \homset(\Fm, \cA), \;
    h[\Gamma]\subseteq F \text{ entails } h(\phi)\in F; \] and we let
  $\Fi_{\vdash}(\cA)$ be the set of all deductive filters of $\vdash$
  on $\cA$.
\end{definition}

\begin{remark}\label{remark:deductive_filters}
  We notice that $\Fi_{\vdash}(\cA)$ is closed under arbitrary
  intersections and thus that it forms a complete lattice under the
  subset ordering. Moreover, if $F$ is a deductive filter and $h$ is
  an endomorphism of $\cA$, then $h^{-1}(F)$ is also a deductive
  filter. We denote by $\Fi^+_{\vdash}(\cA)$ the lattice expansion
  \mbox{$(\Fi_{\vdash}, \subseteq, \{ h^{-1} : h \in \enset(\cL)
    \})$}. We refer the reader to \cite[\S2.3]{Font.2016} for proofs
  of these facts.
\end{remark}

\begin{notation}\label{Q-congruence}
  Given any algebra $\cA$ we let $\Con(\cA)$ be the set of all
  congruences $\theta$ over $\cA$, and $\Con_{\bQ}(\cA)$ be the set of
  all $\bQ$-congruences of $\cA$, i.e., those congruences $\theta$
  over $\cA$ such that $\cA/\theta\in\bQ$.
\end{notation}

\begin{remark}\label{remark:lattice_congruences}
  Similarly to the case of deductive filters, $\Con_{\bQ}(\cA)$ is
  closed under arbitrary intersections and therefore it is also a
  lattice under the subset ordering. Additionally, it is possible to
  verify that $\Con_{\bQ}(\cA)$ is closed under inverse endomorphisms
  of $\cA$. We then write $\Con^+_{\bQ}(\cA)$ for this lattice
  expansion, i.e.,
  $\Con^+_{\bQ}(\cA)=(\Con_{\bQ}(\cA), \subseteq, \{h^{-1} : h \in
  \enset(\cA)\})$.
\end{remark}

\noindent Finally, we introduce syntactical and semantical theories as
follows.

\begin{definition}
  If $\vdash$ is a standard logic, then we denote by $\Th(\vdash)$ the
  set of all \emph{(syntactic) theories} over $\vdash$, i.e., all sets
  $\Gamma \subseteq \Fm$ such that $\Gamma \vdash \phi$ entails
  $\phi\in \Gamma$. If $\bQ$ is a quasivariety, then
  $\Th(\models_{\bQ})$ denotes the set of \emph{(semantical) theories}
  over $\bQ$, i.e., the sets of equations $\Theta\subseteq \Eq$ such
  that $\Theta\models_{\bQ} \alpha\approx \beta$ entails
  $\alpha\approx \beta \in \Theta$.
\end{definition}

\begin{remark}
  We notice that $\Th(\vdash)$ forms a lattice under the subset
  relation, and that it is additionally closed under inverse
  substitutions. We refer by $\Th^+(\vdash)$ to the lattice expansion
  $(\Th(\vdash), \subseteq, \{ \sigma^{-1} : \sigma\in \subst(\cL) \}
  )$. Similarly, also $\Th(\models_{\bQ})$ forms a lattice under the
  subset ordering and is closed under inverse substitutions. We let
  $\Th^+(\models_{\bQ})=(\Th(\models_{\bQ}), \subseteq, \{ \sigma^{-1}
  : \sigma\in \subst(\cL) \} )$. These properties follow from the fact
  that the syntactic theories over $\vdash$ are exactly the deductive
  filters of $\vdash$ over $\Fm$, and that, under the identification
  of equations and pairs of formulas, the semantic theories over
  $\models_{\bQ}$ are the $\bQ$-congruences of $\Fm$.
\end{remark}

Blok and Pigozzi's isomorphism theorem for standard logics provides a
criterion to determine if a logic is algebraized by a quasivariety
based on their associated lattices of filters and congruences. For a
proof of the following theorem we refer the reader to \cite[\S
3.5]{Font.2016}.

\begin{fact}[Isomorphism Theorem]\label{standard.isomorphism.theorem}
  Let $\vdash$ be a standard logic and $\bQ$ a quasivariety, then the
  following are equivalent:
  \begin{enumerate}[(1)]
  \item $\vdash$ is algebraizable with equivalent algebraic semantics
    $\bQ$;
  \item $\Fi^+_{\vdash}(\cA) \cong \Con^+_{\bQ}(\cA)$, for any algebra
    $\cA$;
  \item $\Th^+(\vdash) \cong \Th^+(\models_{\bQ})$.
  \end{enumerate}
\end{fact}

\begin{remark}\label{remark:iso:alg}
  We clarify what are the underlying witnesses in the previous
  theorem (cf.~\cite[p.~150]{Font.2016}).  On the one hand, let $\vdash$ be an algebraizable logic
  with equivalent algebraic semantics $(\bQ,\tau,\Delta)$. Then the
  associated isomorphism $\Fi^+_{\vdash}(\cA) \cong \Con^+_{\bQ}(\cA)$
  is given by the following map from filters to congruences:
  \begin{align*}
    \theta_{(-)}:\Fi_{\vdash}(\cA) & \longrightarrow \Con_{\bQ}(\cA)\\
    G &\longmapsto \theta_G\coloneqq \{(a,b)\in A^2 :  \Delta^A(a,b)\subseteq G   \}
  \end{align*}
  and the following map from congruences to filters
  \begin{align*}
    F_{(-)}:\Con_{\bQ}(\cA) & \longrightarrow \Fi_{\vdash}(\cA) \\
    \eta &\longmapsto F_{\eta}\coloneqq \{ a\in A : \tau^A(a)  \subseteq \eta  \}.
  \end{align*}
  which can be shown to be inverse of each other.  On the other hand,
  suppose $\Omega:\Th^+(\vdash) \cong \Th^+(\models_{\bQ})$ is an
  isomorphism. Then the two transformers $\tau$ and $\Delta$ are
  defined as follows:
  \begin{align*}
    \tau(x)&=\sigma_x(\Omega(\Cn_{\vdash}(x)))\\
    \Delta(x,y)&= \sigma_{x,y}( \Omega^{-1}(\Cn_{\bQ}(x\approx y)   )  ).
  \end{align*}
  where $\Cn_{\vdash}$ and $\Cn_{\bQ}$ denote respectively the closure
  consequence operators on the logic $\vdash$ and the quasivariety
  $\bQ$, $\sigma_x$ is the substitution sending every variable to $x$,
  and $\sigma_{x,y}$ is the substitution sending every variable but
  $y$ to $x$ and $y$ to itself.
\end{remark}

\subsection{Loosely Algebraizable Weak Logics}\label{sec.3.2.2}

We prove in this section a (partial) version of the isomorphism
theorem for loosely algebraizable weak logics. To this end, we start
by introducing a version of deductive filters and congruences relative
to core semantics.

\begin{definition}
  For any expanded algebra $\cA$ and weak logic $\Vdash$, we say that
  $F\subseteq \cA$ is a \emph{core filter} of $\Vdash$ over $\cA$ if:
  \[
    \Gamma\Vdash \phi \implies \forall h\in \homset^c(\Fm, \cA), \;
    h[\Gamma]\subseteq F \text{ entails } h(\phi)\in F;
  \]
  and we denote the set of core-filters of $\cA$ with respect to
  $\Vdash$ by $\Fi_{\Vdash}(\cA)$.
\end{definition}

\begin{remark}
  If $A$ is an expanded algebra and $\bQ$ a quasivariety of algebras
  in the same algebraic signature of $A$, then we write $\Con(\cA)$
  and $\Con_{\bQ}(\cA)$ also for the congruences and the
  $\bQ$-congruences of $A$, respectively. Notice that, if $A$ is an
  expanded algebra and $\theta$ a congruence of its algebraic reduct,
  then the structure $A/\theta$ is simply the quotient of the
  algebraic reduct of $A$ by $\theta$ with
  $\core(A/\theta)=\core(A)/\theta$. This corresponds to viewing
  $A/\theta$ as a strong homomorphic image of $A$
  (cf. Definition~\ref{maps:expanded_algebras} and
  Remark~\ref{remark:def.homomorphism}).
\end{remark}

\begin{definition} Given an expanded algebra $\cA$ and a quasivariety
  $\bQ$ of expanded algebras, a congruence $\theta\in \Con(\cA)$ is
  said to be a \emph{core $\bQ$-congruence} if $\cA/\theta\in \bQ$. We
  write $\Con^c_{\bQ}(\cA)$ for the set of all core $\bQ$-congruences
  over $\cA$.
\end{definition}

\begin{remark}\label{remark:core_filters_and_congruences}
  Similarly to Remark~\ref{remark:deductive_filters} and
  Remark~\ref{remark:lattice_congruences}, it is possible to verify
  that both $\Fi_{\Vdash}(\cA)$ and $\Con^c_{\bQ}(\cA)$ are closed
  under arbitrary intersections, and thus form complete lattices under
  the subset ordering.
\end{remark}

\begin{lemma}\label{core_filters_are_filters}
  Let $\bQ$ be a quasivariety of expanded algebras and let $A$ be a
  core-generated expanded algebra, then:
  \begin{enumerate}[(1)]
  \item if $G\subseteq A$ is a core filter of $\Vdash$, then it is a
    deductive filter of its schematic fragment $\Vdash_{\mrm{s}}$;
  \item if $\theta$ is a core $\bQ$-congruence of $A$, then it is also a
    $\bQ\restriction\cL$-congruence.
  \end{enumerate}
\end{lemma}
\begin{proof}
  Clause (2) follows immediately from the definition of core
  $\bQ$-congruence. We thus consider clause (1). Let $G\subseteq A$ be
  a core filter of $\Vdash$, and let $\Gamma\Vdash_{\mrm{s}} \phi$,
  $h[\Gamma]\subseteq G$ for some $h\in\homset(\Fm,\cA)$. Now, since
  $\cA$ is core-generated, by
  Lemma~\ref{assignment_vs_core_assignment} we can find a substitution
  $\sigma$ and a core assignment $g$, such that $g(\sigma(x))=h(x)$
  for all $x\in \atoms$. In particular, it follows that
  $g[\sigma[\Gamma]]\subseteq G$. Now, since
  $\Gamma\Vdash_{\mrm{s}} \phi$, it follows by uniform substitution
  that $\sigma[\Gamma]\Vdash_{\mrm{s}} \sigma(\phi)$ and therefore
  $\sigma[\Gamma]\Vdash \sigma(\phi)$. Now, since $g$ is a core
  assignment and $G$ a core filter with
  $g[\sigma[\Gamma]]\subseteq G$, it follows that
  $g(\sigma(\phi))\in G$, whence $h(\phi)\in G$. This shows that $G$
  is a deductive filter with respect to $\Vdash_{\mrm{s}}$.
\end{proof}

\begin{lemma}\label{lemma_strong_endo}
  Let $\bQ$ be a quasivariety of expanded algebras, $A$ be a
  core-generated expanded algebra and let $h$ be a strong endomorphism
  of $A$, then the following hold:
  \begin{enumerate}[(1)]
  \item if $F$ is a core filter of $\Vdash$, then $h^{-1}(F)$ is also
    a core filter;
  \item if $\theta\in \Con^c_{\bQ}(\cA)$ then
    $h^{-1}(\theta)\in \Con^c_{\bQ}(\cA)$.
  \end{enumerate}
  Therefore, $\Fi_{\Vdash}(\cA)$ and $\Con^c_{\bQ}(\cA)$ are two
  lattices commuting with all inverse strong endomorphisms of $\cA$.
\end{lemma}
\begin{proof}
  We already pointed out in
  Remark~\ref{remark:core_filters_and_congruences} that
  $\Fi_{\Vdash}(\cA)$ and $\Con^c_{\bQ}(\cA)$ are (complete) lattices
  under the subset ordering, whence we simply prove (1) and (2).

  \smallskip
  \noindent We prove (1). Suppose $\Gamma\Vdash \phi$ and let
  $ g\in \homset^c(\Fm, \cA) $ be such that
  $g[\Gamma]\subseteq h^{-1}(F)$. Then $h(g[\Gamma])\subseteq F$ and,
  since $h\circ g\in \homset^c(\Fm, \cA)$, it follows that
  $h(g[\phi])\subseteq F$ and thus $g(\phi)\in h^{-1}(F)$.

  \smallskip
  \noindent We prove (2). The fact that $h^{-1}(\theta)\in \Con(\cA)$
  follows from Remark~\ref{remark:lattice_congruences}. We next show
  that $A/h^{-1}(\theta)\in \bQ$. Consider a quasiequation
  $\bigwedge_{i\leq n} \varepsilon_i \approx \delta_i \to \alpha
  \approx \beta \in \coretheories{\bQ}$. Since $A/\theta\in \bQ$ we
  have that
  $A/\theta\models^c \bigwedge_{i\leq n} \varepsilon_i \approx
  \delta_i \to \alpha \approx \beta$. Now, let
  $g\in\homset^c(\Fm, A) $ be such that
  $g(\varepsilon_i)/h^{-1}(\theta) = g(\delta_i)/h^{-1}(\theta)$ for
  all $i\leq n$, then it follows that
  $h\circ g(\varepsilon_i)/\theta = h\circ g(\delta_i)/\theta$ for all
  $i\leq n$. Since $h$ is a strong endomorphism of $A$ we obtain that
  $h\circ g(\alpha)/\theta = h\circ g(\beta)/\theta$. It then follows
  that $ g(\alpha)/h^{-1}(\theta) = g(\beta)/h^{-1}(\theta)$. This
  shows that
  $A/h^{-1}(\theta)\models^c \bigwedge_{i\leq n} \varepsilon_i \approx
  \delta_i \to \alpha \approx \beta$ and so
  $A/h^{-1}(\theta)\models^c \coretheories{\bQ}$. Finally, notice that
  since $A$ is core generated, also $A/h^{-1}(\theta)$ is core
  generated, and thus we have from Proposition~\ref{thm:maltsev} that
  $A/h^{-1}(\theta) \in \bQ$. We conclude that
  $h^{-1}(\theta)\in \Con^c_{\bQ}(\cA)$.
\end{proof}

As expanded algebras augment standard algebras with an additional core
predicate, we now expand $\Fi^+_{\vdash}(\cA) $ and
$ \Con^+_{\bQ}(\cA)$ in a similar way to capture core filters and
congruences.

\begin{definition}\label{expanded_lattice}
  Let $\Vdash$ be a weak logic and let $\bQ$ be a quasivariety of
  expanded algebras, for any core-generated expanded algebra $A$, we
  write $\overline{\Fi}_{\Vdash}(\cA) $ for the structure
  \[ (\Fi_{\Vdash_{\mrm{s}}}(\cA), \; \Fi_{\Vdash}(\cA), \;\subseteq,
    \; \{\sigma^{-1}: \sigma\in\strongend(\cA)\} \] and we write
  $\overline{\Con}_{\bQ}(\cA) $ for the structure
  \[ (\Con_{\bQ\restriction\cL}(\cA), \; \Con^c_{\bQ}(\cA), \;
    \subseteq, \; \{\sigma^{-1}: \sigma\in\strongend(\cA)\}), \] where
  $\strongend(A)$ is the set of all strong endomorphisms of the
  expanded algebra $A$ (recall
  Definition~\ref{maps:expanded_algebras}).
\end{definition}

\begin{remark} The key intuition in the context of weak logics is to
  consider \emph{together} arbitrary filters and core filters, and to
  consider \emph{together} arbitrary $\bQ$-congruences and core
  $\bQ$-congruences. This is possible because of
  Lemma~\ref{core_filters_are_filters}, which makes sure that core
  filters are deductive filters, and that core congruences are
  congruences. Then, by Lemma~\ref{lemma_strong_endo}, we additionally
  have that strict endomorphisms also preserve core filters and core
  congruences.
\end{remark}

\begin{theorem}\label{algebraizabilit_to_iso}
  Let $\Vdash$ be a weak logic and suppose it is loosely algebraized
  by $(\bQ,\tau,\Delta)$, then for every core-generated expanded
  algebra $\cA$ there is an isomorphism
  $\Omega:\overline{\Fi}_{\Vdash}(\cA) \cong
  \overline{\Con}_{\bQ}(\cA)$.
\end{theorem}
\begin{proof}
  Firstly, we define the following map from filters to congruences:
  \begin{align*}
    \theta_{(-)}:\Fi_{\Vdash_{\mrm{s}}}(\cA) & \longrightarrow \Con_{\bQ\restriction\cL}(\cA)\\
    G &\longmapsto \theta_G\coloneqq \{(a,b)\in A^2 :  \Delta^A(a,b)\subseteq G   \}
  \end{align*}
  and the following map from congruences to filters
  \begin{align*}
    F_{(-)}:\Con_{\bQ\restriction\cL}(\cA) & \longrightarrow \Fi_{\Vdash_{\mrm{s}}}(\cA) \\
    \eta &\longmapsto F_{\eta}\coloneqq \{ a\in A : \tau^A(a)  \subseteq \eta  \}.
  \end{align*}
  Now, notice that since $\Vdash$ is loosely algebraized by
  $(\bQ,\tau,\Delta)$, then it follows from
  Proposition~\ref{characterisationalg_loose_1} that also
  $\Vdash_{\mrm{s}}$ is algebraized by $(\bQ,\tau,\Delta)$. Since the
  two maps above are exactly the same that occur in the proof of the    isomorphism theorem for standard logics
  (cf. Remark~\ref{remark:iso:alg}), it follows immediately that they
  describe an isomorphism
  $\Fi^+_{\Vdash_{\mrm{s}}}(\cA) \cong \Con^+_{\bQ}(\cA)$. Thus, it
  suffices to show that $\theta_{(-)}$ sends core filters to core
  congruences, and that $F_{(-)}$ sends core congruences to core
  filters.

  \smallskip
  \noindent First, let $G$ be a core filter over $A$, we show that
  $\theta_G$ is a core-$\bQ$-congruence of $A$. Consider $A/\theta_G$
  and notice that this is well-defined because $\theta_G$ is a
  congruence. Now, consider the expansion of $A/\theta_G$ defined by
  letting $\core(A/\theta_G)=\core(A)/\theta_G$, we claim that
  $(A/\theta_G, \core(A)/\theta_G)\in \bQ$.  Let
  $\bigwedge_{i\leq n} \varepsilon_i \approx \delta_i \models^c_{\bQ}
  \alpha \approx \beta$ and suppose $h\in\homset^c(\Fm,A/\theta_G)$ is such
  that $h(\epsilon_i)=h(\delta_i)$ for all $i\leq n$. By
  algebraizability we have that
  $\bigcup_{i\leq n}\Delta(\varepsilon_i \approx \delta_i) \Vdash
  \Delta(\alpha,\beta)$ and so $h(\Delta(\alpha,\beta))\subseteq
  G$. It follows that $(h(\alpha),h(\beta))\in\theta_G$ and so
  $A/\theta_G\models^c \alpha\approx \beta$. This shows that
  $(A/\theta_G, \core(A)/\theta_G)\models^c \bigwedge_{i\leq n}
  \varepsilon_i \approx \delta_i \to \alpha\approx \beta$. Since $\bQ$
  is a core-generated quasivariety of expanded algebras, it follows
  from Proposition~\ref{thm:maltsev} that
  $(A/\theta_G, \core(A)/\theta_G)\in \bQ$. We conclude that the map
  $\theta_{(-)}$ sends core filters to core-$\bQ$-congruences.

  \smallskip
  \noindent Next, we show that for any core-$\bQ$-congruence $\theta$
  over $A$ the filter $F_\theta$ is a core filter. Suppose
  $\Gamma \Vdash \phi$ and $ h\in \homset^c(\Fm, \cA) $ is such that
  $h[\Gamma]\subseteq F_\theta$. Then, we have by definition that
  $\tau^A(c)\subseteq \theta$ for every $c\in h[\Gamma]$. By
  algebraizability, notice that we have
  $\tau(\Gamma)\models^c_{\bQ}\tau(\phi)$, thus we obtain
  $h(\tau(\phi))=\tau^A(h(\phi))\subseteq \theta$ and so
  $h(\phi)\in F_\theta$. This shows that $F_\theta$ is a core
  deductive filter and completes the proof.
\end{proof}

Exactly as in the standard setting, it is possible to specialise the
previous isomorphism to the specific case of filters and congruences
of the (expanded) algebra of formulas $\Fm$. We first define syntactic
and semantic theories in the setting of weak logics.

\begin{definition} \label{def:core_theories} If $\Vdash$ is a weak
  logic then we denote by $\Th(\Vdash)$ the set of all
  \emph{(syntactic) theories} over $\Vdash$, i.e., all sets
  $\Gamma \subseteq \Fm$ such that if $\Gamma\Vdash\phi$ then
  $\phi\in \Gamma$. If $\bQ$ is a core-generated quasivariety of
  expanded algebras, then $\Th(\corevDash_{\bQ})$ denotes the set of
  \emph{(semantical) core theories} over $\bQ$, i.e., those sets of
  equations $\Theta \subseteq \Eq$ such that
  $\Theta \corevDash_{\bQ} \alpha \approx \beta$ entails
  $\alpha \approx \beta \in \Theta$.
\end{definition}

\begin{remark}\label{remark:core_theories_filters}
  It is straightforward to verify that, for $\Vdash$ a weak logic, the
  syntactic theories from \ref{def:core_theories} are exactly the core
  filters of $\Vdash$ over $(\Fm,\atoms)$. Similarly, if $\bQ$ is a
  quasivariety of expanded algebras, then the semantical core theories
  over $\bQ$ are exactly the core congruences of $(\Fm,\atoms)$. Since
  $(\Fm,\atoms)$ is a core-generated algebra, we obtain as in
  Definition~\ref{expanded_lattice} an expanded lattice
  $\overline{\Th}(\Vdash) $ defined by
  \[ (\Th(\Vdash_{\mrm{s}}), \; \Th(\Vdash), \; \subseteq, \;
    \{\sigma^{-1}: \sigma\in\atomsubst(\cL)\}) \] and an expanded
  lattice $\overline{\Th}(\corevDash_{\bQ}) $ defined by
  \[ (\Th(\models_{\bQ}), \; \Th(\corevDash_{\bQ}), \; \subseteq, \;
    \{\sigma^{-1}: \sigma\in\atomsubst(\cL)\}).  \] This clearly
  corresponds to the fact that theories over $\Vdash$ are also
  theories in $\Vdash_{\mrm{s}}$, and theories over $\models^c_{\bQ}$
  are also theories over $\models_{\bQ}$. Additionally, since the
  strong endomorphisms of $(\Fm,\atoms)$ are exactly the atomic
  substitutions, it follows from Lemma~\ref{lemma_strong_endo} that
  inverse atomic substitutions preserve elements in $\Th(\Vdash)$ and
  $\Th(\corevDash_{\bQ})$.
\end{remark}

\begin{corollary}\label{loose_corollary}
  Let $\Vdash$ be a weak logic and suppose it is loosely algebraized
  by $(\bQ,\tau,\Delta)$, then there is an isomorphism
  $\Omega:\overline{\Th}(\Vdash) \cong
  \overline{\Th}(\corevDash_{\bQ})$.
\end{corollary}
\begin{proof}
  It follows immediately from \cref{algebraizabilit_to_iso} and
  Remark~\ref{remark:core_theories_filters}.
\end{proof}

\begin{remark}\label{remark:problems}
  A key aspect of Blok and Pigozzi's isomorphism theorem
  (Fact~\ref{standard.isomorphism.theorem}) is that one can recover
  the two transformers $\tau$ and $\Delta$ from the isomorphism
  between $\Fi^+_{\vdash}(\cA)$ and $\Con^+_{\bQ}(\cA)$, thus showing
  algebraizability. We indicate however that this fact is not as clear
  in the setting of loose algebraizability. In fact, in the proof of
  the isomorphism theorem (cf. Remark~\ref{remark:iso:alg} and \cite[3.5]{Font.2016}) the two
  transformers are defined as follows:
  \begin{align*}
    \tau(x)&=\sigma_x(\Omega(\Cn_{\vdash}(x)))\\
    \Delta(x,y)&= \sigma_{x,y}( \Omega^{-1}(\Cn_{\bQ}(x\approx y)   )  );
  \end{align*}
  where $\Cn_{\vdash}$ and $\Cn_{\bQ}$ denote respectively the closure
  consequence operators on the logic $\vdash$ and the quasivariety
  $\bQ$, while $\sigma_x$ is a substitution sending every variable to
  $x$, and $\sigma_{x,y}$ is one sending every variable but $y$ to $x$
  and $y$ to itself. If one simply mimics this proof in the setting of
  weak logics and lets, for example,
  \[ \tau(x)=\sigma_x(\Omega(\Cn_{\Vdash}(x))), \] where
  $\Cn_{\Vdash}$ is the closure consequence operator of the weak
  logic $\Vdash$, then crucially it could be that structurality
  fails. In fact, we could have that $x\Vdash \delta(x)$ but also
  $\phi\not\Vdash \delta(\phi)$, witnessing the failure of uniform
  substitution. Let $\sigma$ be the substitution sending every
  variable to $\phi$, then we have
  \begin{align*}
    \delta(\phi)\in \Cn_{\Vdash}(\sigma( \Cn_{\Vdash}(x))) \text{ and }\delta(\phi)\notin \Cn_{\Vdash}(\sigma(x)),
  \end{align*}
  and therefore we obtain
  $\Omega(\Cn_{\Vdash}(\sigma( \Cn_{\Vdash}(x))))\neq
  \Omega(\Cn_{\Vdash}(\sigma(x))).$ Clearly, we have that
  $\Omega(\Cn_{\Vdash}(\sigma(x)))= \Omega(\Cn_{\Vdash}(\phi))$,
  simply by definition of $\sigma$. Assume additionally that $\Omega$
  commutes with closure operators and substitution, i.e., that
  \[ \Omega(\Cn_{\Vdash}(\sigma( \Gamma))) = \Cn^c_{\bQ}(
    \sigma(\Omega (\Cn_{\Vdash}( \Gamma) ) ) ) \] for all sets of
  formulas $\Gamma$ and substitutions $\sigma$ (and where
  $\Cn^c_{\bQ}$ refers the closure operator associated to the core
  semantics of a quasivariety of core-generated expanded algebras
  $\bQ$). Then we obtain:
  \begin{align*}
	\Omega(\Cn_{\Vdash}(\sigma( \Cn_{\Vdash}(x)))) &= \Cn^c_{\bQ}( \sigma(\Omega (\Cn_{\Vdash}( \Cn_{\Vdash}(x))   )   )  ) \\
	&= \Cn^c_{\bQ}( \sigma(\Omega (\Cn_{\Vdash}(x))   )   )\\
	& = \Cn^c_{\bQ}( \sigma\sigma_x(\Omega (\Cn_{\Vdash}(x))   )   )\\
	& = \Cn^c_{\bQ}( \tau(\phi))
  \end{align*}
  and therefore
  $\Cn^c_{\bQ}( \tau(\phi))\neq \Omega(\Cn_{\Vdash}(\phi))$. This
  indicates that the usual proof of the isomorphism theorem does not
  work in the setting of loosely algebraizable weak logics. We leave
  it as a pointer for future works to establish if there is a version
  of the isomorphism theorem in the setting of loosely algebraizable
  weak logics.
\end{remark}

\subsection{Strictly Algebraizable Weak Logics}\label{sec.3.2.3}

In contrast to Remark~\ref{remark:problems}, in the case of strictly
algebraizable weak logics, we can use the fact that they are finitely
representable to derive a full version of the isomorphism theorem. We
first recall some preliminary definitions, which we employed
informally already in Remark~\ref{remark:problems}.

\begin{definition}[Closure
  Operators] \label{def:weak_closure_operators} Let $\cA$ be an
  expanded algebra, $\Vdash$ a weak logic, and $\bQ$ a core-generated
  quasivariety of expanded algebras. We denote by $ \Cn_{\Vdash} $ the
  closure operator on $\cA$ defined by letting $\Cn_{\Vdash}(X)$ be
  the smallest core filter of $A$ with respect to $\Vdash$ that
  contains $X\subseteq A$. Similarly, we denote by $\Cn^c_{\bQ}$ the
  closure operator on $A^2$ defined by letting $ \Cn^c_{\bQ} (R)$ be
  the smallest congruence in $\Con^c_{\bQ}(\cA)$ containing
  $R\subseteq A^2$.
\end{definition}

\begin{remark}\label{remark:theories_filters}
  Recall that, in the expanded term algebra $(\Fm,\atoms)$, the core
  filters coincide exactly with the syntactical theories and the core
  congruences coincide with the semantical theories. In this case we
  recover the more intuitive definition of the closure operators
  $\Cn_{\Vdash} $ and $\Cn^c_{\bQ}$, namely, for $\Theta\subseteq \Eq$
  and $\Gamma\subseteq \Fm$:
  \begin{align*}
    \Cn_{\Vdash}(\Gamma)&=\{\phi\in \Fm : \Gamma\Vdash \phi   \}\\
    \Cn^c_{\bQ} (\Theta)&=\{\alpha\approx\beta \in \Eq : \Theta\models^c_{\bQ} \alpha\approx\beta   \}.
  \end{align*}

\end{remark}

\noindent In the context of strictly algebraizable weak logics, it is
also convenient to work with the following notions of
$\Lambda$-filters and $\Sigma$-congruences.

\begin{definition}
  Let $\Lambda$ be a set of formulas in $\Fm$, then for any algebra
  $\cA$ and standard logic $\vdash$, we say that $F\subseteq \cA$ is a
  \emph{$\Lambda$-filter} over $\cA$ with respect to $\vdash$ if
  \[
    \Gamma\cup \atomsubst[\Lambda]\vdash \phi \implies \forall h\in
    \homset(\Fm, \cA), \; h[\Gamma]\subseteq F \text{ entails }
    h(\phi)\in F.
  \]
  We denote the set of $\Lambda$-filters of $\cA$ with respect to
  $\vdash$ by $\Lambda\Fi_{\vdash}(\cA)$. Similarly, we let
  $\Lambda\Th(\vdash) $ be the set of all syntactical theories
  $\Gamma$ over $\vdash$ such that
  $\atomsubst[\Lambda] \subseteq \Gamma$.
\end{definition}

\begin{definition}
  Given an expanded algebra $\cA$, a quasivariety of expanded algebras
  $\bQ$ and a finite set of equations in one variable $\Sigma$, a
  congruence $\theta\in \Con_{\bQ\restriction\cL}(\cA)$ is said to be
  a \emph{$\Sigma$-congruence} if $\Sigma(a)\subseteq \theta$ for all
  $a\in\core(\cA)$. We write $\Sigma\Con_{\bQ\restriction\cL}(\cA)$
  for the set of all $\bQ\restriction\cL$-congruences of $A$ which are
  also $\Sigma$-congruences. Similarly, we let
  $\Sigma\Th(\models_{\bQ})$ denote the the set of all semantical
  theories $\Theta$ such that $ \atomsubst[\Sigma]\subseteq \Theta$.
\end{definition}

\begin{remark}
  It is clear from the definition of $\Lambda$-filters and the
  monotonicity of $\vdash$ that the $\Lambda$-filters of $\vdash$ on
  $A$ are a special kind of deductive filters of $\vdash$. Similarly,
  $\Sigma$-congruences are clearly a special kind of congruences. The
  set of $\Lambda$-filters as well as that of $\Sigma$-congruences are
  closed under arbitrary intersections and hence form complete
  lattices under the subset relation, but, crucially, they are not
  necessarily closed under arbitrary inverse endomorphisms.
\end{remark}

\begin{lemma}\label{lambda_filter_lemma}
  Let $\Vdash$ be a weak logic. If $\Vdash$ is represented by a set
  of formulas $\Lambda$ then
  $\Fi_{\Vdash}(\cA)=\Lambda\mathsf{Fi}_{\Vdash_{\mrm{s}}}(\cA)$ for
  any core-generated expanded algebra $A$. Moreover, if $\Th(\Vdash)=\Lambda\Th(\Vdash_{\mrm{s}})$, then $\Vdash$ is represented by $\Lambda$.
\end{lemma}
\begin{proof}
  We prove the first claim. Suppose $\Vdash$ is
  represented by $\Lambda$. Let $F\in \Fi_{\Vdash}(\cA)$,
  $\Gamma\cup \atomsubst[\Lambda]\Vdash_{\mrm{s}}\phi$ and
  $h\in \homset(\Fm, \cA) $ be such that $h[\Gamma]\subseteq F$. By
  the representability of $\Vdash$ we have that
  $\Gamma\Vdash\phi$. Since $\cA$ is core-generated, by
  Lemma~\ref{assignment_vs_core_assignment} there are a
  core-assignment $g\in \homset^c(\Fm, \cA) $ and a substitution
  $\sigma\in\subst(\cL)$ such that
  $h[\Gamma\cup \{\phi \}]=g[\sigma(\Gamma\cup \{\phi\})]$. It follows
  that $g[\sigma(\Gamma)]\subseteq F$ and so, since
  $\Gamma\Vdash\phi$, we obtain from the definition of core filters
  that $h(\phi)=g(\sigma(\phi))\in F$. This shows that
  $\Fi_{\Vdash}(\cA)\subseteq\Lambda\mathsf{Fi}_{\Vdash_{\mrm{s}}}(\cA)$.

  \smallskip
  \noindent Conversely, let
  $F\in \Lambda\mathsf{Fi}_{\Vdash_{\mrm{s}}}(\cA)$. Suppose
  $\Gamma\Vdash\phi$ and let $h\in \homset^c(\Fm, \cA) $ be such that
  $h[\Gamma]\subseteq F$. By representability we obtain that
  $\Gamma\cup \atomsubst[\Lambda]\Vdash_{\mrm{s}} \phi $, thus it
  follows from $F\in \Lambda\mathsf{Fi}_{\Vdash_{\mrm{s}}}(\cA)$ that
  $h(\phi)\in F$.  This shows that
  $\Lambda\mathsf{Fi}_{\Vdash_{\mrm{s}}}(\cA)\subseteq
  \Fi_{\Vdash}(\cA)$ and thus proves that
  $\Fi_{\Vdash}(\cA)=\Lambda\mathsf{Fi}_{\Vdash_{\mrm{s}}}(\cA)$.

  \medskip
  \noindent We now prove the second claim. Suppose we have
  that $\Th(\Vdash)=\Lambda\Th(\Vdash_{\mrm{s}})$. Then,
  for any $\psi\in \Fm$ we have that:
  \begin{align*}
    \Cn_{\Vdash}(\psi) &= \Cn_{\Vdash_{\mrm{s}}}(\Cn_{\Vdash}(\psi) \cup \atomsubst[\Lambda]); \\
    \Cn_{\Vdash_{\mrm{s}}}(\psi\cup \atomsubst[\Lambda]) &= \Cn_{\Vdash}(\Cn_{\Vdash_{\mrm{s}}}(\psi \cup \atomsubst[\Lambda])).
  \end{align*}
  Therefore, it follows from the first equality (and monotonicity of
  the consequence operators) that
  $\Cn_{\Vdash_{\mrm{s}}}(\psi\cup \atomsubst[\Lambda])\subseteq
  \Cn_{\Vdash}(\psi)$, and it follows from the second equality that
  $\Cn_{\Vdash}(\psi) \subseteq \Cn_{\Vdash_{\mrm{s}}}(\psi\cup
  \atomsubst[\Lambda])$. This yields that
  $\Cn_{\Vdash}(\psi) = \Cn_{\Vdash_{\mrm{s}}}(\psi\cup
  \atomsubst[\Lambda])$. Finally, we derive the following equivalences:
  \begin{flalign*}
    \Gamma \Vdash\phi & \Longleftrightarrow \Cn_{\Vdash}(\phi)\subseteq \Cn_{\Vdash}(\Gamma)  \\
                      & \Longleftrightarrow \Cn_{\Vdash_{\mrm{s}}}(\phi\cup\atomsubst[\Lambda])\subseteq \Cn_{\Vdash_{\mrm{s}}}(\Gamma\cup\atomsubst[\Lambda]) \\
                      & \Longleftrightarrow \Gamma \cup \atomsubst[\Lambda] \Vdash_{\mrm{s}} \phi,
  \end{flalign*}
  \noindent which prove the representability of $\Vdash$ via
  $\Lambda$.
\end{proof}

\begin{lemma}\label{sigma_congruences_lemma}
  Let $\bQ$ be a core-generated quasivariety with core defined by a
  set of equations $\Sigma$, then $\core(\cA)=\Sigma(\cA)$ entails
  $\Con^c_{\bQ}(\cA)=\Sigma\Con_{\bQ\restriction\cL}(\cA)$ for any
  core-generated expanded algebra $A$.
\end{lemma}
\begin{proof}
  Suppose that $\core(\cA)=\Sigma(\cA)$ for all $A\in \bQ$. First, let
  $\theta\in \Con^c_{\bQ}(\cA) $, then it follows that
  $\cA/\theta\in \bQ$ and so that
  $\core(\cA/\theta)= \Sigma[\cA/\theta]$. This yields
  $\Sigma(a)\subseteq \theta$ for all $a\in \core(\cA)$, showing that
  $\theta\in \Sigma\Con_{\bQ\restriction\cL}(\cA) $. Conversely,
  suppose $\theta\in \Sigma\Con_{\bQ\restriction\cL}(\cA)$, then in particular
  $\cA/\theta\in \bQ\restriction \cL$. Since $\bQ$ has core defined by the
  set of equations $\Sigma$, the unique expansion of $\cA/\theta$ which belongs to $\bQ$ is $(\cA/\theta, \Sigma(\cA/\theta))\in \bQ$. It follows that $\core(\cA/\theta)=\Sigma(\cA/\theta)$ and so that $\theta \in \Con^c_{\bQ}(\cA)$.
\end{proof}

We use the previous lemmas to obtain the following isomorphism theorem
for strictly algebraizable weak logics. This result displays how the
setting of strict algebraizability better reproduces the original
approach from Blok and Pigozzi in the context of logics without
uniform substitutions.

\begin{theorem}\label{full.isomorphism_2}
  Let $\Vdash$ be a weak logic and $\bQ$ a core-generated quasivariety
  of expanded algebras with core defined by a finite set of equations
  $\Sigma$. The following statements are equivalent.
  \begin{enumerate}[(1)]
  \item $\Vdash$ is strictly algebraized by
    $(\bQ,\Sigma, \tau,\Delta)$.
  \item There is a finite set of formulas $\Lambda$ such that, for
    every core-generated expanded algebra $\cA$, it holds that
    $\Fi_{\Vdash}(\cA)=\Lambda\mathsf{Fi}_{\Vdash_{\mrm{s}}}(\cA)$ and
    $\Con^c_{\bQ}(\cA)=\Sigma\Con_{\bQ\restriction\cL}(\cA)$. Moreover,
    there is an isomorphism
    $\Omega:\overline{\Fi}_{\Vdash}(\cA) \cong
    \overline{\Con}_{\bQ}(\cA)$.
  \item There is a finite set of formulas $\Lambda$ such that
    $\Th(\Vdash)=\Lambda\Th(\Vdash_{\mrm{s}})$ and
    $\Th(\corevDash_{\bQ})=\Sigma\Th(\models_{\bQ})$. Moreover, there
    is an isomorphism
    $\Omega:\overline{\Th}(\Vdash) \cong
    \overline{\Th}(\corevDash_{\bQ})$.
	\end{enumerate}
\end{theorem}
\begin{proof}
  Direction from (1) to (2) follows immediately from
  \cref{algebraizabilit_to_iso} and Lemmas \ref{lambda_filter_lemma},
  \ref{sigma_congruences_lemma}. Direction from (2) to (3) follows
  from Remark~\ref{remark:core_theories_filters}. It remains to show
  that (3) entails (1).

\smallskip
\noindent Firstly, notice that the isomorphism
$\Omega:\overline{\Th}(\Vdash) \cong \overline{\Th}(\corevDash_{\bQ})$
induces an isomorphism
$\Th^+(\Vdash_{\mrm{s}}) \cong \Th^+(\models_{\bQ})$, thus by the
standard isomorphism theorem it immediately follows that
$\Vdash_{\mrm{s}}$ is algebraized by $(\bQ,\tau,\Delta)$, where
\begin{align*}
  \tau(x)&=\sigma_x(\Omega(\Cn_{\Vdash_{\mrm{s}}}(x)))\\
  \Delta(x,y)&= \sigma_{x,y}( \Omega^{-1}(\Cn_{\bQ\restriction\cL}(x\approx y)   )  ),
\end{align*}
as in Remark~\ref{remark:iso:alg}. In particular, $\sigma_x$ is the
substitution sending every variable to $x$, and $\sigma_{x,y}$ is the
substitution sending every variable but $y$ to $x$ and $y$ to
itself. Now, since $\Omega$ is an isomorphism, the least syntactic
theory in $\Vdash$ must be mapped to the least theory in core
semantics, i.e.,
$\Omega(\Cn_{\Vdash}(\emptyset))=\Cn^c_{\bQ}(\emptyset)$. By
assumption, we have that
$\Cn_{\Vdash}(\emptyset)=\Cn_{\Vdash_{\mrm{s}}}(\atomsubst[\Lambda])$ and
$\Cn^c_{\bQ}(\emptyset)=\Cn_{\bQ\restriction\cL}(\atomsubst[\Sigma])$, thus we
obtain that
$\Omega(\Cn_{\Vdash_{\mrm{s}}}(\atomsubst[\Lambda]))=\Cn_{\bQ\restriction\cL}(\atomsubst[\Sigma])$. Let
$\sigma_\delta$ the substitution sending $x$ to $\delta$ and
$\sigma_{\alpha,\beta}$ the substitution sending $x$ to $\alpha$ and
$y$ to $\beta$. Then we obtain that:
\begin{align*}
  \sigma_x(\Omega(\Cn_{\Vdash_{\mrm{s}}}(\bigcup_{\delta\in\atomsubst[\Lambda]}\sigma_{\delta}(x)))) &=  \sigma_x(\Cn_{\bQ\restriction\cL}(\bigcup_{\alpha\approx\beta\in\atomsubst[\Sigma]}\sigma_{\alpha,\beta}(x\approx y)))\\ &\subseteq \Cn_{\bQ\restriction\cL}(\sigma_x(\bigcup_{\alpha\approx\beta\in\atomsubst[\Sigma]}\sigma_{\alpha,\beta}(x\approx y)))
\end{align*}
as $\sigma_x$ is an atomic substitution. It follows that
$\Sigma\models_{\bQ\restriction\cL} \tau(\Lambda)$. Similarly, we have:
\begin{align*}
  \sigma_{x,y}( \Omega^{-1}(\Cn_{\bQ\restriction\cL}(\bigcup_{\alpha\approx\beta\in\atomsubst[\Sigma]}\sigma_{\alpha,\beta}(x\approx y))   )  ) &= \sigma_{x,y}(\Cn_{\Vdash_{\mrm{s}}}(\bigcup_{\delta\in\atomsubst[\Lambda]}\sigma_{\delta}(x))) \\ &\subseteq  \Cn_{\Vdash_{\mrm{s}}}(\sigma_{x,y}(\bigcup_{\delta\in\atomsubst[\Lambda]}\sigma_{\delta}(x))),
\end{align*}
which shows that $\Lambda\Vdash_{\mrm{s}} \Delta(\Sigma)$. By the
algebraizability of $\Vdash_{\mrm{s}}$ and $\Sigma\models_{\bQ\restriction\cL} \tau(\Lambda)$ this entails that
$\tau(\Lambda)\equiv_{\bQ\restriction\cL}\Sigma$ and
$\Delta(\Sigma)\Vdash_{\mrm{s}}\Lambda$,
$\Lambda \Vdash_{\mrm{s}}\Delta(\Sigma)$. Also, from the
assumption that $\Th(\Vdash)=\Lambda\Th(\Vdash_{\mrm{s}})$ and
Lemma~\ref{lambda_filter_lemma} we also obtain that $\Lambda$
witnesses the finite representability of $\Vdash$. We finally conclude
from Corollary~\ref{characterisationalg_loose_5} that $\Vdash$ is
strictly algebraized by $(\bQ,\Sigma, \tau,\Delta)$.
\end{proof}

\section{Aside: Matrix Semantics}\label{sec.5}

In the previous sections we have focused on loosely and strictly
algebraizable weak logics, thus extending the notion of
algebraizability to the setting of logics without uniform
substitution. However, in the setting of abstract algebraic logic,
algebraizability makes only for one of several properties, and it
could be seen as the concept carving out the most well-behaved family
of logical systems. On the converse direction, in this section we work
towards an increased level of generality and we study the matrix
semantics of arbitrary weak logics. In fact, while not every logic is
algebraizable, it can be shown that every logic admits a matrix
semantics (see e.g., \cite[Th. 4.16]{Font.2016}). Furthermore,
Dellunde and Jansana~\cite{dellunde1996some} provided a
characterisation of the class of matrices of a (possibly infinitary)
logic in terms of some model-theoretic results for first-order logic
without equality. We show in this section that similar results can be
proved in the context of weak logics. We prove in \cref{sec:5.1} that
every weak logic is complete with respect to a suitable class of
so-called bimatrices, and we show in \cref{sec:5.2} that Dellunde and
Jansana's results are still applicable in our setting.

\subsection{Completeness of Matrix Semantics}\label{sec:5.1}

We briefly recall the matrix semantics for standard
logics. Intuitively, the idea is to work with first-order structures
with a predicate $\truth(\cA)$ which encodes the ``truth set'' of the
algebra $\cA$.

\begin{definition}
  A \emph{(logical) matrix} of type $\cL$ is a pair
  $(\cA ,\truth(\cA))$ where $\cA$ is an $\cL$-algebra and
  $\truth(\cA) \subseteq \dom(\cA)$.
\end{definition}

\noindent Matrices induce a consequence relation over propositional
formulas analogously as classes of algebras do. However, notice that
here we work directly with propositional formulas (i.e., terms in the
language $\cL$) and not with equations. This corresponds to the fact
that these logics do not necessarily correspond to quasiequational
theories over classes of algebras.

\begin{notation}
  For notational convenience, in this section we denote logics and
  weak logics by $L, L_0,L_1,\dots$. We then write
  $\Gamma\vdash_{L}\phi$ if $(\Gamma,\phi)\in L$. If $L$ is a weak
  logic, then we write $\Gamma\Vdash_{L}\phi$ if $(\Gamma,\phi)\in L$.
\end{notation}

\begin{definition}\label{def:consequence_matrices}
  Let $\bK$ be a class of $\cL$-matrices and let
  $\Gamma\cup\{ \phi \}$ be a set of propositional formulas, then we
  let
  \begin{align*}
    \Gamma \models_{\bK } \phi  \Longleftrightarrow
    & \text{ for all } \cA \in \bK, \; h\in \homset(\Fm,\cA), \\
    & \text{ if } h[\Gamma]\subseteq \truth(\cA) \text{, then } h(\phi)\in \truth(\cA).
  \end{align*}
  Given a logic $L$, we say that $(\cA,\truth(\cA))$ is a \emph{model}
  of $L$ and write $(\cA,\truth(\cA)) \models L$ if, for every
  $\Gamma \cup \{\phi\} \subseteq \cL$, $\Gamma \vdash_L \phi$ entails
  $\Gamma \models_{\{A\}} \phi$. We often abbreviate this latter
  expression by $\models_{\{A\}}\Gamma$.
\end{definition}

We refer the reader to \cite[\S 4]{Font.2016} for a detailed study of
matrix semantics in the context of standard propositional logics. In
particular, \cite[Thm.~4.16]{Font.2016}) states that \emph{every logic
  is complete with respect to a class of matrices}. We show that one
can obtain the same result in the setting of weak logics. First, we
extend the matrix semantics to the setting of weak logics by
introducing a further predicate, i.e., by viewing them as structures
in an algebraic language $\cL$ augmented by two unary predicates, one
for the core set of $A$ and one for the truth set of $A$.

\begin{definition}
  The tuple $(\cA,\truth(\cA),\core(\cA))$ is a \emph{(logical)
    bimatrix} of type $\cL$ if $\cA$ is a $\cL$-algebra,
  $\truth(\cA)\subseteq \dom(\cA)$ and
  $\core(\cA) \subseteq \dom(\cA)$.
\end{definition}

\begin{notation}
  As in the case of expanded algebras, we write $\homset^c(\Fm,\cA)$
  for the set of all assignments $h: \Fm \to \cA$ such that
  $h[\atoms] \subseteq \core(\cA)$.
\end{notation}

\noindent Bimatrices induce a consequence relation analogous to that
of expanded algebras by restricting attention to assignments over core
elements.

\begin{definition}\label{def:consequence_bimatrices}
  Let $\bK$ be a class of $\cL$-bimatrices and let
  $\Gamma\cup\{ \phi \}$ be a set of propositional formulas, then we
  let
  \begin{align*}
    \Gamma \models^c_{\bK } \phi  \Longleftrightarrow
    & \text{ for all } \cA \in \bK, \; h\in \homset^c(\Fm,\cA), \\
    & \text{ if } h[\Gamma]\subseteq \truth(\cA) \text{, then } h(\phi)\in \truth(\cA).
  \end{align*}
  Given a weak logic $L$, we say that $(\cA,\truth(\cA),\core(\cA))$
  is a \emph{model} of $L$, and write
  ${(\cA, \truth(\cA),\core(\cA)) \models L}$ if, for every
  $\Gamma \cup \{\phi\} \subseteq \cL$, $\Gamma \vdash_L \phi$ entails
  $\Gamma \models^c_{\{A\} } \phi$. For set of formulas $\Gamma$ and a
  bimatrix $A$, we write $A\models^c \Gamma$ if
  $\models^c_{\{A\}}\Gamma$.
\end{definition}

Thus, the main intuition behind bimatrices is the same of expanded
algebras: we add a new predicate specifying the core of the matrix, in
order to consider only the assignments sending atomic formulas to
elements of the core. As shown by the following proposition,
bimatrices give rise to several weak logics. As we stressed already
before (cf.~Remark~\ref{remark:finitary}), the finitary requirement in
the following definition is not necessary \emph{per se}, but we need
it as we are focusing on finitary logical systems.

\begin{definition}
  Let $\bK$ be a class of $\cL$-bimatrices, then $\logic(\bK)$ is the
  set of all pairs $(\Gamma,\phi)$ with
  $\Gamma\cup \{\phi\}\subseteq \Fm$ such that
  $\Gamma_0 \corevDash_{\bK} \phi$ for some finite
  $\Gamma_0\subseteq \Gamma$. We say that a weak logic $L$ is complete
  with respect to $\bK$ if $L=\logic(\bK)$.
\end{definition}

\begin{proposition}
  Let $\bK$ be a class of bimatrices, then $\logic(\bK)$ is a weak
  logic.
\end{proposition}
\begin{proof}
  Let $L=\logic(\bK)$, then the Conditions (1), (2) and (3) from the
  definition of consequence relation \ref{def:consequence_relation}
  are immediately valid by
  Definition~\ref{def:consequence_bimatrices}. Additionally, if
  $\Gamma\Vdash_{L}\phi$, then by definition there is some finite
  $\Gamma_0\subseteq \Gamma$ such that $\Gamma_0 \corevDash_\bK \phi$,
  and therefore $\Gamma_0\Vdash_{L}\phi$, which shows that $L$ is also
  finitary.

  \smallskip
  \noindent Finally, we show that $L$ is closed under atomic
  substitutions. Suppose $\Gamma \corevDash_\bK \phi$ but
  $\sigma[\Gamma] \not\corevDash_\bK \sigma(\phi)$ for some atomic
  substitution $\sigma$. Then there are a bimatrix $M\in \bK$ and a
  core assignment $h:\Fm\to M$ such that
  $h[\sigma[\Gamma]]\subseteq \truth(M)$ but
  $h(\sigma(\phi))\notin \truth(M)$. Now, since $h$ is a core
  assignment and $\sigma$ an atomic substitution, it follows that
  $h\circ \sigma\in \homset^c(\Fm,\cA)$, contradicting
  $ \Gamma \corevDash_\cM \phi $.
\end{proof}

By the previous proposition we have that every class of bimatrices
determines a weak logic. Additionally, we can also show that every
weak logic is complete with respect to a class of bimatrices.

\begin{definition}
  For every weak logic $\Vdash$ and every set of propositional
  formulas $\Gamma\subseteq \Fm$, we let $M^\Gamma_\Vdash$ be the
  bimatrix with domain $\dom(M^\Gamma_\Vdash) = \Fm$ and predicates
  $\truth(M^\Gamma_\Vdash) = \Cn_{\Vdash}(\Gamma)$ and
  $\core(M^\Gamma_\Vdash) = \atoms$. We let $\mathbf{M}_\Vdash$ be the
  class of all bimatrices $M^\Gamma_\Vdash$ for
  $\Gamma \subseteq \Fm$.
\end{definition}

\begin{theorem}\label{matrix.completeness}
  Every weak logic $\Vdash$ is complete with respect to the class
  $\mathbf{M}_\Vdash$.
\end{theorem}
\begin{proof}
  We need to show that $\Vdash \; = \logic(\bK)$. Firstly, suppose
  towards contradiction that $\Gamma \Vdash \phi$ but
  $\Gamma \not\models^c_{\mathbf{M}_\Vdash} \phi$. Then there is a
  bimatrix $M^\Delta_\Vdash$ and a core assignment
  $h: \formulas \to M^\Delta_\Vdash$ such that
  $h[\Gamma] \subseteq \truth(M^\Delta_\Vdash)$ and
  $h(\phi) \notin \truth(M^\Delta_\Vdash)$. Consider now the
  substitution $\sigma$ defined letting $\sigma(x)=h(x)$ for all
  $x\in \atoms$. Notice in particular that this is well-defined
  because the domain of $M^\Delta_\Vdash$ is $\Fm$. Since $h$ is a
  core assignment and $\core(M^\Gamma_\Vdash) = \atoms$, it follows
  that $\sigma$ is an atomic substitution, thus we obtain that
  $\sigma(\Gamma)\Vdash \sigma(\phi)$. Now, since
  $h[\Gamma] \subseteq \truth(M^\Delta_\Vdash)$ it follows in
  particular that $\Delta\Vdash \sigma(\Gamma)$ and so by transitivity
  $\Delta\Vdash \sigma(\phi)$. Since
  $\sigma(\phi)=h(\phi)\notin \truth(M^\Delta_\Vdash)$, this
  contradicts the definition of $\truth(M^\Delta_\Vdash)$. It follows
  that $\Gamma \models^c_{\mathbf{M}_\Vdash} \phi$.

  \smallskip
  \noindent Conversely, suppose
  $\Gamma \models^c_{\mathbf{M}_\Vdash} \phi$ and let
  $\id_{\Fm}: \formulas \to M^\Gamma_\Vdash$ be the identity map. Then
  clearly $\id_{\Fm}[\atoms] \subseteq \atoms$ and
  $\id_{\Fm}[\Gamma] \subseteq \truth(M^\Gamma_\Vdash)$. Since
  $\Gamma \models^c_{\mathbf{M}_\Vdash} \phi$ we then obtain
  $\phi=\id_{\Fm}(\phi)\in \truth(M^\Gamma_\Vdash)$, showing
  $\phi \in \Cn_{\Vdash}(\Gamma)$ and thus $\Gamma\Vdash\phi$.
\end{proof}

\subsection{Connections to Model Theory without Equality}\label{sec:5.2}

In the previous section we have established that every class of
bimatrices defines a weak logic and, conversely, that every weak logic
is complete with respect to a class of bimatrices. Here we next
consider what is exactly the class of \emph{all} bimatrices defined by
a weak logic $L$, i.e., the class of all bimatrices $\cM$ for which
$\cM\models^c \Gamma$ entails $\cM\models^c \phi$ whenever
$\Gamma\Vdash_{L} \phi$.  In the standard context this issue was first
considered by Czelakowski \cite{czelakowski_reduced_1980}, who
characterized the class of all matrices which are models of a
logic. Here we follow however the later work of Dellunde and Jansana
in \cite{dellunde1996some}, which provided a novel proof of
Czelakowski's result by employing the fact that (finitary)
propositional logics can be translated into Horn theories
\emph{without equality}. We start by reviewing this translation, which
is essentially a generalisation of what we already considered in
Remark~\ref{remark:first-order-core}.

\begin{proposition}\label{translation} Let $\cL$ be an algebraic
  language, $\Gamma\cup \phi\subseteq \Fm$ and $|\Gamma|<\aleph_0$,
  then we can translate the consequence relations from
  \ref{def:consequence_matrices} and \ref{def:consequence_bimatrices}
  as follows:
  \begin{enumerate}[(a)]
  \item let $\bK$ be a class of matrices then:
    \[ \Gamma \models_{\bK} \phi \; \Longleftrightarrow \; \bK \models
      \forall x_0,\dots,\forall x_n \Big ( \bigwedge_{\gamma\in
        \Gamma} \truth(\gamma(\xbar)) \to \truth(\phi(\xbar) )\Big
      ); \]

  \item let $\bK$ be a class of bimatrices, then:
    \[ \Gamma \models^c_{\bK} \phi \; \Longleftrightarrow \; \bK
      \models \forall x_0,\dots,\forall x_n \Big (
      \bigwedge_{\gamma\in \Gamma} \truth(\gamma(\xbar)) \land
      \bigwedge_{i\leq n} \core(x_i) \to \truth(\phi(\xbar))\Big )
      . \]
  \end{enumerate}
\end{proposition}
\begin{proof}
  This follows immediately from the definition of $\models_{\bK}$ and
  $\corevDash_{\bK}$.
\end{proof}

It follows from the previous proposition that standard logics in the
language $\cL$ can be encoded by Horn theories in
$\cL\cup \{\truth\}$, while weak logics can be encoded in
$\cL\cup \{\truth,\core\}$. This motivates the following definitions.

\begin{notation}
  Let $\Gamma\cup\{\phi\}\subseteq \Fm$, then we write
  $\Phi(\Gamma,\phi)$ for the first-order formula
  \[ \forall x_0,\dots,\forall x_n \Big ( \bigwedge_{\gamma\in \Gamma}
    \truth(\gamma(\xbar)) \land \bigwedge_{i\leq n} \core(x_i) \to
    \truth(\phi(\xbar))\Big ).\]
\end{notation}

\begin{definition}
  Let $\Vdash$ be a weak logic, then we let $\horn(\Vdash) $ be the
  Horn theory obtained by letting $\Phi(\Gamma,\phi)\in \horn(\Vdash)$
  whenever $\Gamma\Vdash \phi$. We write $\Mod(\Vdash)$ for the class
  of structures $\Mod(\horn(\Vdash))$.
\end{definition}

While Czelakowski's original approach in
\cite{czelakowski_reduced_1980} was specifically tailored to logical
matrices, Dellunde and Jansana considered arbitrary model classes
axiomatized by Horn theories without equality, thus making it possible
to apply their results to the setting of bimatrices and expanded
algebras. We recall from \ref{notation:class_operators} that the
operator $\mathbb{H}_{\mrm{s}}$ refers to the closure under
\emph{strict} homomorphic images, and thus
$\mathbb{H}^{-1}_{\mrm{s}}(\bK)$ is the class of all structures $A$
that are preimages of some strict homomorphism of some structure in
$\bK$. We stress that Dellunde and Jansana's results (from
\cite{dellunde1996some}) hold in the setting \emph{without} the
equality symbol. In particular, when in the rest of this section we
consider Horn formulas, we always restrict attention to Horn formulas
not containing the equality symbol. For clarity, we introduce the
following notation.

\begin{notation}
  We write $\cL^{-}$ to refer to a fixed first-order language without
  the equality symbol $\approx$. Also, we write $\cL^{-}$ for the
  collections of all (first-order) formulas in this language (which clearly do not
  contain the equality symbol).
\end{notation}

\begin{fact}[Dellunde, Jansana]\label{dellunde.jansana.1}
  Let $\bK$ be a class of $\cL^{-}$-structures, then the following are
  equivalent:
  \begin{enumerate}[(1)]
  \item $\bK$ is axiomatised by strict universal Horn formulas in
    $\cL^{-}$;
  \item $\bK$ is closed under
    $\mathbb{H}_{\mrm{s}}^{-1},\mathbb{H}_{\mrm{s}},\mathbb{S},\mathbb{P},
    \mathbb{P}_{\mrm{U}} $ and contains a trivial structure;
  \item
    $\bK=\mathbb{H}_{\mrm{s}}^{-1}\mathbb{H}_{\mrm{s}}\mathbb{S}\mathbb{P}
    \mathbb{P}_{\mrm{U}}(\bK_0) $ for some class $\bK_0$ of
    $\cL$-structures containing a trivial structure.
  \end{enumerate}
\end{fact}

As we mentioned in Fact~\ref{birkhoff_maltsev}, the validity of Horn
formulas is always preserved under the operators
$\mathbb{I},\mathbb{S},\mathbb{P}, \mathbb{P}_{\mrm{U}} $. From the
previous theorem it follows that strict Horn formulas in a language
without equality are also preserved under the operators
$ \mathbb{H}_{\mrm{s}}^{-1} $ and $ \mathbb{H}_{\mrm{s}} $. From this
fact and the previous theorem we immediately obtain as a corollary the
following characterization of the class of models $\Mod(\Vdash)$,
where $\Vdash=\logic(\bK) $ for some class of bimatrices $\bK$.

\begin{corollary}
  Let $\bK$ be a class of bimatrices and let $T=\horn(\logic(\bK) )$,
  then:
  \begin{align*}
    \Mod(T)=\mathbb{H}^{-1}_{\mrm{s}}\mathbb{H}_{\mrm{s}}\mathbb{S}\mathbb{P} \mathbb{P}_{\mrm{U}}(\bK').
  \end{align*}
  \noindent where $\bK'$ is $\bK$ together with some trivial
  bimatrices.
\end{corollary}

\noindent Notice that, by \cref{matrix.completeness}, every weak logic
$\Vdash$ is complete with respect to a class of bimatrices
$\mathbf{M}_\Vdash$ and so the previous theorem applies to all weak
logics $\Vdash$ and provides us with a characterization of the class
of all bimatrices defined by $\Vdash$.

Importantly, however, one can see that the class of all matrices
$\Mod(T)$ above does not meet the intuition about the ``right''
semantics of a (weak) logic, and will contain several pathological
examples. For example, in the case of the standard logic $\cpc$, the
class $\Mod(\cpc)$ is \emph{not} the class of Boolean algebras, but
rather the class $\mathbb{H}_{\mrm{s}}^{-1}(\mathbf{BA})$, where here
we identify $\mathbf{BA}$ with the class of expansions $(B,\{1\})$
such that $B$ is a Boolean algebra (for example, the so-called
``benzene ring'' is a matrix model of $\cpc$, but it is not a Boolean
algebra, cf. \cite[Ex.~4.79]{Font.2016}). In order to identify the
non-pathological models of a propositional logic, Dellunde and Jansana
focused in \cite{dellunde1996some} to the so-called \emph{reduced
  structures}. Since they work at the general level of model theory
without equality, it is again straightforward to adapt their results
to our current setting of bimatrices.

\begin{definition}
  Let $\cL$ be any first-order language, $M$ a $\cL$-structure and
  $X\subseteq M$, then we write $\cL(X)$ for the set of all
  $\cL$-formulas with parameters in $X$.
\end{definition}

\begin{definition}
  Let $\cM$ be a first-order structure, $\cD \subseteq \cM$ and
  $\bar{a}\in M^{<\omega}$, we let the \emph{type without equality of
    $\bar{a}$ over $\cD$ in $\cM$} be the following set of
  equality-free formulas:
  \[\text{tp}^-_\cM(\bar{a}/D) = \{ \phi(x)\in \cL^-(\cD) :
    \cM\models\phi(\bar{a}) \}\]
  \noindent Then, the \emph{Leibniz congruence} is the relation
  $\sim_*$ on $ \cA $ defined by letting, for $a,b \in \cM$:
  \[a\sim_* b \Longleftrightarrow \text{tp}^-_\cM(a/\cM)=
    \text{tp}^-_\cM(b/\cM).\] A model $M$ is \emph{reduced} if the
  Leibniz congruence $\sim_*$ over $M$ is the identity.
\end{definition}

\begin{remark}\label{remark:Leibniz_equality}
  As shown in \cite{dellunde1996some}, $\sim_*$ is the largest
  non-trivial congruence relation on $\cM$, meaning that it is the
  largest congruence $\theta$ over its algebraic reduct such that, for
  any relation symbol $R\in\cL$, if $(a_i,b_i)\in \theta$ for all
  $1\leq i \leq n$ then $A\models R(a_1,\dots,a_n)$ if and only if
  $A\models R(b_1,\dots,b_n)$. As we stressed already in
  \ref{remark:def.homomorphism}, the projections induced by
  congruences respecting this conditions are always strict
  homomorphisms. Additionally, since $\sim_*$ is the greatest such
  congruence, any strict homomorphism from $\cM/\sim_*$ onto some
  further structure $N$ must be the identity.
\end{remark}

\begin{notation}
  Let $M$ be an $\cL^{-}$-structure and let $\sim_*$ be the Leibniz
  congruence over $M$, then we write $M^*$ for the quotient structure
  $\cM/\sim_*$. For any class operator $\mathbb{O} $, we let
  $\mathbb{O}^*(C):=\{\cA^* : \cA\in \mathbb{O}(C)\}$, and for any
  class $\bK$ we let $\bK^*=\{ A^* : A\in \bK\}$. Finally, we let
  $\Mod^*(T)=\{A^* : A\in \Mod(T) \}$.
\end{notation}

\begin{remark}
  By Fact~\ref{dellunde.jansana.1}, the validity of (strict) universal
  Horn formulas in a language without equality is closed by
  $\mathbb{H}_{\mrm{s}}$, and thus $\Mod^*(T)\subseteq \Mod(T)$.
\end{remark}

\noindent From Fact~\ref{dellunde.jansana.1}, one immediately obtain
the following theorem \cite[Thm.~18]{dellunde1996some} and
corollary. Notice that our formulation differs from the original one
as Dellunde and Jansana assume that closure operators are already
closed under isomorphic copies.

\begin{fact}[Dellunde, Jansana]\label{dellunde.jansana.2}
  Let $\bK$ be a class of reduced $\cL^{-}$-structures, then the
  following are equivalent:
  \begin{enumerate}[(1)]
  \item $\bK$ is the class of reduced models of a $\cL^{-}$-universal
    Horn theory;
  \item $\bK$ is closed under the operators
    $\mathbb{I}^*, \mathbb{S}^*,\mathbb{P}^*, \mathbb{P}^*_{\mrm{U}}
    $;
  \item
    $\bK=\mathbb{I}^*\mathbb{S}^*\mathbb{P}^*
    \mathbb{P}^*_{\mrm{U}}(\bK_0) $ for some class $\bK_0$ of
    $\cL^{-}$-structures.
  \end{enumerate}
\end{fact}
\begin{corollary}
  Let $\bK$ be a class of reduced bimatrices, $T=\horn(\logic(\bK))$,
  then:
  \[\Mod^*(T)=\mathbb{I}^*\mathbb{S}^*\mathbb{P}^*
    \mathbb{P}^*_{\mrm{U}}(\bK).\]
\end{corollary}

Thus, since every weak logic is complete with respect to a class of
bimatrices, the previous corollary provides us with a characterization
of the class of reduced bimatrices defined by any weak logic
$\Vdash$. We conclude our abstract study of matrix semantics by
showing that, in the case of a loosely algebraizable weak logic
$\Vdash$, its reduced core-generated bimatrices coincide with the
core-generated expanded algebras of its equivalent algebraic
semantics.  We first define these notions and notice that
Proposition~\ref{th:subst-core} easily extends to this setting.

\begin{definition}
  Let $M$ be a bimatrix, we say that $M$ is core-generated if
  $M=\langle \core(M)\rangle$. We then write $\Mod_{\mrm{CG}}(\Vdash)$
  for the subclass of $\Mod(\Vdash)$ consisting only of core-generated
  structures. We then let
  $\Mod^*_{\mrm{CG}}(\Vdash)=(\Mod(\Vdash))^*$.
\end{definition}

\begin{proposition}\label{th:subst-core-matrix} Let $A$ be a
  core-generated bimatrix, then
  $\sigma(\Theta) \corevDash_{\{ A \}} \sigma(\varepsilon \approx
  \delta)$ for all $\sigma \in \subst(\cL)$ holds if and only if
  $ \Theta \models_{\{A\}} \varepsilon \approx \delta$.
\end{proposition}
\begin{proof}
  This follows by reasoning as in
  Lemma~\ref{assignment_vs_core_assignment} and
  Proposition~\ref{th:subst-core}.
\end{proof}

\noindent We can then characterise the reduced, core-generated
bimatrices of a loosely algebraizable weak logics $\Vdash$ by directly
applying the standard version of this result. We recall the following
classical result from \cite[Thm. 4.60]{Font.2016}.

\begin{fact}\label{standard_alg_bridge}
  Let $\vdash$ be an algebraizable standard logic with equivalent
  algebraic semantics $(\bQ,\tau,\Delta)$. Then
  $(A,\truth(A))\in \Mod^*(\vdash)$ if and only if $A\in \bQ$ and
  $\truth(\cA)=\{a\in\cA : \cA\models \tau^A(a)\}$.
\end{fact}

\begin{corollary}
  Let $\Vdash$ be a loosely algebraizable weak logic with equivalent
  algebraic semantics $(\bQ,\tau,\Delta)$. Then
  $(A,\truth(A),\core(A))\in \Mod_{\mrm{CG}}^*(\Vdash)$ if and only if
  $(A,\core(A))\in \coregen{\bQ}$ and
  $\truth(\cA)=\{a\in\cA : \cA\models \tau^A(a)\}$.
\end{corollary}
\begin{proof}
Firstly, notice that if $\Vdash$ is loosely algebraized by $(\bQ,\tau,\Delta)$, then by  \ref{characterisationalg_loose_1} we have that $\Vdash_{\mrm{s}}$ is also algebraized by $(\bQ,\tau,\Delta)$.

\smallskip
\noindent Suppose now that $(A,\truth(A),\core(A))\in \Mod_{\mrm{CG}}^*(\Vdash)$, then it follows by Proposition~\ref{th:subst-core-matrix} that $(A,\truth(A))\in \Mod(\Vdash_{\mrm{s}})$. We claim that $(A,\truth(A))$ is also reduced. If $a,b\in A$ are distinct elements, then since $\Vdash_{\mrm{s}}$ is algebraized by $(\bQ,\tau,\Delta)$ we have that $A\models a\not\approx b$ is equivalent to $A\not\models \Delta^A(a,b)$. But then we have that $\Delta(x,a)\in \text{tp}^-_\cM(a/A)$ but $\Delta(x,a)\notin\text{tp}^-_\cM(b/A)$, showing that $(A,\truth(A))$ is reduced. It thus follows that $(A,\truth(A))\in \Mod^*(\Vdash_{\mrm{s}})$. So, by Fact~\ref{standard_alg_bridge} and the fact that $A$ is core-generated, we then conclude that $(A,\core(A))\in \coregen{\bQ}$ and $\truth(\cA)=\{a\in\cA : \cA\models \tau^A(a)\}$.

\smallskip
\noindent Conversely, suppose that $(A,\core(A))\in \coregen{\bQ}$ and
  $\truth(\cA)=\{a\in\cA : \cA\models \tau^A(a)\}$. By  \ref{standard_alg_bridge} and the fact that  $\Vdash_{\mrm{s}}$ is algebraized by $(\bQ,\tau,\Delta)$, we have that $(A,\truth(A))\in \Mod^*(\Vdash_{\mrm{s}})$. Also, since $(A,\core(A))$ is core-generated, it follows from \ref{th:subst-core-matrix} that $(A,\truth(A),\core(A))\in \Mod_{\mrm{CG}}^*(\Vdash)$. This concludes our proof.  
\end{proof}

\section{Application: Inquisitive and Dependence Logics}\label{sec.7}

We turn in this section to one application of the abstract machinery
that we studied, i.e., the case of (propositional) inquisitive and
dependence logics. These logical systems make for two interesting
examples of logics where uniform substitution fails, and which have
also been studied from the algebraic point of view. In particular, an
algebraic semantics for the classical version of inquisitive logic
$\inqB$ was introduced in \cite{grilletti} (although some preliminary
inquiry into the subject was provided already in
\cite{roelofsen2011algebraic}) and further investigated in
\cite{Quadrellaro.2019B}. Such semantics was later generalised in
\cite{quadrellaro.2021} to the intuitionistic logic $\inqI$, and to
both the classical and intuitionistic version of dependence logic
$\inqB^\otimes$ and $\inqI^\otimes$. Since these logical systems do
not satisfy the rule of uniform substitution, it has so far been an
open question whether they are in any sense unique --- as a matter of
fact, another (similar) semantics to $\inqI$ was suggested
independently in \cite{Puncochar2021-PUNIHA}. The notion of
algebraizability of weak logics that we have introduced in this
article provides us with a framework to make sense of this
question. In this section, we build on these previous works and relate
them to the notion of algebraizability from the present article. Our
main result is that the classical versions of inquisitive and
dependence logic $\inqB$ and $\inqB^\otimes$ are strictly
algebraizable, while their intuitionistic versions $\inqI$ and
$\inqI^\otimes$ are only loosely so.

\subsection{Inquisitive and Dependence Logic}\label{subsec:inquisitive_logic}

We introduce in this section the intuitionistic and classical
propositional variants of inquisitive and dependence logic. We refer
the reader to \cite{Miglioli1989-PIESRO,Ciardelli2011-CIAIL} for the
original presentation of classical propositional inquisitive logic, to
\cite{Yang2016-YANPLO} for classical propositional dependence logic,
and to \cite{ciardelli2020} for their intuitionistic versions. We make
explicit in the following remark the language in which we formulate
these logics. Notice that while our presentation follows essentially
\cite{ciardelli2020}, our notation is as in \cite{quadrellaro.2021}.

\begin{context}\label{language:inquisitive_dependence}
  We let $\langInt$ be the propositional language
  $\langInt=\{ \land, \lor, \rightarrow, \bot, \top \}$, i.e.,
  $\langInt$ is simply the usual language of propositional
  intuitionistic logic. With some abuse of notation, we denote by
  $\langInt$ also the set of all propositional formulas recursively
  built from $\atoms$ in this syntax. Also, we let $\langInqI$ be the
  propositional language
  $\langInqI=\{ \land, \lor, \otimes, \rightarrow, \bot,\top \}$,
  namely the expansion of $\langInt$ by a novel \emph{tensor
    disjunction} operation $\otimes$. We denote by $\langInqI$ also
  the set of all propositional formulas recursively built from
  $\atoms$ in the syntax $\langInqI$.
\end{context}

\begin{notation}
  We define the inquisitive operation
  $?\phi\coloneqq \phi\lor \neg \phi$. We treat negation as a defined
  operation and we define it by $\neg \phi \coloneqq \phi \to
  \bot$. Also, we write $\phi\leftrightarrow \psi$ as an abbreviation
  for $\phi\to \psi \land \psi\to \phi$.
\end{notation}

\begin{definition}
  A formula of $\langInt^\otimes$ is \emph{standard} if it does not
  contain the symbol $\vee$. We write $\langIntsta$ for the signature
  $\langIntsta=\{ \land, \to, \bot,\top \}$ and for the set of
  formulas determined by it. Similarly, we write $\langIntsta^\otimes$
  for the signature
  $\langIntsta^\otimes=\{ \land, \to, \otimes, \bot,\top \}$ and for
  its corresponding set of formulas.
\end{definition}

We recall briefly the standard semantics of inquisitive and dependence
logic, i.e., their \emph{team} (or \emph{state})
\emph{semantics}. Intuitively, while in the standard semantics of
classical propositional logic a formula is evaluated by a truth-table,
i.e., by an assignment of atomic variables into $\{0,1\}$, in the
classical version of team semantics it is evaluated by a \emph{set} of
such assignments. Similarly, while in the standard semantics of
intuitionistic propositional logic a formula is evaluated at a node of
a poset, in the intuitionistic version of team semantics it is
evaluated by a \emph{set} of such nodes. We make this idea precise via
the following definitions from \cite{ciardelli2020}.

\begin{definition}
  A \textit{Kripke frame} is a partial order $\mathfrak{F}=(W, R)$. A
  \textit{Kripke model} is a pair $\mathfrak{M}=(\mathfrak{F},V)$,
  where $\mathfrak{F}$ is an Kripke frame and
  $V:W\rightarrow \wp(\atoms)$ is a valuation of atomic formulas such
  that, if $p\in V(w)$ and $wRv$, then $p\in V(v)$. We say that a
  Kripke frame (resp. model) $\mathfrak{F}=(W, R)$ is \emph{classical}
  if $R$ is the identity relation.
\end{definition}

\begin{definition}
  Let $\mathfrak{M}=(W,R,V)$ be an Kripke model. A \textit{team} is a
  subset $t\subseteq W$ of the set of possible worlds. A team $s$ is
  an \textit{extension} of a team $t$ if $s\subseteq R[t]$.
\end{definition}

\begin{remark}
  Crucially, an element in a Kripke model $(\mathfrak{F},V)$ is
  essentially a classical assignment, and we write $w(p)=1$ if and
  only if $p\in V(w)$ (when the underlying valuation $V$ is clar from the context). This reflects the main underlying intuition
  that teams are essentially sets of assignments.
\end{remark}

\noindent We can next define  the team semantics of the formulas in $ \langInqI $, and thus of the classical and inquisitive variant of inquisitive and dependence logic.

\begin{definition}\label{def:team_semantics}
  Let $\mathfrak{M}=(W,R,V)$ be a Kripke model. The notion of a
  formula $\phi\in\langInqI$ being \textit{true in a team}
  $t\subseteq W$ is defined as follows:
  \begin{equation*}
    \begin{array}{l @{\hspace{1em}\Longleftrightarrow\hspace{1em}} l}
      \mathfrak{M},t\models p & {}\forall w\in t \ ( w(p)=1),  \\
      \mathfrak{M},t\models \bot &  t=\varnothing,\\
      \mathfrak{M},t\models \psi \lor \chi & \mathfrak{M},t\models \psi \text{ or } \mathfrak{M},t\models \chi,\\
      \mathfrak{M},t\models \psi \land \chi & \mathfrak{M},t\models \psi \text{ and } \mathfrak{M},t\models \chi,\\
      \mathfrak{M},t\models \psi \otimes \chi & \exists s,r\subseteq t \text{ such that } s\cup r = t \text{ and } \mathfrak{M},s\models \psi, \mathfrak{M},r\models \chi, \\
      \mathfrak{M},t\models \psi \rightarrow \chi &  \forall s \ ( \text{if }s\subseteq R[t] \text{ and } \mathfrak{M},s\models\psi \text{ then } \mathfrak{M},s\models \chi ).\\
    \end{array}
  \end{equation*}
  If $\Gamma$ is a set of formulas, then we write
  $\mathfrak{M},t\models\Gamma$ if $\mathfrak{M},t\models \phi$ for
  all $\phi\in \Gamma$. Also, we write $\mathfrak{M}\models \phi$ if
  $\mathfrak{M},t\models\phi$ for all $t\subseteq W$ and we write
  $\mathfrak{F}\models \phi$ if $(\mathfrak{F},V)\models \phi$ for all
  valuations $V$ over $W$.
\end{definition}

\begin{definition} We define the following logical systems: \label{definition:inquisitive_logic}
  \begin{enumerate}[(a)]
  \item the system $\Vdash_\inqI$ of intuitionistic inquisitive logic
    is the consequence relation in the language $\langInt$ over the
    class of all Kripke models $\mathfrak{M}$, namely,
    \begin{flalign*}
      \Gamma \Vdash_\inqI \phi \; \Longleftrightarrow   & \;	\mathfrak{M},t \models \Gamma \text{ entails } 	\mathfrak{M},t\models \phi \\
                                                        &\text{ for all Kripke models } \mathfrak{M} \text{ and teams } t\subseteq \mathfrak{M};
    \end{flalign*}
  \item the system $\Vdash_\inqB$ of classical inquisitive logic is
    the consequence relation in the language $\langInt$ over the class
    of all classical Kripke models $\mathfrak{M}$, namely,
    \begin{flalign*}
      \Gamma \Vdash_\inqB \phi \; \Longleftrightarrow   &\; \mathfrak{M},t\models \Gamma \text{ entails } \mathfrak{M},t\models \phi \\
                                                        &\text{ for all classical Kripke models } \mathfrak{M} \text{ and teams } t\subseteq \mathfrak{M};
    \end{flalign*}
  \item the system $\Vdash_{\inqI^\otimes}$ of intuitionistic
    dependence logic is the consequence relation in the language
    $\langInqI$ over the class of all Kripke models $\mathfrak{M}$,
    namely,
    \begin{flalign*}
      \Gamma \Vdash_{\inqI^\otimes} \phi \; \Longleftrightarrow \; &\;  \mathfrak{M},t\models \Gamma \text{ entails } \mathfrak{M},t \models \phi \\
                                                                   & \text{ for all Kripke models } \mathfrak{M}  \text{ and teams } t\subseteq \mathfrak{M};
    \end{flalign*}
  \item the system $\Vdash_{\inqB^\otimes}$ of classical dependence
    logic is the consequence relation in the language $\langInqI$ over
    the class of all classical Kripke models $\mathfrak{M}$, namely,
    \begin{flalign*}
      \Gamma \Vdash_{\inqB^\otimes} \phi \; \Longleftrightarrow \; & \mathfrak{M},t\models \Gamma \text{ entails } \mathfrak{M},t \models \phi \\
                                                                   & \text{for all classical Kripke models } \mathfrak{M} \text{ and teams } t\subseteq \mathfrak{M}.
    \end{flalign*}
  \end{enumerate}
\end{definition}

\begin{remark}\label{compactness_remark}
  Notice that in \cite{ciardelli2020} the previous systems were
  introduced simply as the sets of validities of the consequence
  relations from
  Definition~\ref{definition:inquisitive_logic}. However this does not
  really make a difference, since these logics are all finitary and
  satisfy the deduction theorem. In particular, the fact that each of
  the previous consequence relations is finitary follows from
  \cite[Cor. 4.22]{ciardelli2020}, while the deduction theorem is
  essentially \cite[Prop. 4.3]{ciardelli2020}. It then follows that,
  for any $L\in \{ \inqB,\inqB^\otimes,\inqI, \inqI^\otimes \}$,
  \begin{align*}
    \Gamma \Vdash_{L} \phi \; \Longleftrightarrow \; \Vdash_{L} \bigwedge_{\psi\in \Gamma_0}\psi \to  \phi \text{ for some finite } \Gamma_0\subseteq \Gamma,
  \end{align*}
  showing that each consequence relation $\Vdash_{L} $ is determined
  by its set of validities.
\end{remark}

The next proposition makes it explicit that these logics are proper
examples of weak logics. The failure of uniform substitution in these
logics has been pointed out since their introduction in
\cite{Ciardelli2011-CIAIL,Yang2016-YANPLO,ciardelli2020}. We provide
details of the following proposition (which essentially develops
Example~\ref{example.inquisitive}) for completeness of exposition,
and since the lack of uniform substitution provides the key motivation
of our abstract work.

\begin{proposition}\label{inquisitive_weak}
  Let $L\in \{ \inqB,\inqB^\otimes,\inqI, \inqI^\otimes \}$, then
  $\Vdash_{L}$ is a weak logic, and it is not closed
  under uniform substitution.
\end{proposition}
\begin{proof}
 Let $L\in \{ \inqB,\inqB^\otimes,\inqI, \inqI^\otimes \}$, then the fact that $\Vdash_{L}$ satisfies Conditions (1)--(3) from Definition~\ref{def:consequence_relation} follows from the definition of the semantic consequence relation $\models$ from Definition~\ref{def:team_semantics}. Condition (4) follows from Remark~\ref{compactness_remark} (thus in particular from \cite[Cor. 4.22]{ciardelli2020}). Finally, the fact that $\Vdash_{L}$ is closed under atomic substitution is immediate to verify using the team semantics from above.

 \smallskip
 \noindent It remains to show that $\Vdash_{L}$ is not closed under
 uniform substitution. We prove this rigorously only for $\inqB$ and
 $\inqB^\otimes$ and mention an example that applies also to $\inqI$
 and $ \inqI^\otimes$ below in Remark~\ref{example_non_US}. Suppose
 thus $L\in \{ \inqB,\inqB^\otimes \}$, we formalize the model that we
 sketched in Example~\ref{example.inquisitive}. Let
 $\mathfrak{F}=(\{a,b,c,d\},=)$ be a classical Kripke frame, let
 $V:\atoms\to \wp(\mathfrak{F})$ be a valuation such that
 $V(p)=\{a,b\}$ and let $\mathfrak{M}=(\mathfrak{F},V)$.

 \smallskip
 \noindent Firstly, we use the team semantics from
 \ref{def:team_semantics} to show that
 $\mathfrak{M},\{b,d\}\models \neg \neg ?p$. Suppose on the contrary
 that there is a team $t \subseteq \{b,d\}$ such that
 $\mathfrak{M},t\models \neg ?p$. Then, we have in particular that
 either $\mathfrak{M},\{b\}\not \models ?p$ or
 $\mathfrak{M},\{d\}\not \models ?p$, which are clearly false. We thus
 conclude that $\mathfrak{M},\{b,d\}\models \neg \neg
 ?p$. Additionally, we notice that
 $\mathfrak{M},\{b,d\}\not \models ?p$. In fact,
 $\mathfrak{M},\{b,d\}\not \models p$ since $d(p)=0$, and
 $\mathfrak{M},\{b,d\}\not \models \neg p$ since $b(p)=1$.  This shows
 that $\neg \neg ?p \not \Vdash_{\inqB} ?p $ and
 $\neg \neg ?p \not \Vdash_{\inqB^\otimes} ?p $.

 \smallskip
 \noindent However, suppose that $(\mathfrak{M},t)\models \neg \neg p$
 for some classical Kripke model $\mathfrak{M}=(W,R,V)$ and
 $t\subseteq W$. Then for every $w\in t$ we have
 $(\mathfrak{M},\{w\})\not\models \neg p$, which entails that $w(p)=1$
 for all $w\in t$, and thus that $(\mathfrak{M},t)\models p$. In other
 words, we proved that $\neg \neg p \Vdash_{L} p$ in any of the logics
 $L\in \{ \inqB,\inqB^\otimes\}$. We conclude that $\inqB$ and
 $\inqB^\otimes$ are not closed under the substitution $p\mapsto ?p$.
\end{proof}

\begin{remark}\label{example_non_US}
  We mention a counterexample to uniform substitution that applies to
  all $L\in \{ \inqB,\inqB^\otimes,\inqI, \inqI^\otimes \}$. One can
  verify that:
  \[ \Vdash_{L}(p \rightarrow (q \lor r)) \rightarrow ((p \rightarrow
    q)\lor (p \rightarrow r)), \]
  \noindent but the result of the substitution $p\mapsto q\lor r$ is
  not a validity of these logics:
  \[ \not\Vdash_{L} ((q\lor r) \rightarrow (q\lor r)) \rightarrow
    (((q\lor r) \rightarrow q)\lor ((q\lor r) \rightarrow r)). \]
  \noindent We refer the reader to \cite[4.5]{ciardelli2020} for an
  explanation of this fact and a lengthier discussion of the failure
  of uniform substitution in propositional inquisitive and dependence
  logic.
\end{remark}

\subsection{Strict Algebraizability of $\inqB$ and $\inqB^\otimes$}\label{algebraizability_inquisitive}

We prove in this section the strict algebraizability of the classical
versions of inquisitive and dependence logic.  We recall some
important facts and definitions.

\begin{definition}
  A \textit{Heyting algebra} $\cH$ is a bounded distributive lattice
  with an operation $\to$ such that, for all $a,b,c\in \cH$:
  \[ a\land b \leq c \; \Longleftrightarrow \; a \leq b \to c, \]
  where $\leq$ is the lattice order.  The negation of $a\in H$ is
  $\neg a:= a\to \bot$. An element $a\in\cH$ is \textit{regular} if
  $a = \neg \neg a $. We write $\cH_\neg$ for the subset of regular
  elements of $\cH$, and we say that a Heyting algebra $\cH$ is
  \textit{regular} if $\cH=\langle \cH_\neg \rangle$. We say that a
  variety $\vari$ of Heyting algebras is \emph{regularly generated} if
  it is core-generated with $\core(H)=H_\neg$ for all $H\in \vari$.
\end{definition}

\begin{remark}\label{intermediate_logics_algebraizable}
  We recall that an \textit{intermediate logic} is a standard logic
  $\intlog$ such that $\ipc\subseteq\intlog\subseteq \cpc$, where
  $\ipc$ denotes the intuitionistic propositional calculus, and $\cpc$
  the classical propositional calculus. Every intermediate logic
  $\intlog$ is algebraized by the variety of Heyting algebras
  $\mrm{Var}(\intlog)$ defined by
  \begin{align*}
    \mrm{Var}(\intlog)=\{H\in \mathbf{HA} : H\models \phi \approx 1 \; \text{ for all } \phi\in\langInt \text{ such that} \vdash_\intlog \phi \}
  \end{align*}
  and by the transformers $(\tau, \Delta)$ given by
  $\tau(\phi)=(\phi\approx 1)$ and
  $\Delta(\alpha,\beta)=\alpha\leftrightarrow \beta$
  (cf. \cite[Ex. 3.34]{Font.2016}).
\end{remark}

The following two definitions are less standard. First, we recall the
$\dna$-logics from Example~\ref{example.inquisitive}. Negative
variants were originally considered in \cite{Miglioli1989-PIESRO}.

\begin{definition}
  Let $\intlog$ be an intermediate logic, then its negative variant
  $\intlog^\neg$ is the set of formulas
  $\intlog^\neg=\{ \phi[\neg p_0,\dots,\neg p_n/p_0,\dots,p_n] :
  \phi\in \langInt \text{ and } \vdash_\intlog \phi \}$. A
  $\dna$-logic is the negative variant of some intermediate logic.
\end{definition}

\noindent Additionally, we recall Medvedev's logic of finite problems,
which was originally introduced by Medvedev in \cite{Medvedev}. This
is an intermediate logic defined in terms of validity in a specific
class of Kripke frames. We direct the reader to \cite[\S
2.2]{Zakharyaschev.1997} for details about the Kripke semantics of
intermediate logics.

\begin{definition}
  We recall that the intermediate logic $\ml$, known as Medvedev's
  logic of finite problems, is the logic of all Kripke frames of the
  form $(\wp^+(s),\supseteq)$, where $s$ is a finite set and
  $\wp^+(s)=\wp(s){\setminus}\{\emptyset \}$.
\end{definition}

\noindent The next fact collects together some results by Ciardelli
\cite{Ciardelli:09thesis} on the schematic variant of $\inqB$, and the
characterization of regularly generated varieties from
\cite{Quadrellaro.2019B}. From these it is then immediate to prove the
strict algebraizability of $\inqB$.

\begin{fact}\label{well-known-facts} $\;$
  \begin{enumerate}[(1)]
  \item $\mathtt{ML}^\neg=\inqB$ and $\Schm(\inqB)=\mathtt{ML}$.
  \item $\mrm{Var}(\mathtt{ML})$ is regularly generated.
  \item Let $\tau(\phi)=(\phi\approx 1)$ for all $\phi\in \Fm$, then
    for all $\Gamma \cup \{ \phi \} \subseteq \langInt$,
    $\Gamma \Vdash_{\inqB} \phi $ is equivalent to
    $ \tau(\Gamma) \corevDash_{\mrm{Var}(\mathtt{ML})} \tau(\phi)$,
    where $\core(H)=H_\neg$ for all $H\in \mrm{Var}(\mathtt{ML})$.
  \end{enumerate}
\end{fact}
\begin{proof}
  Considering Clause (1), both $\mathtt{ML}^\neg=\inqB$ and
  $\Schm(\inqB)=\mathtt{ML}$ were proved by Ciardelli in
  \cite{Ciardelli:09thesis}. Clause (2) follows immediately from (1)
  together with Proposition 4.17 from \cite{Quadrellaro.2019B}.
  Clause (3) is essentially the main algebraic completeness result
  for inquisitive logic shown in \cite{grilletti} (and later in
  \cite{Quadrellaro.2019B}). However, in these articles completeness is
  formulated with respect to the intermediate logics $\mathtt{KP}$ and
  $\mathtt{ND}$, thus we explain how to obtain completeness with
  respect to $\mrm{Var}(\mathtt{ML})$ as stated in (3). Firstly, by
  Theorem 3.32 in \cite{Quadrellaro.2019B} we have that
  $\Gamma \models_{\inqB} \phi $ is equivalent to
  $ \Gamma \corevDash_{\mrm{Var}(\mathtt{ND})} \phi$, where
  $\mathtt{ND}$ is a specific intermediate logic and $\core(H)=H_\neg$
  for all $H\in \mrm{Var}(\mathtt{ND})$. Additionally, it follows from
  \cite[Thm. 5.9]{Quadrellaro.2019B} that $\mathtt{ND}$ and
  $\mathtt{ML}$ contain the same regularly generated Heyting
  algebras. It then follows from Proposition~\ref{thm:maltsev} that
  $\Gamma \Vdash_{\inqB} \phi $ is equivalent to
  $\Gamma \corevDash_{\mrm{Var}(\mathtt{ML})} \phi$.
\end{proof}

\begin{theorem}\label{alg.inqb}
  $\inqB$ is strictly algebraizable.
\end{theorem}
\begin{proof}
  Let $\tau(\phi)= \phi\approx 1$, $\Delta(x,y)=x\leftrightarrow y$
  and $\Sigma=\{ x\approx \neg \neg x \}$. We prove that
  $(\mrm{Var}(\mathtt{ML}), \Sigma, \tau,\Delta)$ strictly algebraizes $\inqB$.
  Firstly, by Fact~\ref{well-known-facts}(2) above we have that
  $\mrm{Var}(\mathtt{ML})$ is core-generated with core defined by
  $\Sigma$. Then, by Proposition~\ref{facts.weak.algebraizability} it
  suffices to show that
  $(\mrm{Var}(\mathtt{ML}), x\approx\neg\neg x, \phi\approx 1,
  x\leftrightarrow y)$ satisfies \ref{loose_Alg}(\ref{l.alg1}) and
  \ref{loose_Alg}(\ref{l.alg4}). By Fact~\ref{well-known-facts}(3),
  Condition \ref{loose_Alg}(\ref{l.alg1}) follows immediately.
  Moreover, for all $H\in \mrm{Var}(\mathtt{ML}) $ and $x,y\in\cH$, we
  have that $x=y$ if and only if $x\leq y$ and $y\leq x$. This is
  equivalent to $H\models x\to y \approx 1 $ and
  $H\models y\to x\approx 1$. It follows that
  $x\approx y\equiv^c_{\mrm{Var}(\mathtt{ML})} \{x\to y\approx 1, y\to
  x\approx 1 \}$, showing that \ref{loose_Alg}(\ref{l.alg4}) holds.
  It follows that $\inqB$ is strictly algebraizable.
\end{proof}

To extend the previous theorem to dependence logic, we firstly need to
introduce a suitable notion of dependence algebras. We refer the
reader to \cite[p. 57]{Burris.1981} for the definition of subdirectly
irreducible algebra.

\begin{definition}\label{def:classical_dependence_algebra}
  A $\inqB^\otimes$-algebra $\cA$ is a structure in the signature
  $\langInqI$ satisfying the following conditions:
\begin{enumerate}[(a)]
\item
  $\cA{ \upharpoonright } \{\lor,\land,\to, \bot \} \in
  \mrm{Var}(\mathtt{ML}) $,
\item
  $\cA_\neg{ \upharpoonright } \{\otimes,\land,\to, \bot \} \in \BA$,
\item
  $ \cA \models \forall x \forall y \forall z \; (x \otimes (y \lor z)
  \approx (x\otimes y) \lor (x\otimes z)) $,
\item
  $ \cA \models \forall x \forall y \forall z \forall k \; ((x\to z)
  \to (y \to k) \approx (x\otimes y) \to (z\otimes k)) $;
\end{enumerate}
we then let $ \mathsf{InqBAlg^\otimes} $ be the variety of all
$\inqB^\otimes$-algebras and $ \mathsf{InqBAlg^\otimes_{FRSI}} $ be
its subclass of finite, regular and subdirectly irreducible elements.
\end{definition}

\begin{remark}\label{dependence_algebra_remark}
  Our definition is essentially from \cite[2.2]{quadrellaro.2021},
  with the difference that here we assume that the equations hold in
  the full algebra and not only in the subalgebra generated by the
  core. Since our results deal with core semantics and core-generated
  structure, this does not affect the validity of the results from
  \cite{quadrellaro.2021}.
\end{remark}

\noindent The previous class of algebras was shown to provide a sound
and complete semantics for $\inqB^\otimes$. We recall the following
fact from \cite[2.15, 3.20]{quadrellaro.2021} and use it to show that
$\inqB^\otimes$ is strictly algebraizable. Notice that, since
$\cA{ \upharpoonright } \{\lor,\land,\to, \bot \} \in
\mrm{Var}(\mathtt{ML}) $, for $A\in \mathsf{InqBAlg^\otimes}$, it
follows that the notion of regular elements and the subset $A_\neg$
are well defined in this context.

\begin{fact}\label{complete.dependence}
  Let $\tau(\phi)=(\phi\approx 1)$ for all $\phi\in \Fm$, then for all
  $\Gamma \cup \{ \phi \} \subseteq \langInt^\otimes$,
  $\Gamma \Vdash_{\inqB^\otimes} \phi $ if and only if
  $ \tau(\Gamma) \corevDash_{\mathsf{InqBAlg^\otimes_{FRSI}}}
  \tau(\phi)$, where $\core(A)=A_\neg$ for all
  $A\in\mathsf{InqBAlg^\otimes_{FRSI}}$.
\end{fact}

\begin{theorem}
  $\inqB^\otimes$ is strictly algebraizable.
\end{theorem}
\begin{proof}
  This follows from Fact~\ref{complete.dependence} analogously to
  Theorem~\ref{alg.inqb}, by choosing for equivalent algebraic
  semantics the variety
  $ \mathbb{V}(\mathsf{InqBAlg^\otimes_{FRSI}})$, the transformers
  $\tau(\phi):= \phi\approx 1$ and
  $\Delta(x,y)=x\rightarrow y \land y\rightarrow x$, and the finite
  set of formulas defining the core
  $\Sigma=\{ x\approx \neg \neg x \}$.
\end{proof}

\subsection{Loose Algebraizability of $\inqI$ and $\inqI^\otimes$}
We show in this section that $\inqI$ and $\inqI^\otimes$ are both
loosely algebraizable, but they are not strictly so. The loose
algebraizability of these logics is essentially a corollary of the
algebraic completeness result from \cite{quadrellaro.2021}. We start
by introducing the algebraic semantics for $\inqI$ and $\inqI^\otimes$
described in \cite{quadrellaro.2021}. As in
Definition~\ref{def:classical_dependence_algebra} we slightly modify
the definitions from \cite{quadrellaro.2021} so that the equations
defining $\inqI$-algebras and $\inqI^\otimes$-algebra are valid in the
entire structure and not only in the subalgebra generated by the
core. As stressed in Remark~\ref{dependence_algebra_remark}, this does
not affect the validity of the completeness results from
\cite{quadrellaro.2021}.

\begin{definition}
  A \textit{Brouwerian semilattice} $B$ is a bounded join-semilattice
  with an additional operation $\to$ such that, for all $a,b,c\in B$:
  \[ a\land b \leq c \Longleftrightarrow a \leq b\to c,\] and we write
  \textbf{BS} for the class of all Browerian semilattices.
\end{definition}

\begin{definition}\label{inqui}
  An $\inqI$-algebra $A$ is an expanded algebra in the signature
  $\langInt$ satisfying the following conditions:
  \begin{enumerate}[(a)]
  \item
    $\cA{ \upharpoonright } \{\lor,\land,\to, \bot \} \in \mathbf{HA}
    $,
  \item
    $\core(A) { \upharpoonright } \{\land,\to, \bot \}\in \mathbf{BS}
    $,
  \item
    $A\models \forall x \forall y \forall a \; (\core(a) \to (a
    \rightarrow (x \lor y) \approx (a \rightarrow x)\lor (a
    \rightarrow y)))$;
  \end{enumerate}
  and we then write $ \mathsf{InqIAlg} $ for the class of all
  $\inqI$-algebras.
\end{definition}

\begin{definition}\label{depe}
  An $\inqI^\otimes$-algebra $A$ is is an expanded algebra in the
  signature $\langInqI$ satisfying the following conditions:
  \begin{enumerate}[(a)]
  \item
    $\cA{ \upharpoonright } \{\lor,\land,\to, \bot \} \in \mathbf{HA}
    $,
  \item
    $\core(A) { \upharpoonright } \{\land,\otimes\to, \bot \}\in
    \mathbf{HA} $,
  \item
    $A\models \forall x \forall y \forall a \; (\core(a) \to (a
    \rightarrow (x \lor y) \approx (a \rightarrow x)\lor (a
    \rightarrow y)))$,
  \item
    $ \cA \models \forall x \forall y \forall z \; (x \otimes (y \lor
    z) \approx (x\otimes y) \lor (x\otimes z)) $,
  \item
    $ \cA \models \forall x \forall y \forall z \forall k \; ((x\to z)
    \to (y \to k) \approx (x\otimes y) \to (z\otimes k)) $;
  \end{enumerate}
  and we then write $ \mathsf{InqIAlg^\otimes} $ for the class of all
  $\inqI^\otimes$-algebras.
\end{definition}

We next recall the algebraic completeness result for intuitionistic
inquisitive and dependence logic. This follows immediately from
\cite[Prop.~2.15]{quadrellaro.2021}. From this, it is immediately to
show that both $\inqI$ and $\inqI^\otimes$ are loosely algebraizable.

\begin{fact}\label{complete.dependence_2}
  Let $\tau:\Fm\to\Eq$ be defined by $\tau(x)=(x\approx 1)$, then the
  following completeness results hold:
\begin{enumerate}[(1)]
\item For all $\Gamma \cup \{ \phi \} \subseteq \langInt$,
  $\Gamma \Vdash_{\inqI} \phi $ if and only if
  $ \tau(\Gamma)\corevDash_{\mathsf{InqIAlg}} \tau(\phi).$
\item For all $\Gamma \cup \{ \phi \} \subseteq \langInt^\otimes$,
  $\Gamma \Vdash_{\inqI^\otimes} \phi $ if and only if
  $ \tau(\Gamma) \corevDash_{\mathsf{InqIAlg^\otimes}}\tau(\phi).$
\end{enumerate}
\end{fact}

\begin{theorem}\label{loose.alg.inqi}
  $\inqI$ and $\inqI^\otimes$ are loosely algebraizable.
\end{theorem}
\begin{proof}
  Consider the logic $\inqI$. Let $\mathbf{Q}$ be the quasivariety of
  expanded algebras generated by $\coregen{\mathsf{InqIAlg}}$, i.e.,
  $ \mathbf{Q}= \mathbb{Q}(\coregen{\mathsf{InqIAlg}}) $,
  $\tau(\phi)=(\phi\approx 1)$ and
  $\Delta(\alpha,\beta)=(\alpha\to\beta)\land(\beta\to\alpha)$. Notice
  that, since every
  $H\in\bQ{ \upharpoonright } \{\lor,\land,\to, \bot \}$ is a Heyting
  algebra, it follows that
  $\tau(\Delta(\alpha,\beta))=\alpha\leftrightarrow \beta \approx 1$
  is equivalent to $\alpha\approx \beta$. This establishes
  \ref{loose_Alg}(\ref{l.alg4}). Moreover, we have by the algebraic
  completeness \ref{complete.dependence_2}(1) and by
  Proposition~\ref{preservation.thm} that
  \[ \Gamma\Vdash_{\inqI} \phi \Longleftrightarrow
    \tau(\Gamma)\models^c_{\bQ} \tau(\phi), \] which establishes
  \ref{loose_Alg}(\ref{l.alg1}). It follows from
  Proposition~\ref{facts.weak.algebraizability} that
  $(\bQ,\tau,\Delta)$ loosely algebraizes $\inqI$. The loose
  algebraizability of $\inqI^\otimes$ follows analogously.
\end{proof}

While the loose algebraizability of $\inqI$ and $\inqI^\otimes$
follows straightforwardly from \cite{quadrellaro.2021}, the fact that
they are not strictly algebraizabile is substantially more subtle. In
principle, to show that a logic is not strictly algebraizable one
should test infinitely many finite sets of equations to see whether
any of them defines the core of the algebras in the corresponding
quasivariety of expanded algebras, which is clearly a no-go. However,
the failure of strict algebraizability in $\inqI$ and $\inqI^\otimes$
is witnessed by a simpler patter. Namely, we can find two inquisitive
(dependence) algebras $H,K$ which have the same algebraic reduct, but
different cores. Clearly, if a logic is strictly algebraizable this
cannot happen, as the core is uniquely determined by a finite set of
equations $\Sigma$. We make this intuition precise in the following
proof.

\begin{figure}
	\begin{tikzpicture}[scale=1.00]
	\pgfmathsetmacro{\NODESIZE}{1.5pt}
		\node[highlight] at ( 0,0) {};
		\node[highlight] at ( 0,1) {};
		\node[highlight] at ( 0,2) {};

		\node[dot,label={180:$0 \;$ } ] (0) at ( 0,0) {};
		\node[dot,label={180:$s\;$ }] (s) at (0,1) {};
		\node[dot,label={180:$1\;$ }] (1) at (0,2) {};

		\draw (0) -- (s);
		\draw (s) -- (1);
	-- caption trick
	\node[draw=none] (label) at  (0,-1) {};
	\node [] at (label.south)  {$H_0$};
	\end{tikzpicture} 
	\hspace{5em}
	\begin{tikzpicture}[scale=1.00]
	\pgfmathsetmacro{\NODESIZE}{1.5pt}
		\node[highlight] at ( 0,0) {};
		\node[highlight] at ( 0,2) {};

		\node[dot,label={180:$0\;$}] (0) at ( 0,0) {};
		\node[dot,label={180:$s\;$}] (s) at (0,1) {};
		\node[dot,label={180:$1\;$}] (1) at (0,2) {};

	\draw (0) -- (s);
		\draw (s) -- (1);
	-- caption trick
	\node[draw=none] (label) at  (0,-1) {};
	\node [] at (label.south)  {$H_1$};
	\end{tikzpicture}
	\hspace{5em}
	\begin{tikzpicture}[scale=1.00]
	\pgfmathsetmacro{\NODESIZE}{1.5pt}
		\node[highlight] at ( 0,0) {};
		\node[highlight] at ( -1,1) {};
		\node[highlight] at ( 0,3) {};
    \node[highlight] at ( 1,1) {};
		\node[dot,label={180:$0\;$}] (0) at ( 0,0) {};
		\node[dot,label={180:$a\;$}] (a) at (-1,1) {};
		\node[dot,label={0:$\;b$}] (b) at ( 1,1) {};
		\node[dot,label={180:$s\;$}] (s) at (0,2) {};
		\node[dot,label={180:$1\;$}] (1) at (0,3) {};

		\draw (0) -- (a);
		\draw (0) -- (b);
		\draw (a) -- (s);
		\draw (b) -- (s);
		\draw (s) -- (1);
	-- caption trick
	\node[draw=none] (label) at  (0,-1) {};
	\node [] at (label.south)  {$H_2$};
	\end{tikzpicture} \caption{The inquisitive (dependence) algebras from \cref{inqI.nonalgebraizable}.}
	\label{examples_inquisitive_algebra}
\end{figure}

\begin{theorem}\label{inqI.nonalgebraizable}
  $\inqI$ and $\inqI^\otimes$ are not strictly algebraizable.
\end{theorem}
\begin{proof}
  We prove this for $\inqI$, as the proof easily adapts to
  $\inqI^\otimes$. Firstly, notice that if $(\bQ,\Sigma,\tau,\Delta)$
  strictly algebraizes $\inqI$, then clearly $(\bQ,\tau,\Delta)$ is
  also a witness of its loose algebraizability. It then follows from
  \cref{loose.alg.inqi} and \cref{uniquenessalg_1} that
  $\bQ= \mathbb{Q}(\coregen{\mathsf{InqIAlg}})$. Now, by directly
  checking the definition of $\inqI$-algebras from
  Definition~\ref{inqui} (or, alternatively, by applying the
  categorical duality between finite Kripke frames and finite,
  well-connected, core-generated, $\inqI$-algebras from
  \cite{quadrellaro.2021}), one can verify that the expanded algebras
  $H_0$ and $H_2$ from \cref{examples_inquisitive_algebra} (with
  circles indicating which elements are in the core) belong to
  $\coregen{\mathsf{InqIAlg}}$. Furthermore, since $H_1\leq H_2$, it
  also follows that $H_1\in
  \mathbb{Q}(\coregen{\mathsf{InqIAlg}})$. Crucially, $H_0$ and $H_1$
  have the same algebraic reduct but different subsets of core
  elements. In particular, assuming that $(\bQ,\Sigma,\tau,\Delta)$
  strictly algebraizes $\inqI$, we obtain that
  $\core(H_0)=\Sigma(H_0)=\Sigma(H_1)=\core(H_1)$, contradicting the
  fact that $\core(H_0)\neq \core(H_1)$. It follows that $\inqI$ is
  not strictly algebraizable. The same argument also works for
  $\inqI^\otimes$, as it suffices to notice that the core subsets of
  $H_0,H_1,H_2$ can all be expanded to form Heyting algebras with an
  additional disjunction $\otimes$, which additionally satisfies the
  axioms from \ref{depe} (this is done in full generality in \cite[\S
  4.3]{quadrellaro.2021}).
\end{proof}

\begin{corollary}
  $\inqI$ and $\inqI^\otimes$ are not finitely representable.
\end{corollary}
\begin{proof}
  Immediate from \cref{loose.alg.inqi}, \cref{inqI.nonalgebraizable},
  and \cref{characterisationalg_loose}.
\end{proof}

\begin{question}
  We notice that in the proof of \cref{inqI.nonalgebraizable} we were
  in a sense lucky, i.e., we proved the impossibility of strict
  algebraizability by showcasing two expanded algebras with the same
  algebraic reduct and different cores. Must such situation always
  happen whenever a weak logic is loosely algebraizable, but not
  strictly so? More precisely, consider the following properties of a
  weak logic $\Vdash$:
\begin{enumerate}[(i)]
    \item $\Vdash$ is loosely algebraized by a tuple $(\bQ,\tau,\Delta)$;
    \item $\Vdash$ is not strictly algebraizable;
    \item if $A\restriction \cL = B\restriction \cL$ then $\core(A)=\core(B)$.
\end{enumerate}
\emph{Is it possible to find a weak logic $\Vdash$ which satisfies the
  properties (i)-(iii) from above?}
\end{question}

\begin{question}
  We notice that the argument for the loose algebraizability of
  $\inqI$ and $\inqI^\otimes$ easily generalises to the entire class
  of intermediate inquisitive (and dependence) logics from
  \cite{quadrellaro.2021}. Now, as mentioned in
  Fact~\ref{well-known-facts}, the schematic fragment of $\inqB$ is
  $\mathtt{ML}$. Moreover, as pointed out in
  \cite{ferguson2023structural}, it is a corollary of the main result
  from \cite{grilletti2022medvedev} that the schematic fragment of
  $\inqI$ is $\mathtt{ML}$ as well. In other words, all the
  intermediate inquisitive logics $\mathtt{L}$ have the same schematic
  fragment $\mathtt{ML}$.  It thus follows from our
  Theorem~\ref{characterisationalg_loose} that the strictly
  algebraizable intermediate inquisitive logics are exactly those
  obtained by closure under modus ponens of sets of the form
  $\mathtt{ML} \cup \atomsubst[\Lambda]$, for some finite
  set of formulas $\Lambda\subseteq \langInt$ in one variable. Is it possible to
  refine this characterisation? For example, \emph{can one provide a
    semantical characterisation of the intermediate inquisitive and
    dependence logics which are strictly algebraizable?}  We leave
  this and the previous question as pointers for future
  investigations.
\end{question}

\section*{Acknowledgement}

The authors would like to thank Tommaso Moraschini for useful discussions and pointers to the literature, Fan Yang and Fredrik Nordvall Forsberg for their valuable comments and remarks to an early draft of this article, and one anonymous reviewer for their careful reading of the article and their many very useful suggestions. 

\section*{Funding}

The second author was supported by grant 336283 of the Academy of Finland and Research Funds of the University of Helsinki, and by an InDam postdoc scholarship.

\bibliographystyle{abbrv}
\bibliography{references}

@article{casanovas1996elementary,
	title={On elementary equivalence for equality-free logic},
	author={Casanovas, Enrique and Dellunde, Pilar and Jansana, Ramon},
	journal={Notre Dame Journal of Formal Logic},
	volume={37},
	number={3},
	pages={506--522},
	year={1996},
	publisher={Duke University Press}
}

@article{nicolau2024polyatomic,
	title={Polyatomic logics and generalized {B}lok--{E}sakia theory},
	author={Almeida, Rodrigo Nicolau},
	journal={Journal of Logic and Computation},
	volume={34},
	number={5},
	pages={887--935},
	year={2024},
	publisher={Oxford University Press}
}

@article{Medvedev,
	Author = {J.T. Medvedev},
	Journal = {Soviet Mathematics Doklady},
	Number = {1},
	Pages = {227--230},
	Title = {Finite Problems},
	Volume = {3},
	Year = {1962}}

@book{Rasiowa1974-RASAAA-2,
	publisher = {Amsterdam, Netherlands: Warszawa, Pwn - Polish Scientific Publishers},
	year = {1974},
	author = {Helena Rasiowa},
	title = {An Algebraic Approach to Non-Classical Logics}
}

@article{quadrellaro.2021,
title = {On intermediate inquisitive and dependence logics: An algebraic study},
journal = {Annals of Pure and Applied Logic},
volume = {173},
number = {10},
year = {2022},
issn = {0168-0072},
doi = {https://doi.org/10.1016/j.apal.2022.103143},
author = {Davide Emilio Quadrellaro}
}

@article{Tarski1986-TARWAL,
	title = {What Are Logical Notions?},
	pages = {143--154},
	author = {Alfred Tarski},
	doi = {10.1080/01445348608837096},
	number = {2},
	year = {1986},
	journal = {History and Philosophy of Logic},
	publisher = {Taylor \& Francis},
	volume = {7}
}

@article{czelakowski_reduced_1980,
	title = {Reduced products of logical matrices},
	volume = {39},
	issn = {1572-8730},
	url = {https://doi.org/10.1007/BF00373095},
	doi = {10.1007/BF00373095},
	number = {1},
	journal = {Studia Logica},
	author = {Czelakowski, Janusz},
	month = mar,
	year = {1980},
	pages = {19--43},
}

@book{Blok1989-BLOAL,
Author = {Blok, W. J. and Pigozzi, D.},
Title = {Algebraizable logics},
FSeries = {Memoirs of the American Mathematical Society},
Series = {Memoirs of the American Mathematical Society},
ISSN = {0065-9266},
Volume = {396},
ISBN = {978-0-8218-2459-7; 978-1-4704-0816-9},
Year = {1989},
Publisher = {American Mathematical Society (AMS)},
Language = {English},
DOI = {10.1090/memo/0396},
zbMATH = {4085644},
Zbl = {0664.03042}
}

@article{dellunde1996some,
	title={Some characterization theorems for infinitary universal Horn logic without equality},
	author={Dellunde, Pilar and Jansana, Ramon},
	journal={The Journal of Symbolic Logic},
	volume={61},
	number={4},
	pages={1242--1260},
	year={1996},
	publisher={Cambridge University Press}
}

@article{ciardelli2020,
	author = "Ciardelli, Ivano and Iemhoff, Rosalie and Yang, Fan",
	doi = "10.1215/00294527-2019-0033",
	journal = "Notre Dame Journal of Formal Logic",
	month = "01",
	number = "1",
	pages = "75--115",
	publisher = "Duke University Press",
	title = "Questions and Dependency in Intuitionistic Logic",
	url = "https://doi.org/10.1215/00294527-2019-0033",
	volume = "61",
	year = "2020"
}

@article{Holliday2013-HOLIDA-5,
	doi = {10.1007/s11229-013-0278-0},
	author = {Wesley H. Holliday and Tomohiro Hoshi and Thomas F. Icard},
	pages = {31--55},
	journal = {Synthese},
	year = {2013},
	volume = {190},
	number = {1},
	title = {Information Dynamics and Uniform Substitution}
}

@InProceedings{10.1007/978-3-642-24130-7_6,
	author="Holliday, Wesley H.
	and Hoshi, Tomohiro
	and Icard, Thomas F.",
	editor="van Ditmarsch, Hans
	and Lang, J{\'e}r{\^o}me
	and Ju, Shier",
	title="Schematic Validity in Dynamic Epistemic Logic: Decidability",
	booktitle="Logic, Rationality, and Interaction",
	year="2011",
	publisher="Springer Berlin Heidelberg",
	address="Berlin, Heidelberg",
	pages="87--96",
	isbn="978-3-642-24130-7"
}

@article{buss1990modal,
	title={The modal logic of pure provability.},
	author={Buss, Samuel R},
	journal={Notre Dame Journal of Formal Logic},
	volume={31},
	number={2},
	pages={225--231},
	year={1990},
	publisher={Duke University Press}
}

@article{Miglioli1989-PIESRO,
	author = {Pierangelo Miglioli and Ugo Moscato and Mario Ornaghi and Silvia Quazza and Gabriele Usberti},
	number = {4},
	publisher = {Duke University Press},
	doi = {10.1305/ndjfl/1093635238},
	pages = {543--562},
	title = {Some Results on Intermediate Constructive Logics},
	journal = {Notre Dame Journal of Formal Logic},
	year = {1989},
	volume = {30}
}

@book{Zakharyaschev.1997,
	year = {1997},
	author = {Alexander Chagrov and Michael Zakharyaschev},
	title = {Modal Logic},
	publisher = {Clarendon Press},
	address = {Oxford}
}

@book{Font.2016,
	year = {2016},
	title = {Abstract Algebraic Logic},
	subtitle = {An Introductory Textbook},
	publisher = {College Publication},
	address =   {London},
	author = {Font, Josep Maria}
}

@inproceedings{roelofsen2011algebraic,
	title={Algebraic foundations for inquisitive semantics},
	author={Roelofsen, Floris},
	booktitle={International Workshop on Logic, Rationality and Interaction},
	pages={233--243},
	year={2011},
	organization={Springer}
}

@InProceedings{grillq,
author="Grilletti, Gianluca
and Quadrellaro, Davide Emilio",
editor="{\"O}zg{\"u}n, Ayb{\"u}ke
and Zinova, Yulia",
title="Lattices of Intermediate Theories via {R}uitenburg's Theorem",
booktitle="Language, Logic, and Computation",
year="2022",
publisher="Springer International Publishing",
address="Cham",
pages="297--322",
isbn="978-3-030-98479-3"
}

@article{Quadrellaro.2019B,
	title={An {A}lgebraic {A}pproach to {I}nquisitive and {DNA}-Logics},
	author={Nick Bezhanishvili and Gianluca Grilletti and Davide Emilio Quadrellaro},
  journal={The Review of Symbolic Logic},
volume={15},
number={4},
pages={950--990},
year={2022},
publisher={Cambridge University Press}
}

@article{ferguson2023structural,
  title={Inquisitive split and structural completeness},
  author={Ferguson, Thomas and Pun{\v{c}}och{\'a}{\v{r}}, V{\'\i}t},
  journal={Mathematical Structures in Computer Science},
  volume={35},
  year={2025},
  publisher={Cambridge University Press}
  }

@article{Puncochar2016-PUNAGO,
	volume = {45},
	journal = {Journal of Philosophical Logic},
	publisher = {Springer Verlag},
	title = {A Generalization of Inquisitive Semantics},
	author = {V\'{i}t Pun\v{c}och\'{a}{\v{r}}},
	pages = {399--428},
	year = {2016},
	number = {4}
}

@article{Puncochar2021-PUNIHA,
	pages = {995--1017},
	year = {2021},
	number = {5},
	author = {V{\'{\i}}t Pun{\v{c}}och{\'{a}}{\v{r}}},
	journal = {Studia Logica},
	title = {Inquisitive {H}eyting Algebras},
	volume = {109},
	publisher = {Springer Verlag}
}

@article{Yang2016-YANPLO,
	doi = {10.1016/j.apal.2016.03.003},
	author = {Fan Yang and Jouko V\"a\"an\"anen},
	volume = {167},
	journal = {Annals of Pure and Applied Logic},
	pages = {557--589},
	title = {Propositional Logics of Dependence},
	number = {7},
	year = {2016}
}

@article{Ciardelli2011-CIAIL,
	volume = {40},
	year = {2011},
	title = {Inquisitive Logic},
	journal = {Journal of Philosophical Logic},
	pages = {55--94},
	doi = {10.1007/s10992-010-9142-6},
	author = {Ivano Ciardelli and Floris Roelofsen},
	publisher = {Springer},
	number = {1}
}

@book{Burris.1981,
	year = {1981},
	author = {Stanley N. Burris and H.P. Sankappanavar},
	title = {A Course in Universal Algebra},
	publisher = {Springer}
}

@inproceedings{grilletti2022medvedev,
 Author = {Grilletti, Gianluca},
 Title = {Medvedev logic is the logic of finite distributive lattices without top element},
 BookTitle = {Advances in modal logic. Vol. 14. Proceedings of the 14th conference (AiML 2022), Rennes, France, August 22--25, 2022},
 ISBN = {978-1-84890-413-2},
 Pages = {451--466},
 Year = {2022},
 Publisher = {London: College Publications},
 Language = {English},
 zbMATH = {7668110},
 Zbl = {1531.03067}
}

@incollection{grilletti,
	title = {Algebraic and Topological Semantics for Inquisitive Logic Via Choice-Free Duality},
year = {2019},
pages = {35--52},
booktitle = {Logic, Language, Information, and Computation. WoLLIC 2019. Lecture Notes in Computer Science, Vol. 11541},
author = {Nick Bezhanishvili and Gianluca Grilletti and Wesley H. Holliday},
publisher = {Springer}
}

@article{moraschini2018logical,
	title={A logical and algebraic characterization of adjunctions between generalized quasi-varieties},
	author={Moraschini, Tommaso},
	journal={The Journal of Symbolic Logic},
	volume={83},
	number={3},
	pages={899--919},
	year={2018},
	publisher={Cambridge University Press}
}

@unpublished{Ciardelli:09thesis,
	author               = {Ivano Ciardelli},
	title                = {Inquisitive semantics and intermediate logics},
	year                 = {2009},
	note                 = {MSc Thesis, University of Amsterdam}
}

@book{chang1990model,
   title={Model theory},
   author={Chang, Chen Chung and Keisler, Howard  Jerome},
   year={1990},
   publisher={Elsevier}
}

\end{document}